\documentclass[11pt,a4paper]{article}

\usepackage{amsmath}
\usepackage{amssymb}
\usepackage{amsthm}
\usepackage{geometry}
\usepackage{float}
\usepackage{hyperref}
\usepackage{mathtools}
\usepackage{tikz}
\usetikzlibrary{arrows.meta,decorations.markings,positioning}

\tikzset{
 rhp diagram/.style={line cap=round,every node/.style={font=\small,inner sep=2pt}},
 rhp contour/.style={line width=0.75pt,postaction={decorate},
  decoration={markings,mark=at position 0.52 with
   {\arrow{Latex[length=2.2mm,width=1.6mm]}}}},
 rhp circle/.style={line width=0.75pt,postaction={decorate},
  decoration={markings,
   mark=at position 0.25 with {\arrow{Latex[length=2.2mm,width=1.6mm]}},
   mark=at position 0.75 with {\arrow{Latex[length=2.2mm,width=1.6mm]}}}},
 rhp guide/.style={black!45,line width=0.45pt,densely dashed}
}

\hypersetup{hidelinks}

\newtheorem{Theorem}{Theorem}[section]
\newtheorem{Proposition}[Theorem]{Proposition}
\newtheorem{Lemma}[Theorem]{Lemma}
\newtheorem{Corollary}[Theorem]{Corollary}
\newtheorem{RHProblem}[Theorem]{RH problem}
\newtheorem{rhp}[Theorem]{RH problem}
\theoremstyle{definition}
\newtheorem{Assumption}{Assumption}
\theoremstyle{remark}
\newtheorem{Remark}[Theorem]{Remark}
\newtheorem*{Remark*}{Remark}

\newcommand{\R}{\mathbb{R}}
\newcommand{\Prob}{\mathbb{P}}
\newcommand{\dd}{\mathrm{d}}
\newcommand{\ii}{\mathrm{i}}
\newcommand{\Ai}{\mathrm{Ai}}
\newcommand{\Bes}{\mathrm{Bes}}
\newcommand{\Boh}{\mathcal{O}}

\title{Tail Asymptotics for Gaps in Unitary Ensembles and the Airy Point Process}
\author{Junwen Liu\thanks{Department of Mathematics, The University of Hong Kong.
Email: \texttt{jwliu26@hku.hk}}}
\date{}

\begin{document}

\maketitle

\begin{abstract}
We derive tail asymptotics for the gap between the $k$th and $(k+1)$st
largest particles, with $k$ fixed, in unitary ensembles and in the Airy point
process. For the finite-$N$ ensembles, we consider the probabilities of a
rescaled gap larger than $a_N$ and of a rescaled gap smaller than $b_N$,
in the joint limits $N\to\infty$, $a_N\to\infty$ with
$a_N=o(N^{2/3})$, and $N\to\infty$, $b_N\downarrow0$, respectively.
For the Airy point process, the right and left tails are analyzed as the
gap size tends to infinity and zero, respectively.
More precisely, in both models, we obtain logarithmic right-tail
asymptotics up to the constant term and cubic left-tail asymptotics with
coefficient $k/6$. As an application, we also derive the limiting conditional laws of the leading particles given the existence of an abnormally large or small gap.
In particular, for $k=1$, large-gap conditioning leads to a transition from
Tracy--Widom to Gaussian fluctuations of the largest particle, similar to
the BBP transition.
Further, among any fixed number of leading gaps, we also determine the limiting probabilities of which gap attains a large maximum or a small minimum.

\end{abstract}

\section{Introduction}
\label{sec:introduction}

We consider the unitary ensemble of $N\times N$ Hermitian matrices $M$
with density proportional to
$e^{-N\operatorname{Tr}V(M)}$, where $V$ is a one-cut regular polynomial
potential as specified in Assumption~\ref{ass:potential} below. The ordered eigenvalues
$\lambda_1>\cdots>\lambda_N$ have joint density
\begin{align}
 \frac1{\mathcal Z_N}
 \prod_{1\leq i<j\leq N}(\lambda_i-\lambda_j)^2
 \exp\left(-N\sum_{i=1}^NV(\lambda_i)\right),
\label{eq:unitary-density}
\end{align}
where $\mathcal Z_N$ is the normalization constant.
For $V(x)=x^2/2$, the ensemble reduces to the Gaussian unitary ensemble
(GUE). In this case, the centered and rescaled largest eigenvalue
converges to the Tracy--Widom distribution $F_2$, which admits an
Airy-kernel Fredholm determinant representation and a representation in
terms of a Painlev\'e II transcendent~\cite{TW94}.
For regular analytic potentials with one-cut support and sufficient growth
at infinity, the same limiting law follows from Riemann--Hilbert (RH)
asymptotics for the associated orthogonal polynomials~\cite{DKMVZ99};
an alternative proof was given by Pastur and Shcherbina~\cite{PS03}.

Small deviation estimates quantify the probabilities of deviations of
the largest eigenvalue from the spectral edge on scales between the
Tracy--Widom fluctuation scale and the macroscopic scale;
see Ledoux~\cite{Ledoux07} for a review.
For every $\beta>0$, Ram\'{\i}rez, Rider and Vir\'ag~\cite{RRV11}
established the leading logarithmic asymptotics of both tails of the
limiting Tracy--Widom $\beta$ distribution, and Dumaz and
Vir\'ag~\cite{DV13} subsequently refined the right-tail estimate.
For Hermite and Laguerre $\beta$-ensembles, Ledoux and
Rider~\cite{LR10} established nonasymptotic bounds for both tails,
while Baslingker et al.~\cite{BBBK25} obtained upper and lower moderate
deviation bounds with optimal constants in the exponents.
Paquette and Zeitouni~\cite{PZ17} proved two-sided right-tail bounds
for the GUE, and Eichelsbacher, Kriecherbauer and Sch\"uler~\cite{EKS16}
obtained precise right-tail asymptotics for the largest eigenvalue in
a class of unitary invariant ensembles.
For real symmetric and complex Hermitian Wigner matrices,
Erd\H{o}s and Xu~\cite{EX23} established upper bounds for both tails,
and Bao et al.~\cite{BCEHK25} obtained upper and lower bounds for
both tails with optimal constants in the exponents.
For fixed $k\ge2$, Baslingker~\cite{Baslingker26} obtained upper and
lower bounds for the right tail of the $k$th largest eigenvalue,
with the optimal exponential constant, in certain Hermite and
Laguerre $\beta$-ensembles.

A related soft-edge statistic is the gap between the two largest
eigenvalues. For the GUE, its reciprocal determines the limiting
distribution of the rescaled halting time of the Toda algorithm for
computing the largest eigenvalue~\cite{DeiftTrogdon}.
Perret and Schehr~\cite{PerretSchehr} obtained a representation for the
limiting density of this gap after soft-edge rescaling, using a
Painlev\'e XXXIV Lax pair. They then derived formal expansions of this
limiting density as the rescaled gap tends to zero or infinity, assuming
a power series for the Lax-pair solution in the former case and an unproved
boundary condition in the latter. At the soft edge of unitary ensembles, Witte, Bornemann and
Forrester~\cite{WBF} expressed the joint density of the two largest
eigenvalues, and hence their gap density, in terms of a Painlev\'e II
transcendent and its associated linear system. More generally, Claeys and
Doeraene~\cite{ClaeysDoeraene} expressed joint distributions of Airy
particles in terms of solutions to coupled Painlev\'e II equations.

For fixed $k\ge1$, a natural question is how likely the gap between the
$k$th and $(k+1)$st largest particles is to be abnormally large or small.
Such gap events can also be viewed dynamically through the Airy line
ensemble~\cite{CorwinHammond14}, whose fixed-time marginal is the Airy point
process and which is stationary and ergodic~\cite{CorwinSun14}.
Heuristically, over long observation times one expects to encounter
abnormally large and small gaps between any fixed pair of adjacent curves.
This motivates studying how the leading particles are distributed under
either condition.
We address these questions for the finite-$N$ ensemble
\eqref{eq:unitary-density} and for the Airy point
process~\cite{Forrester93,TW93}.
For the finite-$N$ ensemble, we study the probability that the
soft-edge rescaled gap exceeds $a_N$ as $N\to\infty$ and $a_N\to\infty$,
with $a_N=o(N^{2/3})$.
We also study the probability that it is smaller than $b_N$ along every
joint sequence $N\to\infty$, $b_N\downarrow0$.
For the Airy point process, we study the probabilities that the gap exceeds
$a$ as $a\to\infty$ and that it is smaller than $b$ as $b\downarrow0$.
We obtain logarithmic right-tail asymptotics up to the
constant term and cubic left-tail asymptotics with coefficient $k/6$
in both models.
When $k=1$, integrating the formal density expansions of Perret and Schehr
gives the same Airy tail formulas as ours, up to the orders stated below.
We establish these asymptotics rigorously.
Our approach is based on Janossy densities and the Riemann--Hilbert
representations of the associated resolvent kernels.

We also determine the limiting conditional laws of the leading particles
in both models.
For fixed $k\ge2$, the centered and rescaled $k$th particle position
converges to a standard Gaussian variable under large-gap conditioning.
In the joint limit, the rescaled distances to the particles above it are
independent of this variable and have the law of an ordered
$(k-1)$-particle Laguerre unitary ensemble with parameter $\alpha=2$.
For every fixed $k\ge1$, small-gap conditioning gives a joint limit for
the first $k$ edge-scaled positions and the gap
divided by $b$ in the Airy model or by $b_N$ in the finite-$N$ ensemble.
Both the large-gap and small-gap limits hold in total variation.
For any fixed $k\ge1$ and $a>0$, the Janossy formula and the Airy
resolvent identities~\cite{TW94} give a Painlev\'e II representation of the
joint law of the leading $k$ Airy particles conditioned on the $k$th gap
exceeding $a$.
For $k=1$, increasing $a$ produces a Tracy--Widom to Gaussian crossover
for the largest particle, similar to the rank-one BBP
transition~\cite{BBP05}.
Section~\ref{subsec:conditional-interpretation} explains these conditional
laws and the cubic left tail in terms of the conditioned particle
configurations.

These tail asymptotics and conditional laws also describe the maximum and
minimum of the first $k$ gaps for any fixed $k\ge1$.
A large maximum is asymptotically attained by the first gap, whereas the
limiting probability that the $j$th gap is the smallest is
$2j/[k(k+1)]$, $1\le j\le k$, conditional on the minimum being small.

\subsection{Main results}
\label{sec:main-results}
\label{sec:static-results}

The macroscopic eigenvalue distribution is described by the equilibrium
measure, the unique minimizer, over
probability measures $\mu$ on $\R$, of the logarithmic energy
\cite{DKMVZ99,RZ}
\begin{equation*}
 \int_{\R}V(x)\,\dd\mu(x)
 -\iint_{\R^2}\log|x-y|\,\dd\mu(x)\,\dd\mu(y).
\end{equation*}
Throughout this paper, we assume that the potential $V$ satisfies the
following one-cut regularity conditions; see \cite{DKMVZ99} and
\cite[Section~2.1]{RZ} for further details.

\begin{Assumption}\phantomsection\label{ass:potential}
\leavevmode
\begin{enumerate}
 \item[(i)] $V$ is a fixed real polynomial of even degree at least two, with
 positive leading coefficient.
 \item[(ii)] Its equilibrium measure has support $[A,B]$, where
  $A<B$, and density
  \begin{equation}
   \rho(x):=\sqrt{(B-x)(x-A)}\,h(x),\qquad A<x<B,
  \label{eq:equilibrium-density}
  \end{equation}
  with $h$ real analytic
 and strictly positive in a neighborhood of $[A,B]$.
 \item[(iii)] There is an Euler--Lagrange constant $\ell$ such that
 \begin{align}
  2\int_A^B\log|x-t|\rho(t)\,\dd t-V(x)&=\ell,
    &&x\in[A,B],\label{eq:equilibrium-EL-support}\\
  2\int_A^B\log|x-t|\rho(t)\,\dd t-V(x)&<\ell,
    &&x\in\R\setminus[A,B].
 \label{eq:equilibrium-EL}
 \end{align}
\end{enumerate}
\end{Assumption}

Define 
\begin{equation}
 \kappa:=\bigl(\pi\sqrt{B-A}\,h(B)\bigr)^{2/3},
\label{eq:kappa-definition}
\end{equation}
and 
\begin{equation}
 \varepsilon_N:=(\kappa N^{2/3})^{-1},
\label{eq:edge-scale}
\end{equation}
the rescaled eigenvalues are
\begin{equation}
 \zeta_{N,j}:=\frac{\lambda_j-B}{\varepsilon_N},\qquad 1\leq j\leq N.
\label{eq:finite-edge-position}
\end{equation}
As $N\to\infty$, the point process formed by the rescaled eigenvalues
converges in distribution to the Airy point process~\cite{DKMVZ99,PS03}.
Denote its ordered particles by
$\zeta_1>\zeta_2>\cdots$ and define the rescaled finite-$N$ gap and the
limiting Airy gap by
\begin{equation}
\Delta_{N,k}:=\zeta_{N,k}-\zeta_{N,k+1},
\label{eq:finite-gap-definition}
\end{equation}
and
\begin{equation}
\Delta_{\mathrm{Ai},k}:=\zeta_k-\zeta_{k+1},
\label{eq:airy-gap-definition}
\end{equation}
respectively.

\pagebreak[3]
\subsubsection*{Right tails and conditional distributions}
\label{sec:right-results}

For $z\in\mathbb C\setminus(-\infty,B]$ and $s>0$, define
\begin{align}
 g(z)&:=\int_A^B\log(z-x)\rho(x)\,\dd x,\notag\\
 \Phi(B+s)&:=V(B+s)+\ell-2g(B+s),
 \label{eq:exterior-phase}\\
 \Gamma(s)&:=\frac{\Phi'(B+s)}{\kappa^{3/2}\sqrt{s}}.\notag
\end{align}
Equations \eqref{eq:equilibrium-density} and
\eqref{eq:equilibrium-EL-support}--\eqref{eq:equilibrium-EL}
give $\Phi(B)=0$, $\Phi(B+s)>0$, and, for sufficiently small
$s>0$,
\[
 \Gamma(s)=2\frac{h(B+s)}{h(B)}
 \sqrt{1+\frac{s}{B-A}}.
\]
This representation extends $\Gamma$ analytically to $s=0$, with
$\Gamma(0)=2$.

\begin{Theorem}
\label{thm:fixed-k-right-tail}
As $a\to\infty$,
\begin{align}
 \log \Prob(\Delta_{\mathrm{Ai},k}>a)
 ={}&-\frac{4k}{3}a^{3/2}
  +\frac{8\sqrt2}{3}k^{3/2}a^{3/4}
  -\frac{48k^2-11}{32}\log a
  +\log C_k+o(1),
\label{eq:fixed-k-airy-right}
\end{align}
where
\begin{align}
 C_k
 :={}&2^{-3k^2-k+1+5/48}e^{\zeta'(-1)}
 \pi^{1/2-k}k^{-5/16}
 \prod_{j=0}^{k-2}j!(j+2)!,
\label{eq:def-CRk}
\end{align}
and $\zeta$ is the Riemann zeta function \cite[Section~25.2]{DLMF}.
For $k=1$, the empty product in \eqref{eq:def-CRk} is understood to be $1$.

For the ensemble \eqref{eq:unitary-density} and every sequence $a_N\to\infty$
with $a_NN^{-2/3}\to0$, we have
\begin{align}
 \log \Prob(\Delta_{N,k}>a_N)
 ={}&-kN\Phi(B+a_N\varepsilon_N)
  +\frac43k^{3/2}\Gamma(a_N\varepsilon_N)^{3/2}a_N^{3/4}
  -\frac{48k^2-11}{32}\log a_N
 \notag\\
  &+\log C_k+o(1).
\label{eq:fixed-k-finite-right}
\end{align}
In the rest of the paper, $a_N$ is as in
Theorem~\ref{thm:fixed-k-right-tail}.
\end{Theorem}

\begin{Remark*}
The finite-$N$ and Airy right-tail probabilities need not be
asymptotically equivalent in this regime.
For the GUE, we have
$[A,B]=[-2,2]$, $\kappa=1$, and
\begin{align*}
 \Phi(2+s)&=\int_0^s\sqrt{t(t+4)}\,\dd t,
 &\Gamma(s)&=\sqrt{4+s}.
\end{align*}
Taking $a_N=N^{1/3}$ gives
\[
 \log\frac{\Prob(\Delta_{N,k}>N^{1/3})}
 {\Prob(\Delta_{\mathrm{Ai},k}>N^{1/3})}
 =-\frac{k}{10}N^{1/6}+o(1).
\]
\end{Remark*}

\begin{Corollary}
\label{cor:fixed-k-conditional-joint}
\label{cor:airy-conditional-joint}
Fix $k\geq2$.
Write $\mathbf s=(s_1,\ldots,s_{k-1})$ and define
\begin{equation}
\begin{aligned}
  \mathcal S_k&:=\{(s_1,\ldots,s_{k-1}):s_1>\cdots>s_{k-1}>0\}.
\end{aligned}
\label{eq:laguerre-chamber}
\end{equation}
Define the density on $\mathbb R\times\mathcal S_k$ by
\begin{equation*}
  \varrho_k(z,\mathbf s):=\frac{e^{-z^2/2}}{\sqrt{2\pi}}\,
  \frac{\displaystyle
  \prod_{1\le i<j\le k-1}(s_i-s_j)^2
  \prod_{i=1}^{k-1}s_i^2 e^{-\sum_{i=1}^{k-1}s_i}}
  {\displaystyle\prod_{j=0}^{k-2}j!(j+2)!}.
\end{equation*}
It is the product of a standard Gaussian density and the density of the
ordered $(k-1)$-particle Laguerre unitary ensemble with parameter $\alpha=2$.

For the Airy point process and $a>0$, define
\begin{equation}
 m_k(a):=a-2\sqrt{2k}\,a^{1/4},\qquad
 \sigma_k(a):=(2k)^{-1/4}a^{-1/8}.
\label{eq:airy-conditional-scales}
\end{equation}
Conditionally on $\{\Delta_{\mathrm{Ai},k}>a\}$, let $\varrho_{a,k}$ denote
the density of
\[
 \left(\frac{\zeta_k-m_k(a)}{\sigma_k(a)},
 2a^{1/2}(\zeta_1-\zeta_k),\ldots,
  2a^{1/2}(\zeta_{k-1}-\zeta_k)\right).
\]
As $a\to\infty$, $\varrho_{a,k}$ converges to $\varrho_k$ in
$L^1(\mathbb R\times\mathcal S_k)$. Under the same conditioning, we also have
\begin{equation}
 \left(\frac{\zeta_1-m_k(a)}{\sigma_k(a)},\ldots,
 \frac{\zeta_k-m_k(a)}{\sigma_k(a)}\right)
 \Longrightarrow (Z,\ldots,Z),\qquad Z\sim\mathcal N(0,1).
\label{eq:airy-common-gaussian}
\end{equation}

For the ensemble \eqref{eq:unitary-density} and every sequence
$a_N\to\infty$ with $a_N=o(N^{2/3})$, define
\begin{equation}
\begin{aligned}
  m_{N,k}&:=a_N-2\sqrt{k\Gamma(a_N\varepsilon_N)}\,a_N^{1/4},
  \qquad
  \sigma_{N,k}:=(k\Gamma(a_N\varepsilon_N))^{-1/4}a_N^{-1/8},
\end{aligned}
\label{eq:finite-conditional-scales}
\end{equation}
and
\begin{align}
  Z_{N,k}&:=\frac{\zeta_{N,k}-m_{N,k}}{\sigma_{N,k}},\notag\\
  S_{N,i}&:=\Gamma(a_N\varepsilon_N)a_N^{1/2}
  (\zeta_{N,i}-\zeta_{N,k}),
  \qquad 1\leq i<k.
\label{eq:finite-conditional-variables}
\end{align}
Conditionally on $\{\Delta_{N,k}>a_N\}$, let $\varrho_{N,k}$ denote the density
of $(Z_{N,k},S_{N,1},\ldots,S_{N,k-1})$. As $N\to\infty$, $\varrho_{N,k}$
converges to $\varrho_k$ in $L^1(\mathbb R\times\mathcal S_k)$.
Under the same conditioning, we also have
\begin{equation}
  \left(\frac{\zeta_{N,1}-m_{N,k}}{\sigma_{N,k}},\ldots,
  \frac{\zeta_{N,k}-m_{N,k}}{\sigma_{N,k}}\right)
 \ \Longrightarrow\ (Z,\ldots,Z),
 \qquad Z\sim\mathcal N(0,1).
\label{eq:fixed-k-common-gaussian}
\end{equation}
Equivalently, the conditional laws with densities $\varrho_{a,k}$ and
$\varrho_{N,k}$ converge in total variation to the law with density $\varrho_k$.
\end{Corollary}

The right-tail asymptotics yield a corresponding result for the maximum
of the first $k$ edge gaps. For
$\star\in\{\mathrm{Ai},N\}$, define
\begin{equation*}
 \Delta_{\star,k}^{\max}:=\max_{1\leq j\leq k}\Delta_{\star,j}.
\end{equation*}

\begin{Corollary}
\label{cor:finite-block-large-gaps}
Conditional on the maximum of the first $k$ Airy gaps exceeding $a$, the
first gap attains this maximum with probability tending to one as
$a\to\infty$:
\begin{equation}
 \Prob(\Delta_{\mathrm{Ai},1}=\Delta_{\mathrm{Ai},k}^{\max}
       \mid\Delta_{\mathrm{Ai},k}^{\max}>a)
 \longrightarrow1.
\label{eq:airy-block-maximum}
\end{equation}
The same conclusion holds for the ensemble \eqref{eq:unitary-density}
as $N\to\infty$:
\begin{equation}
 \Prob(\Delta_{N,1}=\Delta_{N,k}^{\max}
       \mid\Delta_{N,k}^{\max}>a_N)
 \longrightarrow1.
\label{eq:finite-block-maximum}
\end{equation}
\end{Corollary}

\begin{proof}
For each $2\leq j\leq k$, \eqref{eq:fixed-k-airy-right} gives
\[
 \log\frac{\Prob(\Delta_{\mathrm{Ai},j}>a)}
               {\Prob(\Delta_{\mathrm{Ai},1}>a)}
 =-\frac{4(j-1)}3a^{3/2}+O(a^{3/4}).
\]
In the finite-$N$ ensemble, \eqref{eq:exterior-phase} yields
\[
 N\Phi(B+a_N\varepsilon_N)
 =N\kappa^{3/2}\int_0^{a_N\varepsilon_N}\sqrt{t}\,\Gamma(t)\,\dd t
 =\frac43a_N^{3/2}+o(a_N^{3/2}).
\]
Since $\Gamma(a_N\varepsilon_N)\to2$, \eqref{eq:fixed-k-finite-right} likewise implies
\[
 \frac{\Prob(\Delta_{N,j}>a_N)}
      {\Prob(\Delta_{N,1}>a_N)}\longrightarrow0,
 \qquad 2\leq j\leq k.
\]
For either model and $t>0$, the union bound gives
\[
 1-\Prob(\Delta_{\star,1}=\Delta_{\star,k}^{\max}
       \mid\Delta_{\star,k}^{\max}>t)
 \leq\sum_{j=2}^{k}
 \frac{\Prob(\Delta_{\star,j}>t)}{\Prob(\Delta_{\star,1}>t)}.
\]
Using these limits with $(\star,t)=(\mathrm{Ai},a)$
and $(\star,t)=(N,a_N)$ now yields
\eqref{eq:airy-block-maximum} and \eqref{eq:finite-block-maximum}.
\end{proof}

\subsubsection*{Left tails and conditional distributions}
\label{sec:left-results}

\begin{Theorem}
\label{thm:fixed-k-left-tail}
For the Airy point process,
\begin{align}
 \Prob(\Delta_{\mathrm{Ai},k}<b)
 =\frac{k}{6}b^3+O(b^4),
 \qquad b\downarrow0.
\label{eq:fixed-k-airy-left}
\end{align}
Moreover, for the ensemble \eqref{eq:unitary-density} and every sequence
$b_N\downarrow0$ as $N\to\infty$,
\begin{align}
 \Prob(\Delta_{N,k}<b_N)
 =\frac{k}{6}b_N^3[1+o(1)].
\label{eq:fixed-k-finite-left}
\end{align}
\end{Theorem}

The Airy point process has correlation kernel
\begin{equation}
K_{\mathrm{Ai}}(x,y):=\frac{\operatorname{Ai}(x)\operatorname{Ai}'(y)-\operatorname{Ai}'(x)\operatorname{Ai}(y)}{x-y},
\label{eq:airy-kernel-main}
\end{equation}
where $\operatorname{Ai}$ denotes the Airy function \cite[Section~9.2]{DLMF}.
For $\mu\in\R$, let $K_{\mathrm{Ai},\mu}$ be the integral operator with kernel
\eqref{eq:airy-kernel-main} on $L^2((\mu,\infty),\dd x)$, so that
$F_2(\mu)=\mathbb P(\zeta_1\le \mu)=\det(I-K_{\mathrm{Ai},\mu})$
\cite[(1.19)]{TW94}.
Define the resolvent operator
\begin{equation}
 L_\mu:=K_{\mathrm{Ai},\mu}(I-K_{\mathrm{Ai},\mu})^{-1},
\label{eq:airy-resolvent-definition}
\end{equation}
and denote its kernel by $L_\mu(x,y)$.
We write $\partial_1$ and $\partial_2$ for partial derivatives with respect
to the first and second variables, respectively.
Let $\mathbf z=(z_1,\ldots,z_{k-1})$ with
$z_1>\cdots>z_{k-1}>\mu$. When $k=1$, the $\mathbf z$-variables and
integrals over them are omitted. Define
\begin{align}
 \mathcal M_{k,\mu}(\mathbf z)
 &:={}
 \begin{pmatrix}
 [L_\mu(z_i,z_j)]_{i,j=1}^{k-1}
   &[\partial_2L_\mu(z_i,\mu)]_{i=1}^{k-1}
   &[L_\mu(z_i,\mu)]_{i=1}^{k-1}\\
 [\partial_1L_\mu(\mu,z_j)]_{j=1}^{k-1}
   &\partial_1\partial_2L_\mu(\mu,\mu)&\partial_1L_\mu(\mu,\mu)\\
 [L_\mu(\mu,z_j)]_{j=1}^{k-1}
   &\partial_2L_\mu(\mu,\mu)&L_\mu(\mu,\mu)
 \end{pmatrix},
\notag \\
 q_k(\mu)
 &:={}
 F_2(\mu)\int_{z_1>\cdots>z_{k-1}>\mu}
 \det\mathcal M_{k,\mu}(\mathbf z)
 \prod_{j=1}^{k-1}\dd z_j.
\label{eq:def-left-qk}
\end{align}


\begin{Corollary}
\label{cor:fixed-k-left-collision-law}
Define the probability measure $\nu_k$ on
\begin{align*}
 \mathcal C_k
 :={}&
 \begin{cases}
 \{(\mathbf z,\mu)\in\R^k:z_1>\cdots>z_{k-1}>\mu\},&k\geq2,\\
 \R,&k=1
 \end{cases}
\end{align*}
by
\begin{equation}
\begin{aligned}
 \dd\nu_k(\mathbf z,\mu)
 :={}&\frac2k F_2(\mu)\det\mathcal M_{k,\mu}(\mathbf z)
 \prod_{j=1}^{k-1}\dd z_j\,\dd \mu.
\end{aligned}
\label{eq:collision-measure}
\end{equation}
Let $(\mathbf Z,U)$ have law $\nu_k$, where
$\mathbf Z=(Z_1,\ldots,Z_{k-1})$ for $k\geq2$ and $\mathbf Z$ is empty for
$k=1$, and let $T$ be independent of $(\mathbf Z,U)$ with density
$3t^2\mathbf1_{(0,1)}(t)$.
Writing $\mathcal L(X)$ for the law of a random variable $X$, we have,
as $b\downarrow0$,
\begin{align}
 \mathcal L\left(
 \left(\zeta_1,\ldots,\zeta_k,
 \frac{\Delta_{\mathrm{Ai},k}}{b}\right)
 \ \middle|\ \Delta_{\mathrm{Ai},k}<b\right)
 \longrightarrow
 \mathcal L(\mathbf Z,U,T),
\label{eq:airy-left-collision-law}
\end{align}
and for the ensemble \eqref{eq:unitary-density} and every
sequence $b_N\downarrow0$ as $N\to\infty$,
\begin{align}
 \mathcal L\left(
 \left(\zeta_{N,1},\ldots,\zeta_{N,k},
 \frac{\Delta_{N,k}}{b_N}\right)
 \ \middle|\ \Delta_{N,k}<b_N\right)
 \longrightarrow
 \mathcal L(\mathbf Z,U,T),
\label{eq:finite-left-collision-law}
\end{align}
with both limits in total variation.
\end{Corollary}

The left-tail asymptotics and conditional limits yield the corresponding
result for the minimum of the first $k$ edge gaps. For
$\star\in\{\mathrm{Ai},N\}$, define
\begin{equation*}
\begin{aligned}
 \Delta_{\star,k}^{\min}
 &:=\min_{1\leq j\leq k}\Delta_{\star,j},\\
 J_{\star,k}
 &:=\min\{j\in\{1,\ldots,k\}:
       \Delta_{\star,j}=\Delta_{\star,k}^{\min}\}.
\end{aligned}
\end{equation*}

\begin{Corollary}
\label{cor:finite-block-small-gaps}
For the Airy point process, when $b\downarrow0$, we have:
\begin{align*}
 \Prob(\Delta_{\mathrm{Ai},k}^{\min}<b)
 &=\frac{k(k+1)}{12}b^3+o(b^3),
 \\
\intertext{and the index of the minimum satisfies}
 \Prob(J_{\mathrm{Ai},k}=j\mid\Delta_{\mathrm{Ai},k}^{\min}<b)
 &\longrightarrow
 \frac{2j}{k(k+1)},
 \qquad 1\leq j\leq k.
\end{align*}
The same tail asymptotics and conditional index distribution hold for
the ensemble \eqref{eq:unitary-density} and every sequence $b_N\downarrow0$
as $N\to\infty$:
\begin{align*}
 \Prob(\Delta_{N,k}^{\min}<b_N)
 &=\frac{k(k+1)}{12}b_N^3+o(b_N^3),
 \\
 \Prob(J_{N,k}=j\mid\Delta_{N,k}^{\min}<b_N)
 &\longrightarrow
 \frac{2j}{k(k+1)},
 \qquad 1\leq j\leq k.
\end{align*}
\end{Corollary}

\begin{proof}
Set $A_j(b):=\{\Delta_{\mathrm{Ai},j}<b\}$ for $1\leq j\leq k$.
For $i<j$, \eqref{eq:airy-left-collision-law} and
\eqref{eq:fixed-k-airy-left} give
\[
 \Prob(A_i(b)\cap A_j(b))
 =\Prob(A_i(b)\mid A_j(b))\Prob(A_j(b))
 =o(1)\,O(b^3)=o(b^3).
\]
Since $k$ is fixed, the preceding intersection estimate and
inclusion--exclusion give
\begin{align*}
 \Prob(\Delta_{\mathrm{Ai},k}^{\min}<b)
 &=\sum_{j=1}^{k}\Prob(A_j(b))+o(b^3)
 =\frac{k(k+1)}{12}b^3+o(b^3),\\
 \Prob(J_{\mathrm{Ai},k}=j,\Delta_{\mathrm{Ai},k}^{\min}<b)
 &=\Prob(A_j(b))+o(b^3)
 =\frac{j}{6}b^3+o(b^3).
\end{align*}
Taking the ratio gives the conditional index limit.
Both finite-$N$ conclusions follow by the same argument, using
\eqref{eq:fixed-k-finite-left} and \eqref{eq:finite-left-collision-law}.
\end{proof}

\subsection{Interpretation of the conditional laws}
\label{subsec:conditional-interpretation}

The conditional laws describe how the particles are arranged when a gap is
unusually large or small. We first illustrate the large-gap mechanism using
the case $k=1$ in the Airy point process.
A large first gap can be formed by moving the largest point $\zeta_1$
to the right and the remaining points $\zeta_2,\zeta_3,\ldots$ to the left.
If $\zeta_1=a-t$, the condition $\Delta_{\mathrm{Ai},1}>a$
is equivalent to $\zeta_2<-t$, so all points except $\zeta_1$ lie below $-t$.
For $0<t\ll a$, the leading right-tail exponent at $a-t$ is
$-\frac43(a-t)^{3/2}\approx-\frac43a^{3/2}+2\sqrt a\,t$.
Keeping all other points below $-t$ adds the left-tail term $-t^3/12$.
Combining these two effects suggests the rough
probability weight
\[
 e^{-\frac43a^{3/2}+2\sqrt a\,t-\frac{t^3}{12}}=e^{-\frac43a^{3/2}+\frac{8\sqrt2}{3}a^{3/4}}\exp\left[-\frac{t+4\sqrt2\,a^{1/4}}{12}\left(t-2\sqrt2\,a^{1/4}\right)^2\right],
\]
with Gaussian fluctuations of the largest particle arising from
the locally quadratic form of the exponent.

As $a\downarrow0$, the conditioning disappears and the largest point has
its usual Tracy--Widom law. The change from Tracy--Widom to Gaussian
fluctuations resembles the rank-one BBP transition~\cite{BBP05}.
Write
\[
 H_a(x):=\Prob(\zeta_1\le x\mid\Delta_{\mathrm{Ai},1}>a),
\]
and let $B_w$ denote the rank-one BBP distribution function in the same
Airy normalization. For fixed $a>0$ and $w\in\R$, the following asymptotics
hold as $r\to\infty$:\footnote{%
By \eqref{eq:conditional-painleve-density}, the bound
$L_s(x,x)\ge K_{\mathrm{Ai}}(x,x)$ for $x>s$, integrated on
$[x,\infty)$ or $[x-1,x]$, gives lower bounds for the two tails;
removing the gap condition gives upper bounds.
The Airy kernel and Tracy--Widom asymptotics then give the $H_a$ estimates.
The BBP estimates follow from the scalar rank-one representation and its
$x$-equation in \cite{Baik2006}, together with the fixed-parameter
asymptotics in \cite{BR00}.}
\[
 \begin{array}{c|cc}
  &\text{Left tail}&\text{Right tail}\\[2pt]
  H_a&\log H_a(-r)\sim-r^3/12
     &\log(1-H_a(r))\sim-4r^{3/2}/3\\[2pt]
  B_w&\log B_w(-r)\sim-r^3/12
     &\log(1-B_w(r))\sim-2r^{3/2}/3
 \end{array}
\]

\begin{Remark}[Painlev\'e representation]
\label{rem:conditional-painleve}
For every fixed $k\geq1$ and $a>0$, the conditional joint density of the
leading $k$ Airy particles given $\{\Delta_{\mathrm{Ai},k}>a\}$ has the
representation \eqref{eq:conditional-painleve-joint-density} in terms of the
Hastings--McLeod solution of Painlev\'e II.
Let $q$ denote this solution, so that
\begin{equation}
 q''(u)=u q(u)+2q(u)^3,\qquad
 q(u)\sim\Ai(u)\quad(u\to+\infty).
\label{eq:conditional-hm-solution}
\end{equation}
For $b\geq0$, let $f_b$ solve the boundary-value problem~\cite[Section~4.1]{PerretSchehr}
\begin{equation}
 \partial_u^2 f_b(u)=\bigl(u+b+2q(u)^2\bigr)f_b(u),\qquad
 f_b(u)\sim\Ai(u+b)\quad(u\to+\infty).
\label{eq:conditional-linear-equation}
\end{equation}

It follows from \cite[Sections~II--III]{TW94} that, for $s\in\R$ and $b,c>0$,
\begin{equation}
 L_s(s+b,s+c)=\int_s^\infty f_b(u)f_c(u)\,\dd u.
\label{eq:conditional-resolvent-gram}
\end{equation}

For $t_1>\cdots>t_k$, set $s:=t_k-a$ and $b_i:=t_i-s$,
$1\leq i\leq k$, so that $b_k=a$.
The ordered Janossy density on $(s,\infty)$ and
\eqref{eq:conditional-resolvent-gram}, together with \eqref{eq:fk-3-23}, give the joint density of
$(\zeta_1,\ldots,\zeta_k)$, conditional on $\{\Delta_{\mathrm{Ai},k}>a\}$,
with respect to $\dd t_1\cdots\dd t_k$:
\begin{equation}
 p_{k,a}(t_1,\ldots,t_k)
 =\frac{F_2(s)}{\Prob(\Delta_{\mathrm{Ai},k}>a)}
  \det\left[\int_s^\infty f_{b_i}(u)f_{b_j}(u)\,\dd u\right]_{i,j=1}^k.
\label{eq:conditional-painleve-joint-density}
\end{equation}

For $k=1$, define
\begin{equation}
 r_a(s):=\int_s^\infty f_a(u)^2\,\dd u.
\label{eq:conditional-resolvent-integral}
\end{equation}
Equation~\eqref{eq:conditional-painleve-joint-density} reduces to
\begin{equation}
 \Prob(\Delta_{\mathrm{Ai},1}>a)H_a'(x)
 =F_2(x-a)L_{x-a}(x,x)=F_2(x-a)r_a(x-a),
\label{eq:conditional-painleve-density}
\end{equation}
and hence
\begin{equation}
 H_a(x)=
 \frac{\displaystyle\int_{-\infty}^{x-a}F_2(s)r_a(s)\,\dd s}
 {\displaystyle\int_{-\infty}^{\infty}F_2(s)r_a(s)\,\dd s}.
\label{eq:conditional-painleve-cdf}
\end{equation}

\end{Remark}

For fixed $k\ge2$, the same balance of escape and compression costs gives
a Gaussian fluctuation shared by the first $k$ particles.
To explain the relative positions under $\Delta_{\mathrm{Ai},k}>a$ as
$a\to\infty$, write $\zeta_k=a-t$ and $d_i=\zeta_i-\zeta_k$ for $i<k$.
For $0<t\ll a$ and $d_i\ll a$, moving the $i$th particle a distance $d_i$
farther to the right changes its right-tail exponent by
\[
 -\frac43(a-t+d_i)^{3/2}+\frac43(a-t)^{3/2}
 \approx-2\sqrt a\,d_i.
\]
Mutual repulsion among the first $k-1$ points gives the factors
$(d_i-d_j)^2$, and repulsion from the $k$th point gives $d_i^2$.
Combining these with the exponential weights suggests a density
proportional to
\[
 \prod_{1\leq i<j\leq k-1}(d_i-d_j)^2
 \prod_{i=1}^{k-1}d_i^2
 \exp\left[-2\sqrt a\sum_{i=1}^{k-1}d_i\right].
\]
Rescaling the distances by $2\sqrt a$ gives the ordered
$(k-1)$-particle Laguerre density with $\alpha=2$.

As the $k$th and $(k+1)$st particles approach a common location,
their quadratic repulsion contributes a factor
$\Delta_{\mathrm{Ai},k}^2$ to the leading joint density; integrating over
gaps in $(0,b)$ explains the cubic small-gap probability.
The collision location and the higher particles have limiting joint law
$\nu_k$, as in Corollary~\ref{cor:fixed-k-left-collision-law}.
For the general potentials considered here, the preceding corollaries give
the same limiting laws after the stated normalization.

\subsection{Proof strategy}
The eigenvalues form a determinantal point process
\cite[Section~1.2]{DKMVZ99} with correlation kernel $K_N(x,z)$.
For $y\in\R$, let $K_{N,y}$ be the integral operator with kernel $K_N$
on $L^2((y,\infty),\dd x)$, and set
\begin{align}
 D_N(y)&:=\det(I-K_{N,y})=\Prob(\lambda_1\leq y),
 &L_{N,y}&:=K_{N,y}(I-K_{N,y})^{-1},
\label{eq:finite-janossy-data}
\end{align}
where $L_{N,y}(x,z)$ denotes the resolvent kernel. For $k<N$ and
$x_1>\cdots>x_k>y$, the ordered Janossy density
\cite[p.~140]{DVJ1988} specifies the top $k$ eigenvalues at
$(x_1,\ldots,x_k)$ and all remaining eigenvalues below $y$:
\begin{equation}
 J^{(k)}_{N,y}(x_1,\ldots,x_k)
 :=D_N(y)\det[L_{N,y}(x_i,x_j)]_{i,j=1}^k.
\label{eq:fixed-k-janossy-density}
\end{equation}

For $a>0$, set $y=\lambda_k-a\varepsilon_N$,
with $\varepsilon_N$ as in \eqref{eq:edge-scale}.
By \eqref{eq:finite-gap-definition},
\[
 \Delta_{N,k}>a\quad\Longleftrightarrow\quad\lambda_{k+1}<y<\lambda_k.
\]
Integrating the ordered Janossy density
\eqref{eq:fixed-k-janossy-density} gives
\begin{equation}
 \Prob(\Delta_{N,k}>a)
 =\int_{\R}\int_{x_1>\cdots>x_{k-1}>y+a\varepsilon_N}
 J^{(k)}_{N,y}(x_1,\ldots,x_{k-1},y+a\varepsilon_N)
 \prod_{i=1}^{k-1}\dd x_i\,\dd y.
\label{eq:finite-fixed-gap-janossy}
\end{equation}
The same argument for the Airy gap \eqref{eq:airy-gap-definition}
yields
\begin{equation}
\mathbb P(\Delta_{\mathrm{Ai},k}>a)
=\int_{\mathbb R}F_2(\mu)
\int_{x_1>\cdots>x_{k-1}>\mu+a}
\det[L_\mu(x_i,x_j)]_{i,j=1}^k
\prod_{i=1}^{k-1}\mathrm dx_i\,\mathrm d\mu,
\label{eq:fk-3-23}
\end{equation}
where $x_k=\mu+a$ and $L_\mu$ is defined in
\eqref{eq:airy-resolvent-definition}.

Borodin and Soshnikov~\cite[Theorem~1]{BS} give a Christoffel--Darboux
formula for $L_{N,y}$ using polynomials orthogonal on $(-\infty,y]$
with weight $e^{-NV(x)}$. For fixed $\mu,x,z\in\R$ with $x,z>\mu$,
Rider and Zhou~\cite[Theorem~1.1]{RZ} proved
\[
 \varepsilon_NL_{N,B+\varepsilon_N\mu}
 (B+\varepsilon_Nx,B+\varepsilon_Nz)
 \longrightarrow L_\mu(x,z),\qquad N\to\infty.
\]
The tail analysis requires uniform estimates as the truncation point and
particle positions vary with the gap threshold.

For the Airy right tail as $a\to\infty$, write the truncation point as $-\mu$.
We use the RH asymptotics of $L_{-\mu}$ and divided differences to obtain
the leading term of the kernel determinant in \eqref{eq:fk-3-23}, and then
integrate this term over the relative particle positions using the
Laguerre integral formula. Combining the resulting expression with the
Tracy--Widom left-tail expansion~\cite{DIK} gives
the $\mu$-dependent phase $2k\sqrt a\,\mu-\mu^3/12$, whose maximum over
$\mu>0$ is at $\mu=2\sqrt{2k}\,a^{1/4}$. We localize the integral near this
point using kernel bounds and induction on $k$. The induction handles truncation
points far to the left: if the $k$th particle lies sufficiently far to
the left, then either the largest particle does too or one of the
preceding gaps is large. The largest-particle lower tail and the
lower-index gap estimates therefore control this range, leaving a
one-dimensional Laplace integral. For finite $N$, we analyze the orthogonal polynomials truncated at
$y=B-\mu\varepsilon_N$ using the constrained equilibrium measure and the
Airy resolvent model. The endpoint identity \cite[(1.28)]{RZ}
\begin{equation}
 \partial_y\log D_N(y)=L_{N,y}(y,y)
\label{eq:finite-endpoint-identity}
\end{equation}
and its Airy counterpart give the comparison \eqref{eq:det-comparison} by integration
from $y=B$, where \eqref{eq:det-anchor} controls the integration constant.
Together with the kernel determinant asymptotics, this reduces the proof
to the same integrations and truncation argument as in the Airy case.

For the Airy left tail, we use the Gram representation of the Janossy
density to factor out the squared gap. Integrating the leading coefficient
over the other particle positions gives the collision density $q_k$ in
\eqref{eq:def-left-qk}. Integrating over gaps below $b$ then yields the
cubic left-tail asymptotics.
The coefficient is $k/6$, by the mass identity
$\int_{\R}q_k(\mu)\,\dd\mu=k/2$
in Lemma~\ref{lem:collision-density-mass}, which we prove using a parameterized
confluent Fredholm minor and the Airy resolvent identities.
For finite $N$, $L^1$ convergence of the Janossy
density divided by the squared gap, together with a uniform mass bound,
yields the same coefficient along every joint sequence
$N\to\infty$, $b_N\downarrow0$.

We obtain the conditional densities by normalizing the joint particle
densities on the large- or small-gap event by the corresponding gap-tail
probabilities. The same asymptotic estimates show that the rescaled conditional
densities converge pointwise to densities of total mass one, so
Scheff\'e's lemma~\cite{Scheffe1947} gives convergence in total variation.

\subsection{Organization of the paper}

We begin in Section~\ref{sec:rh-analysis} with the RH analysis of the Airy
resolvent model (Section~\ref{subsec:right-airy-model}) and the truncated
orthogonal polynomials (Section~\ref{subsec:right-rhp-comparison}).
For the latter, we formulate the RH problem and use the constrained equilibrium
measure to carry out the RH transformations, expressing the resolvent kernel
in terms of the RH solution. Using this analysis, we prove the Airy and
finite-$N$ right-tail asymptotics in Sections~\ref{subsec:right-airy-tail}
and~\ref{subsec:right-finite-tail}, respectively, and derive the conditional
Gaussian--Laguerre laws in Section~\ref{subsec:right-conditional-limits}.
We then turn to the left tail in Section~\ref{sec:left-proof}, where we
evaluate the total mass of the collision density, prove the Airy and
finite-$N$ small-gap asymptotics, and derive the conditional limit laws for
both models. Appendices~\ref{airy} and~\ref{A} contain the Airy and Bessel
parametrices, respectively, used in the local RH constructions of
Section~\ref{sec:rh-analysis}.

\subsection*{Notation}

\begin{itemize}

\item Unless otherwise specified, $c,C>0$ denote constants whose values may
change between estimates. For positive $f,g$, we write $f=O(g)$ if $f\leq Cg$,
and $f\asymp g$ if $cg\leq f\leq Cg$. Implicit constants may depend on
$V$, $k$, and locally fixed parameters; this dependence is usually suppressed.

\item Unless otherwise specified, logarithms and fractional powers are taken on
their principal branches in $\mathbb C\setminus(-\infty,0]$.

\item For $z_0\in\mathbb C$ and $r>0$, define
\begin{equation*}
	D(z_0, r) := \{ z\in \mathbb{C} \mid |z-z_0|<r \}.
\end{equation*}
When the boundary $\partial D(z_0,r)$ is used as an RH contour, it is
oriented clockwise. We write $K\Subset U$ when $\overline K$ is compact
and contained in $U$.

\item We use the Pauli matrix
\begin{equation*}
	\sigma_3:=
	\begin{pmatrix}
		1 & 0 \\
		0 & -1
	\end{pmatrix},
\end{equation*}
and write $I$ for the identity matrix or operator and $E_{ij}$ for the
matrix with a single $1$ in entry $(i,j)$ and zeros elsewhere.

\item For $\mathbf x=(x_1,\ldots,x_m)$, the Vandermonde product is
\[
 \Delta(\mathbf x)
 :=\prod_{1\leq i<j\leq m}(x_i-x_j).
\]
\end{itemize}

\subsection*{Acknowledgements}

I thank Zhigang Bao for bringing this problem to my attention and for
subsequent discussions and help. I also thank Lun Zhang for helpful
suggestions and discussions on this work, and Guilherme Silva for
introducing me to Janossy densities while he was visiting Lun Zhang at
Fudan University. The author is partially supported by Hong Kong RGC Grant GRF 17304225 and 17303826.

\subsection*{AI claim}

I have used OpenAI models to assist
with some calculations and to polish this manuscript. I take full
responsibility for the content of this paper.

\newpage
\section{Riemann--Hilbert Analysis}
\label{sec:rh-analysis}

\subsection{Asymptotics of the Airy resolvent}
\label{subsec:right-airy-model}

We first analyze the Airy resolvent model; the resulting estimates enter
the Airy gap asymptotics and the finite-$N$ RH construction below.

\subsubsection{Model RH problem and kernel formula}

For $\mu\geq0$, the Airy resolvent kernel $L_{-\mu}$ defined in
\eqref{eq:airy-resolvent-definition} can be expressed in terms of the solution
$M(\,\cdot\,;-\mu)$ of the following Rider--Zhou model RH problem \cite{RZ}.

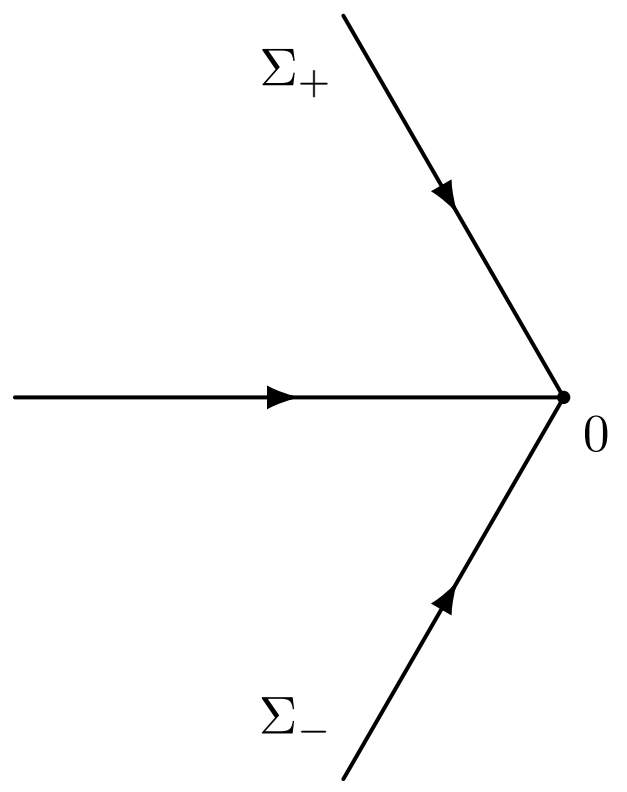
\begin{figure}[htbp]
\centering
\begin{tikzpicture}[rhp diagram,x=1.1cm,y=1.1cm]
 \draw[rhp contour] (-3.3,0)--(0,0);
 \draw[rhp contour] (120:2.65)--(0,0);
 \draw[rhp contour] (-120:2.65)--(0,0);
 \fill (0,0) circle (1.25pt);
 \node[below right=2pt] at (0,0) {$0$};
 \node[left=6pt] at (120:2.25) {$\Sigma_+$};
 \node[left=6pt] at (-120:2.25) {$\Sigma_-$};
\end{tikzpicture}
\caption{Contours $(-\infty,0]$ and $\Sigma_\pm$ in the RH problem for $M$.}
\label{fig:Mminus-contour}
\end{figure}

\begin{RHProblem}\label{rhp:Mminus}
\hfill\break
\begin{enumerate}
    \item[(a)] $M(\,\cdot\,;-\mu)$ is defined and analytic in
    $\mathbb C\setminus((-\infty,0]\cup\Sigma_+\cup\Sigma_-)$,
    where $\Sigma_\pm$ are the oriented rays shown in
    Figure~\ref{fig:Mminus-contour}.
    \item[(b)] For $\zeta\in(-\infty,0)\cup\Sigma_+\cup\Sigma_-$, we have
    \[
        M_+(\zeta;-\mu)=M_-(\zeta;-\mu)
        \begin{cases}
        \begin{pmatrix}0&1\\-1&0\end{pmatrix},
            &\zeta\in(-\infty,0),\\[3mm]
        I+e^{\theta(\zeta;-\mu)}E_{21},
            &\zeta\in\Sigma_+\cup\Sigma_-,
        \end{cases}
    \]
    where
    \begin{equation}
     \theta(\zeta;-\mu):=\frac43\zeta^{3/2}-2\mu\zeta^{1/2}.
    \label{eq:airy-model-phase}
    \end{equation}
    \item[(c)] As $\zeta\to\infty$ with $\zeta\in \mathbb C\setminus((-\infty,0]\cup\Sigma_+\cup\Sigma_-)$,
    \[
        M(\zeta;-\mu)=\left(I+O(\zeta^{-1})\right)
        \mathfrak m(\zeta),
        \qquad
        \mathfrak m(\zeta):=\zeta^{-\sigma_3/4}
        \frac1{\sqrt2}\begin{pmatrix}1&1\\-1&1\end{pmatrix}
        e^{-\pi i\sigma_3/4}.
    \]
    \item[(d)] As $\zeta\to0$, we have
    \[
        M(\zeta;-\mu)=O(\log\zeta).
    \]
\end{enumerate}
\end{RHProblem}

Section~6 of \cite{RZ} proves unique solvability of this problem for every
$\mu>0$. The solution extends continuously to $\mu=0$ by
\cite[Section~7.1, Proposition~7.1 and Lemma~7.2]{RZ}. For $\mu\geq0$, define
the vector
\begin{equation}
    v_\mu(\zeta)
    :=\frac{e^{\pi i/4}}{\sqrt{2\pi}}
      e^{-\theta(\zeta;-\mu)/2}M(\zeta;-\mu)\binom{1}{0}.
\label{eq:airy-recovery-vector}
\end{equation}
The kernel formula in \cite[Theorem~1.1 and Definition~1.2]{RZ} gives, for $x,y>0$,
\begin{align}
 L_{-\mu}(x-\mu,y-\mu)
 &=\frac{\det(v_\mu(x),v_\mu(y))}{x-y},
 \label{eq:fk-2-3}\\
 L_{-\mu}(x-\mu,x-\mu)
 &=\det(v_\mu'(x),v_\mu(x)).
 \label{eq:airy-diagonal-kernel}
\end{align}

\subsubsection{Resolvent estimates}

To analyze the large-$\mu$ behavior, we factor out the leading term at
infinity by defining
\begin{equation}
 G_\mu(\zeta):=M(\zeta;-\mu)\mathfrak m(\zeta)^{-1},
\label{eq:airy-comparison-matrix}
\end{equation}
with $\mu\geq0$ and $\zeta\in\mathbb C\setminus((-\infty,0]\cup\Sigma_+\cup\Sigma_-)$.

\begin{Lemma}
\label{prop:outer-B}
For every fixed $\rho>0$, as $\zeta\to+\infty$,
\begin{align}
    L_{-\mu}(\zeta-\mu,\zeta-\mu)
    &=\frac{e^{-\theta(\zeta;-\mu)}}{8\pi\zeta}
      \left(1+O(\zeta^{-1/2})\right),
\label{eq:model-B-uniform}\\
    \|v_\mu(\zeta)\|^2
    &\leq C e^{-\theta(\zeta;-\mu)}\zeta^{1/2},
\label{eq:model-v-uniform}
\end{align}
uniformly for $\mu\geq0$ and $\zeta>\rho\mu$.
\end{Lemma}

\begin{proof}
Fix $\rho>0$ and $0<r<\min(\rho/2,1/4)$, and first let $\mu\to\infty$.
Set $w:=\zeta/\mu$ and define $\eta(w):=w(1-2w/3)^2$, so that the phase function $\theta$ in
\eqref{eq:airy-model-phase} satisfies
\[
    \theta(\mu w;-\mu)=-2(\mu^3\eta(w))^{1/2},
    \qquad |w|<r.
\]
As in \cite[Section~7.2, proof of Lemma~7.5]{RZ}, we define the local parametrix
\begin{equation}\label{eq:airy-bessel-parametrix}
    P_\mu(w):=\left(\pi\mu\left(1-\frac{2w}{3}\right)\right)^{\sigma_3/2}
    e^{-\pi i\sigma_3/4}\Phi^{(\Bes)}(\mu^3\eta(w))
    e^{\theta(\mu w;-\mu)\sigma_3/2},
\end{equation}
where $\Phi^{(\Bes)}$ is given by \eqref{def:Bespara}.
The corresponding error matrix is
\begin{equation}\label{eq:airy-rh-error-definition}
    \mathcal R_\mu(w):=
    \begin{cases}
      M(\mu w;-\mu)\mathfrak m(\mu w)^{-1},&|w|>r,\\
      M(\mu w;-\mu)P_\mu(w)^{-1},&|w|<r,
    \end{cases}
\end{equation}
with jump matrix
\begin{equation}\label{eq:airy-error-jump}
    J_{\mathcal R,\mu}(s):=
    \begin{cases}
    P_\mu(s)\mathfrak m(\mu s)^{-1},&s\in \partial D(0,r),\\
    \mathfrak m(\mu s)
    \left(I+e^{\theta(\mu s;-\mu)}E_{21}\right)
    \mathfrak m(\mu s)^{-1},&s\in\Gamma_+\cup\Gamma_-,
    \end{cases}.
\end{equation}
\begin{figure}[htbp]
\centering
\begin{tikzpicture}[rhp diagram,x=1.1cm,y=1.1cm]
 \draw[rhp contour] (120:2.8)--(120:1);
 \draw[rhp contour] (-120:2.8)--(-120:1);
 \draw[rhp circle] (1,0) arc[start angle=0,end angle=-360,radius=1];
 \fill (0,0) circle (1.25pt);
 \node[below right=2pt] at (0,0) {$0$};
 \node[left=6pt] at (120:2.35) {$\Gamma_+$};
 \node[left=6pt] at (-120:2.35) {$\Gamma_-$};
 \node[right] at (0.84,0.65) {$\partial D(0,r)$};
\end{tikzpicture}
\caption{Contours $\partial D(0,r)$ and $\Gamma_\pm$ in the RH problem for $\mathcal R_\mu$.}
\label{fig:bessel-error-contour}
\end{figure}
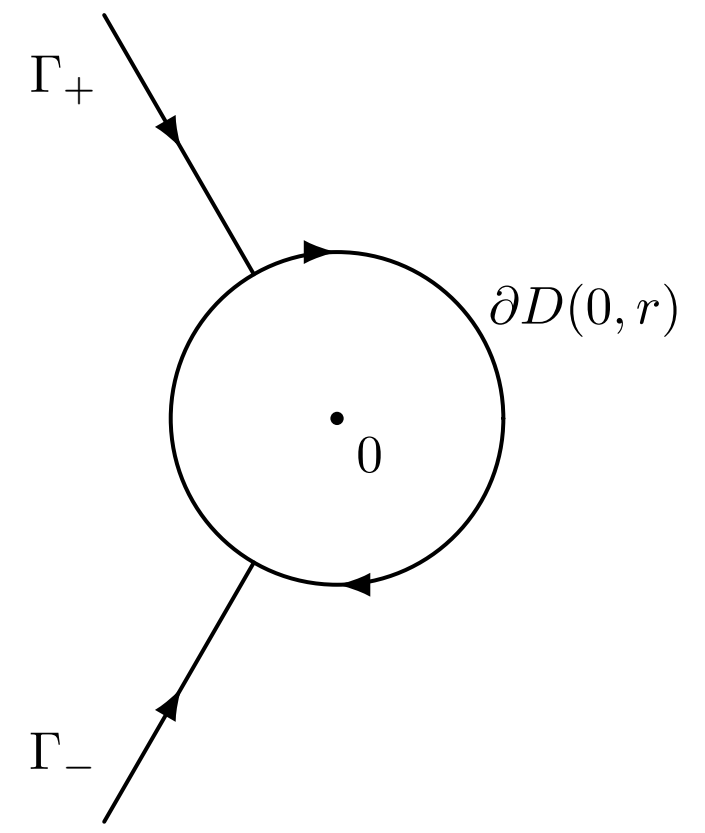

By \eqref{infty:Bes} and \eqref{eq:airy-model-phase},
\[
 J_{\mathcal R,\mu}(s)-I=
 \begin{cases}
 O(\mu^{-1}),&s\in\partial D(0,r),\\
 O(e^{-c\mu^{3/2}|s|^{3/2}}),&s\in\Gamma_+\cup\Gamma_-.
 \end{cases}
\]
Small-norm theory \cite{Deift1999,DZ93} therefore gives
\begin{equation}
 \mathcal R_\mu(w)=I+O(\mu^{-1}),
\label{eq:airy-rh-error-expansion}
\end{equation}
uniformly for $w$ at a fixed positive distance from
$\partial D(0,r)\cup\Gamma_+\cup\Gamma_-$, and
\[
 \|\mathcal R_{\mu,-}(J_{\mathcal R,\mu}-I)\|_
 {L^1(\partial D(0,r)\cup\Gamma_+\cup\Gamma_-)}=O(\mu^{-1}).
\]
Since $G_\mu(\zeta)=\mathcal R_\mu(\zeta/\mu)$, the preceding $L^1$
bound and the Cauchy representation give, for every integer $j\geq1$,
as $|\zeta|\to\infty$,
\begin{equation}
 G_\mu(\zeta)-I=O(|\zeta|^{-1}),
 \qquad
 \partial_\zeta^jG_\mu(\zeta)=O(|\zeta|^{-j-1}),
\label{eq:airy-G-bounds}
\end{equation}
uniformly for $\mu\geq0$ and $|\zeta|\geq\rho\mu$.
For bounded $\mu$, these estimates follow from
RH problem~\ref{rhp:Mminus}(c) and Cauchy's formula.
Using \eqref{eq:airy-comparison-matrix}, we rewrite
\eqref{eq:airy-recovery-vector} as
\begin{equation}
 v_\mu(\zeta)=G_\mu(\zeta)v_\mu^{(0)}(\zeta),
 \qquad
 v_\mu^{(0)}(\zeta):=
 \frac{e^{-\theta(\zeta;-\mu)/2}}{2\sqrt\pi}
 \begin{pmatrix}\zeta^{-1/4}\\-\zeta^{1/4}\end{pmatrix},
\label{eq:airy-leading-vector}
\end{equation}
where $\theta$ is defined in \eqref{eq:airy-model-phase}.
Applying \eqref{eq:airy-diagonal-kernel} to
\eqref{eq:airy-leading-vector} then gives
\begin{equation}
 L_{-\mu}(\zeta-\mu,\zeta-\mu)
 =\frac{e^{-\theta(\zeta;-\mu)}}{8\pi\zeta}
  +\det\bigl(G_\mu^{-1}G_\mu'v_\mu^{(0)},v_\mu^{(0)}\bigr).
\label{eq:airy-Wronskian-recovery}
\end{equation}
By \eqref{eq:airy-G-bounds}, the last term is
$O(e^{-\theta(\zeta;-\mu)}\zeta^{-3/2})$; the bound \eqref{eq:model-v-uniform}
follows from \eqref{eq:airy-leading-vector}.
\end{proof}

\begin{Lemma}\label{lem:wall-Wronskian}
For all $\mu\geq0$, we have
\[
    |L_{-\mu}(-\mu,-\mu)|\leq C(1+\mu^2),
    \qquad
    \|v_\mu(0)\|^2\leq C(1+\mu).
\]
\end{Lemma}

\begin{proof}
By continuity in $\mu$, it suffices to consider large $\mu$.
Substituting \eqref{eq:airy-bessel-parametrix} and
\eqref{eq:airy-rh-error-definition} into \eqref{eq:airy-recovery-vector}
gives, for $\zeta$ near the origin,
\[
    v_\mu(\zeta)=\frac{(\mu-\frac23\zeta)^{1/2}}{\sqrt2}
    \mathcal R_\mu(\zeta/\mu)
    \begin{pmatrix}
      I_0(-\theta(\zeta;-\mu)/2)\\
      -\zeta^{1/2}I_1(-\theta(\zeta;-\mu)/2)
    \end{pmatrix},
\]
where $\mathcal R_\mu$ is the error matrix defined in
\eqref{eq:airy-rh-error-definition}, and
$I_0$ and $I_1=I_0'$ are the modified Bessel functions in \eqref{def:Bespara}.
Equation~\eqref{eq:airy-rh-error-expansion} and Cauchy's formula on a smaller
disk give $\mathcal R_\mu(0)=I+O(\mu^{-1})$ and
$\mathcal R_\mu'(0)=O(\mu^{-1})$, where the prime denotes differentiation
in $w=\zeta/\mu$.
Expanding $I_0$ and $I_1$ at the origin
\cite[Equations~(10.25.2) and~(10.29.3)]{DLMF} and applying
\eqref{eq:airy-diagonal-kernel}, we obtain
\begin{align*}
    \|v_\mu(0)\|^2
    &=\frac\mu2\left\|\mathcal R_\mu(0)\binom10\right\|^2
      =O(\mu),\\
    L_{-\mu}(-\mu,-\mu)
    &=\frac{\mu^2}{4}
    -\frac12\left(\mathcal R_\mu(0)^{-1}\mathcal R_\mu'(0)\right)_{21}\\
    &=\frac{\mu^2}{4}+O(\mu^{-1}).
\end{align*}
This completes the proof.
\end{proof}

\subsection{Asymptotics for truncated orthogonal polynomials}
\label{subsec:right-rhp-comparison}

We now analyze the truncated orthogonal-polynomial RH problem by the
Deift--Zhou steepest descent method~\cite{DZ93}, using the Airy resolvent
model from Section~\ref{subsec:right-airy-model} and the construction of
Rider and Zhou~\cite[Sections~2--4]{RZ}.

\subsubsection{Truncated orthogonal-polynomial RH problem}
\label{subsec:right-janossy-input}

The truncated orthogonal-polynomial RH problem gives a representation of
the Janossy kernel, and hence, by \eqref{eq:finite-endpoint-identity}, of the
logarithmic derivative of the Fredholm determinant.
For $y\in\R$, let
$Y:=Y(y;N)$ solve the orthogonal-polynomial RH problem of
Fokas, Its and Kitaev~\cite{FIK92} for the weight
$e^{-NV(x)}\mathbf1_{(-\infty,y]}(x)$.

\begin{RHProblem}\label{rhp:Y}
\hfill\break
\begin{enumerate}
 \item[(a)] $Y$ is defined and analytic in $\mathbb C\setminus(-\infty,y]$.
 \item[(b)] For $x<y$, we have
 \[
  Y_+(x)=Y_-(x)\begin{pmatrix}1&e^{-NV(x)}\\0&1\end{pmatrix}.
 \]
 \item[(c)] As $z\to\infty$ with $z\in \mathbb C\setminus(-\infty,y]$,
 \[
  Y(z)=\left(I+O(z^{-1})\right)z^{N\sigma_3}.
 \]
 \item[(d)] As $z\to y$, we have
 \[
  Y(z)=\begin{pmatrix}
   O(1)&O(\log(z-y))\\
   O(1)&O(\log(z-y))
  \end{pmatrix}.
 \]
\end{enumerate}
\end{RHProblem}

Let $\pi_j(\cdot;y)$ be the monic degree-$j$ orthogonal polynomial
for this weight, with squared norm
\begin{equation}
 h_j(y):=\int_{-\infty}^y\pi_j(x;y)^2e^{-NV(x)}\,\dd x>0.
\label{eq:truncated-op-norm}
\end{equation}
The unique solution of RH problem~\ref{rhp:Y} is
\begin{equation}
 Y(z)=
 \begin{pmatrix}
  \pi_N(z;y)&
  \displaystyle\frac1{2\pi i}\int_{-\infty}^y
  \frac{\pi_N(x;y)e^{-NV(x)}}{x-z}\,\dd x\\[2ex]
  \displaystyle-\frac{2\pi i}{h_{N-1}(y)}\pi_{N-1}(z;y)&
  \displaystyle-\frac1{h_{N-1}(y)}\int_{-\infty}^y
  \frac{\pi_{N-1}(x;y)e^{-NV(x)}}{x-z}\,\dd x
 \end{pmatrix},
 \qquad z\in\mathbb C\setminus(-\infty,y],
 \label{eq:truncated-op-solution}
\end{equation}
see, for example, \cite[Section~2]{RZ}.
The truncated orthogonal-polynomial representation of the Janossy kernel
\cite[Section~1]{RZ} gives
\[
 L_{N,y}(u,v)
 =e^{-N(V(u)+V(v))/2}
 \sum_{j=0}^{N-1}\frac{\pi_j(u;y)\pi_j(v;y)}{h_j(y)},
 \qquad u,v>y.
\]

\subsubsection{The constrained equilibrium measure}
\label{subsec:right-moving-wall}

Recall the definition of $\varepsilon_N$ in \eqref{eq:edge-scale} and set
\begin{equation}
 \delta:=\mu\varepsilon_N,\qquad y:=B-\delta,\qquad 0\leq\delta\ll1.
\label{eq:moving-wall-variables}
\end{equation}
Let $\rho_\delta$ be the equilibrium density constrained to $(-\infty,y]$,
with support $[A_\delta,y]$ and Euler--Lagrange constant $\ell_\delta$.
By \cite[Lemma~2.2]{RZ}, we have
\begin{equation}
 \rho_\delta(x)=\frac1{2\pi}
 \sqrt{\frac{x-A_\delta}{y-x}}
 \bigl[(y-x)H_\delta(x)+C_\delta\bigr],
 \qquad A_\delta<x<y.
\label{eq:constrained-equilibrium-density}
\end{equation}
Here $H_\delta$ is analytic on a fixed neighborhood of $[A,B]$, and
$A_\delta,H_\delta,C_\delta$ depend analytically on $\delta$, with
\[
 H_0=2\pi h,\qquad C_\delta=\pi h(B)\delta+O(\delta^2).
\]
For $z\in\mathbb C\setminus(-\infty,y]$, define
\[
 g_\delta(z):=\int_{A_\delta}^y\log(z-x)\rho_\delta(x)\,\dd x,
 \qquad
 \Phi_\delta(z):=V(z)+\ell_\delta-2g_\delta(z).
\]

\subsubsection{First transformation: \texorpdfstring{$Y\to T_N$}{Y to T}}

This transformation normalizes the RH problem for $Y$ at infinity and is
defined by
\begin{equation}\label{eq:finite-Y-to-T}
 T_N(z):=e^{-N\ell_\delta\sigma_3/2}Y(z)
 e^{-N(g_\delta(z)-\ell_\delta/2)\sigma_3}.
\end{equation}
In view of the relations
\begin{equation}
\begin{aligned}
 g_{\delta,+}(x)+g_{\delta,-}(x)&=V(x)+\ell_\delta,
 &&A_\delta<x<y,\\
 \Phi_{\delta,+}(x)&=g_{\delta,-}(x)-g_{\delta,+}(x)=-\Phi_{\delta,-}(x),
 &&A_\delta<x<y,\\
 g_{\delta,+}(x)-g_{\delta,-}(x)&=2\pi i,
 &&x<A_\delta,\\
 \Re\Phi_{\delta,+}(x)&=V(x)+\ell_\delta-g_{\delta,+}(x)-g_{\delta,-}(x)>0,
 &&x<A_\delta,
\end{aligned}
\label{eq:finite-phase-boundary-relations}
\end{equation}
and  $g_\delta(z)=\log z+O(z^{-1})$ as
$z\to\infty$, it follows from RH problem~\ref{rhp:Y} that $T_N$ satisfies
the following RH problem.

\begin{RHProblem}\label{rhp:T}
\hfill\break
\begin{enumerate}
 \item[(a)] $T_N$ is defined and analytic in $\mathbb C\setminus(-\infty,y]$.
 \item[(b)] For $x<y$, we have
 \[
  T_{N,+}(x)=T_{N,-}(x)
  \begin{cases}
   \begin{pmatrix}e^{N\Phi_{\delta,+}(x)}&1\\0&e^{N\Phi_{\delta,-}(x)}\end{pmatrix},
       &A_\delta<x<y,\\[3mm]
   I+e^{-N\Phi_{\delta,+}(x)}E_{12},&x<A_\delta.
  \end{cases}
 \]
 \item[(c)] As $z\to\infty$ with $z\in\mathbb C\setminus(-\infty,y]$, we have
 \[
  T_N(z)=I+O(z^{-1}).
 \]
\end{enumerate}
\end{RHProblem}

\subsubsection{Second transformation: \texorpdfstring{$T_N\to S_N$}{T to S}}

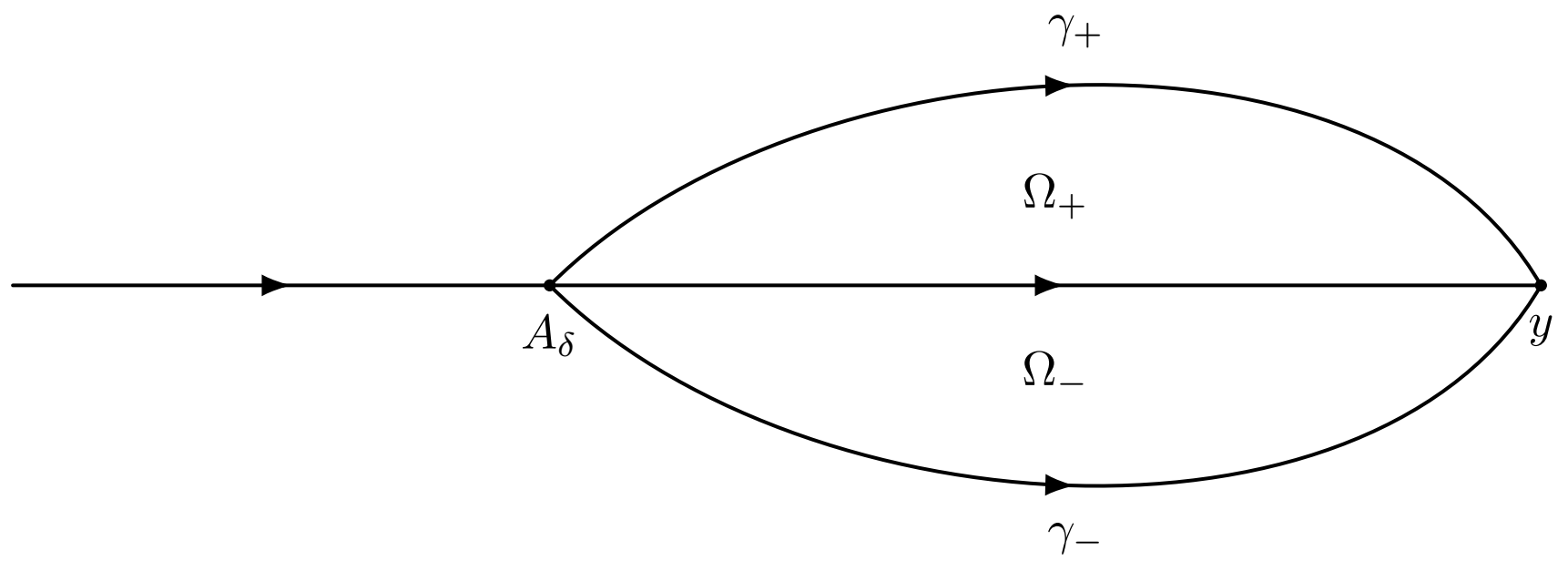
\begin{figure}[htbp]
\centering
\begin{tikzpicture}[rhp diagram,x=1.5cm,y=1.5cm]
 \draw[rhp contour] (-5,0)--(-2.4,0);
 \draw[rhp contour] (-2.4,0)--(2.4,0);
 \draw[rhp contour] (-2.4,0) .. controls (-1.2,1.2) and (1.6,1.38564) .. (2.4,0);
 \draw[rhp contour] (-2.4,0) .. controls (-1.2,-1.2) and (1.6,-1.38564) .. (2.4,0);
 \fill (-2.4,0) circle (1.25pt);
 \fill (2.4,0) circle (1.25pt);
 \node[below=4pt] at (-2.4,0) {$A_\delta$};
 \node[below=4pt] at (2.4,0) {$y$};
 \node[above=3pt] at (0.15,1.03) {$\gamma_+$};
 \node[below=3pt] at (0.15,-1.03) {$\gamma_-$};
 \node at (0.05,0.43) {$\Omega_+$};
 \node at (0.05,-0.43) {$\Omega_-$};
\end{tikzpicture}
\caption{Regions $\Omega_\pm$ and the jump contours in the RH problem for $S_N$.}
\label{fig:S-lens-contour}
\end{figure}

By the phase relations in \eqref{eq:finite-phase-boundary-relations},
the jump matrix in RH problem~\ref{rhp:T}
tends to $I$ exponentially fast as $N\to\infty$ on $(-\infty,A_\delta)$, whereas its diagonal entries on $(A_\delta,y)$
are highly oscillatory. We therefore open a lens around this interval
to replace the oscillatory jump by a constant matrix, at the cost of
introducing additional jumps that tend exponentially to $I$ outside the
endpoint disks. To proceed, we use $\Phi_{\delta,-}=-\Phi_{\delta,+}$ to
factor the jump matrix as
\begin{equation*}
 \begin{pmatrix}e^{N\Phi_{\delta,+}(x)}&1\\0&e^{N\Phi_{\delta,-}(x)}\end{pmatrix}=
 \begin{pmatrix}1&0\\e^{N\Phi_{\delta,-}(x)}&1\end{pmatrix}
 \begin{pmatrix}0&1\\-1&0\end{pmatrix}
 \begin{pmatrix}1&0\\e^{N\Phi_{\delta,+}(x)}&1\end{pmatrix}.
\end{equation*}
Let $\Omega_\pm$ denote the regions shown in
Figure~\ref{fig:S-lens-contour}, and define
\begin{equation}\label{eq:finite-T-to-S}
 S_N(z):=T_N(z)
 \begin{cases}
  I-e^{N\Phi_\delta(z)}E_{21},
       &z\in\Omega_+,\\[1mm]
  I+e^{N\Phi_\delta(z)}E_{21},
       &z\in\Omega_-,\\[1mm]
  I,&\text{elsewhere}.
 \end{cases}
\end{equation}
It is then readily seen from RH problem~\ref{rhp:T} for $T_N$ that $S_N$
satisfies the following RH problem.

\begin{RHProblem}\label{rhp:S}
\hfill\break
\begin{enumerate}
    \item[(a)] $S_N$ is defined and analytic in $\mathbb C\setminus\Sigma_S$,
    where $\Sigma_S:=(-\infty,y]\cup\gamma_+\cup\gamma_-$ is the contour
    shown in Figure~\ref{fig:S-lens-contour}.
    \item[(b)] For $z\in\Sigma_S$, we have
    \[
        S_{N,+}(z)=S_{N,-}(z)J_S(z),
    \]
    where
    \[
        J_S(z):=
        \begin{cases}
        \begin{pmatrix}0&1\\-1&0\end{pmatrix},&A_\delta<z<y,\\[3mm]
        I+e^{N\Phi_\delta(z)}E_{21},&z\in\gamma_+\cup\gamma_-,\\[1mm]
        I+e^{-N\Phi_{\delta,+}(z)}E_{12},&z<A_\delta.
        \end{cases}
    \]
    \item[(c)] As $z\to\infty$ with $z\in\mathbb C\setminus\Sigma_S$, we have
    \[
        S_N(z)=I+O(z^{-1}).
    \]
\end{enumerate}
\end{RHProblem}

\subsubsection{Global parametrix}

As $N\to\infty$, all the jump matrices $J_S(z)$ of $S_N$ given in RH
problem~\ref{rhp:S} tend to the identity matrix exponentially fast for $z$
bounded away from $[A_\delta,y]$.  This leads to the following global parametrix.

\begin{RHProblem}\label{rhp:global-parametrix}
\hfill\break
\begin{enumerate}
 \item[(a)] $P_N^{(\infty)}$ is defined and analytic in
 $\mathbb C\setminus[A_\delta,y]$.
 \item[(b)] For $x\in(A_\delta,y)$, we have
 \[
  P_{N,+}^{(\infty)}(x)=P_{N,-}^{(\infty)}(x)
  \begin{pmatrix}0&1\\-1&0\end{pmatrix}.
 \]
 \item[(c)] As $z\to\infty$ with $z\in\mathbb C\setminus[A_\delta,y]$, we have
 \[
  P_N^{(\infty)}(z)=I+O(z^{-1}).
 \]
\end{enumerate}
\end{RHProblem}

The solution is
\[
    P_N^{(\infty)}(z)=\frac12
    \begin{pmatrix}
      \beta_\delta+\beta_\delta^{-1}&i(\beta_\delta-\beta_\delta^{-1})\\
      -i(\beta_\delta-\beta_\delta^{-1})&\beta_\delta+\beta_\delta^{-1}
    \end{pmatrix},
\]
where $\beta_\delta(z)=\left(\frac{z-A_\delta}{z-y}\right)^{1/4}$.

\subsubsection{Local parametrix near \texorpdfstring{$y$}{y}}
\label{subsubsec:wall-parametrix}

In a sufficiently small fixed neighborhood $D(y,\epsilon)$ of $y$, we
intend to solve the following RH problem, which serves as an approximation
of $S_N$.

\begin{RHProblem}\label{rhp:wall-parametrix}
\hfill\break
\begin{enumerate}
 \item[(a)] $P_N^{(y)}$ is defined and analytic in $D(y,\epsilon)\setminus\Sigma_S$.
 \item[(b)] For $z\in D(y,\epsilon)\cap\Sigma_S$, we have
 \[
  P_{N,+}^{(y)}(z)=P_{N,-}^{(y)}(z)J_S(z).
 \]
 \item[(c)] As $N\to\infty$, $P_N^{(y)}(z)$ satisfies the matching condition
  \begin{equation}
  P_N^{(y)}(z)(P_N^{(\infty)}(z))^{-1}
  =I+O(N^{-1/3}),
  \label{eq:upper-matching-norms}
  \end{equation}
 uniformly for $z\in\partial D(y,\epsilon)$ and $\delta=O(a_N\varepsilon_N)$,
 where $\varepsilon_N$ is defined in \eqref{eq:edge-scale}.
\end{enumerate}
\end{RHProblem}

To use the model $M$ in RH problem~\ref{rhp:Mminus}, we match its phase
to $N\Phi_\delta$, as in \cite[Lemma~4.1]{RZ}.

\begin{Lemma}\label{lem:local-change-of-variables}
For all sufficiently small $\delta\geq0$, there exist $\mu_N=\mu_N(\delta)$
and a conformal map $\zeta_{N,\mu}$ on the disk $D(y,\epsilon)$,
where $y=B-\delta$, such that
\begin{equation}
 N\Phi_\delta(z)=\theta(\zeta_{N,\mu}(z);-\mu_N)
 =\frac43\zeta_{N,\mu}(z)^{3/2}-2\mu_N\zeta_{N,\mu}(z)^{1/2},
\label{eq:exact-normal-form}
\end{equation}
for $z\in D(y,\epsilon)\setminus(-\infty,y]$. Here
\begin{align}
 \zeta_{N,\mu}(z)&:=N^{2/3}(z-y)d_\delta(z),
\label{eq:moving-wall-coordinate}\\
 \mu_N(\delta)&:=N^{2/3}\delta_{\mathrm{eff}}(\delta)
 =\mu+O(\mu^2N^{-2/3}).
\label{eq:muhat-expansion}
\end{align}
$d_\delta(z)$ is jointly analytic in $(\delta,z)$ near $(0,B)$,
with $d_0(B)=\kappa$, where $\kappa$ is defined in
\eqref{eq:kappa-definition}.
The function $\delta_{\mathrm{eff}}$ is analytic near $\delta=0$ and
satisfies $\delta_{\mathrm{eff}}=\kappa\delta+O(\delta^2)$.
\end{Lemma}

\begin{proof}
Writing $z=y+t^2$ with $t>0$, we obtain from \eqref{eq:constrained-equilibrium-density}
\begin{equation}
 \frac{d}{dt}\Phi_\delta(y+t^2)
 =2\sqrt{y+t^2-A_\delta}
 \left[t^2H_\delta(y+t^2)-C_\delta\right].
\label{eq:hard-phase}
\end{equation}
We first seek $d_\delta$ and $\delta_{\mathrm{eff}}$ satisfying
\begin{equation}
 \frac{\Phi_\delta(y+t^2)}{t}
 =\frac43t^2d_\delta(y+t^2)^{3/2}
 -2\delta_{\mathrm{eff}}d_\delta(y+t^2)^{1/2}.
\label{eq:local-coordinate-equation}
\end{equation}
We solve this equation for the even function $d_\delta(y+t^2)$ and
$\delta_{\mathrm{eff}}$ near the positive analytic solution $d_0$
at $\delta=\delta_{\mathrm{eff}}=0$, with $d_0(B)=\kappa$.
Linearizing at this solution and dividing by $2d_0(B+t^2)^{1/2}$ on a fixed disk,
we determine the variation of $\delta_{\mathrm{eff}}$ from the constant term
and that of $d_\delta(y+t^2)$ by dividing the remainder by $t^2$.
Applying the analytic implicit function theorem to even analytic functions
continuous on the closed disk gives $d_\delta(y+t^2)$ and
$\delta_{\mathrm{eff}}$ satisfying \eqref{eq:local-coordinate-equation},
with $\delta_{\mathrm{eff}}$ analytic in $\delta$
and, by evenness, $d_\delta(z)$ jointly analytic in $(\delta,z)$.
Defining $\zeta_{N,\mu}$ by \eqref{eq:moving-wall-coordinate} and setting
$\mu_N=N^{2/3}\delta_{\mathrm{eff}}$, we obtain \eqref{eq:exact-normal-form}
from \eqref{eq:local-coordinate-equation}.
Evaluating \eqref{eq:local-coordinate-equation} at $t=0$ and using \eqref{eq:hard-phase}, we obtain
$\delta_{\mathrm{eff}}=\kappa\delta+O(\delta^2)$;
substituting $\delta=\mu\varepsilon_N$ then gives \eqref{eq:muhat-expansion}.
This completes the proof.
\end{proof}

We can therefore define the local parametrix by
\begin{equation*}
 P_N^{(y)}(z):=E_N^{(y)}(z)
 M(\zeta_{N,\mu}(z);-\mu_N),
\end{equation*}
where
\[
 E_N^{(y)}(z):=P_N^{(\infty)}(z)\frac1{\sqrt{2i}}\begin{pmatrix}i&-i\\1&1\end{pmatrix}
 \zeta_{N,\mu}(z)^{\sigma_3/4}
 .
\]

\begin{Proposition}\label{prop:wall-parametrix}
The matrix $P_N^{(y)}$ solves
RH problem~\ref{rhp:wall-parametrix}.
\end{Proposition}
\begin{proof}
We have
\begin{equation}
 E_N^{(y)}(z)=\frac1{\sqrt{2i}}
 \begin{pmatrix}i&-i\\1&1\end{pmatrix}
 \left[N^{2/3}(z-A_\delta)d_\delta(z)\right]^{\sigma_3/4},
\label{eq:E-factorization}
\end{equation}
which shows that $E_N^{(y)}$ is analytic and invertible
in $D(y,\epsilon)$.
The jump relation in RH problem~\ref{rhp:Mminus} and the phase identity
\eqref{eq:exact-normal-form} then give the required jumps for $P_N^{(y)}$.
On $\partial D(y,\epsilon)$, Lemma~\ref{lem:local-change-of-variables}
and \eqref{eq:airy-G-bounds} give
\[
 G_{\mu_N}(\zeta_{N,\mu})-I=O(N^{-2/3}),
\]
uniformly for $\delta=O(a_N\varepsilon_N)$.
Together with \eqref{eq:E-factorization}, this yields
\[
 \begin{aligned}
 P_N^{(y)}(P_N^{(\infty)})^{-1}-I
 &=E_N^{(y)}\bigl(G_{\mu_N}(\zeta_{N,\mu})-I\bigr)(E_N^{(y)})^{-1}\\
 &=O(N^{1/6}N^{-2/3}N^{1/6})=O(N^{-1/3}),
 \end{aligned}
\]
which proves \eqref{eq:upper-matching-norms}.
\end{proof}

\subsubsection{Local parametrix near \texorpdfstring{$A_\delta$}{A-delta}}
\label{subsubsec:left-parametrix}

\begin{RHProblem}\label{rhp:left-parametrix}
\hfill\break
\begin{enumerate}
 \item[(a)] $P_N^{(A_\delta)}$ is defined and analytic in $D(A_\delta,\epsilon)\setminus\Sigma_S$.
 \item[(b)] For $z\in D(A_\delta,\epsilon)\cap\Sigma_S$, we have
 \[
  P_{N,+}^{(A_\delta)}(z)=P_{N,-}^{(A_\delta)}(z)J_S(z).
 \]
 \item[(c)] As $N\to\infty$, $P_N^{(A_\delta)}(z)$ satisfies the matching condition
 \begin{align}
  P_N^{(A_\delta)}(z)(P_N^{(\infty)}(z))^{-1}
  =I+O(N^{-1}),
  \label{eq:left-airy-matching}
 \end{align}
 uniformly for $z\in\partial D(A_\delta,\epsilon)$ and all sufficiently small $\delta\geq0$.
\end{enumerate}
\end{RHProblem}

The above RH problem can be solved using the Airy parametrix
$\Phi^{(\mathrm{Ai})}$ defined in Appendix~\ref{airy}. To this end, we introduce
\begin{align*}
 \phi_\delta(z)&:=\frac12\int_{A_\delta}^z \Phi_\delta'(s)\,\dd s,
 &f_\delta(z)&:=\left(\frac32\phi_\delta(z)\right)^{2/3},
\end{align*}
where the $2/3$ power is chosen so that $f_\delta$ is analytic in
$D(A_\delta,\epsilon)$ and positive for real $z<A_\delta$ in this disk.
The phases satisfy
\begin{equation}
 \Phi_\delta=2\phi_\delta\mp2\pi i
\label{eq:left-airy-phase-relation}
\end{equation}
in the upper and lower half-disks, respectively.
We then define
\begin{equation}
 P_N^{(A_\delta)}(z):=E_N^{(A)}(z)\sigma_3
 \Phi^{(\mathrm{Ai})}(N^{2/3}f_\delta(z))
 \sigma_3e^{N\phi_\delta(z)\sigma_3},
\label{eq:left-airy-parametrix}
\end{equation}
where
\[
 E_N^{(A)}(z):=P_N^{(\infty)}(z)
 \frac1{\sqrt2}\begin{pmatrix}1&i\\i&1\end{pmatrix}
 \left(N^{2/3}f_\delta(z)\right)^{\sigma_3/4}.
\]

\begin{Proposition}
\label{prop:left-airy-parametrix}
The matrix $P_N^{(A_\delta)}$ solves
RH problem~\ref{rhp:left-parametrix}.
\end{Proposition}
\begin{proof}
We have
\[
 E_N^{(A)}(z)=\frac1{\sqrt2}\begin{pmatrix}1&i\\i&1\end{pmatrix}
 \left(N^{2/3}(z-y)\frac{f_\delta(z)}{z-A_\delta}\right)^{\sigma_3/4},
\]
which shows that $E_N^{(A)}$ is analytic and invertible in
$D(A_\delta,\epsilon)$.
The jump relation follows from 
\eqref{eq:left-airy-phase-relation} together with the jump condition \eqref{eq:airy-jump-matrices}  of the Airy parametrix.
The matching condition follows from the Airy asymptotics \eqref{infty:Ai}.
\end{proof}

\subsubsection{Final transformation}
\label{subsubsec:final-transformation}

Define
\begin{equation}\label{eq:finite-rh-error-definition}
    R_N(z):=
    \begin{cases}
      S_N(z)P_N^{(y)}(z)^{-1},&z\in D(y,\epsilon),\\
      S_N(z)P_N^{(A_\delta)}(z)^{-1},&z\in D(A_\delta,\epsilon),\\
      S_N(z)(P_N^{(\infty)}(z))^{-1},
      &\text{elsewhere}.
    \end{cases}
\end{equation}

\begin{figure}[htbp]
\centering
\begin{tikzpicture}[rhp diagram,x=1.5cm,y=1.5cm]
 \draw[rhp contour] (-5,0)--(-3.05,0);
 \begin{scope}[even odd rule]
  \clip (-3.1,-1.1) rectangle (3.1,1.1)
   (-2.4,0) circle (0.65) (2.4,0) circle (0.65);
  \draw[rhp contour] (-2.4,0) .. controls (-1.2,1.2) and (1.6,1.38564) .. (2.4,0);
  \draw[rhp contour] (-2.4,0) .. controls (-1.2,-1.2) and (1.6,-1.38564) .. (2.4,0);
 \end{scope}
 \draw[rhp circle] (-1.75,0) arc[start angle=0,end angle=-360,radius=0.65];
 \draw[rhp circle] (3.05,0) arc[start angle=0,end angle=-360,radius=0.65];
 \fill (-2.4,0) circle (1.25pt);
 \fill (2.4,0) circle (1.25pt);
 \node[below=4pt] at (-2.4,0) {$A_\delta$};
 \node[below=4pt] at (2.4,0) {$y$};
 \node[above] at (-2.4,0.76) {$\partial D(A_\delta,\epsilon)$};
 \node[above] at (2.4,0.76) {$\partial D(y,\epsilon)$};
 \node[above=3pt] at (0.15,1.03) {$\gamma_+$};
 \node[below=3pt] at (0.15,-1.03) {$\gamma_-$};
\end{tikzpicture}
\caption{Contour $\Sigma_R$ in the RH problem for $R_N$.}
\label{fig:right-rhp-lens-error}
\end{figure}
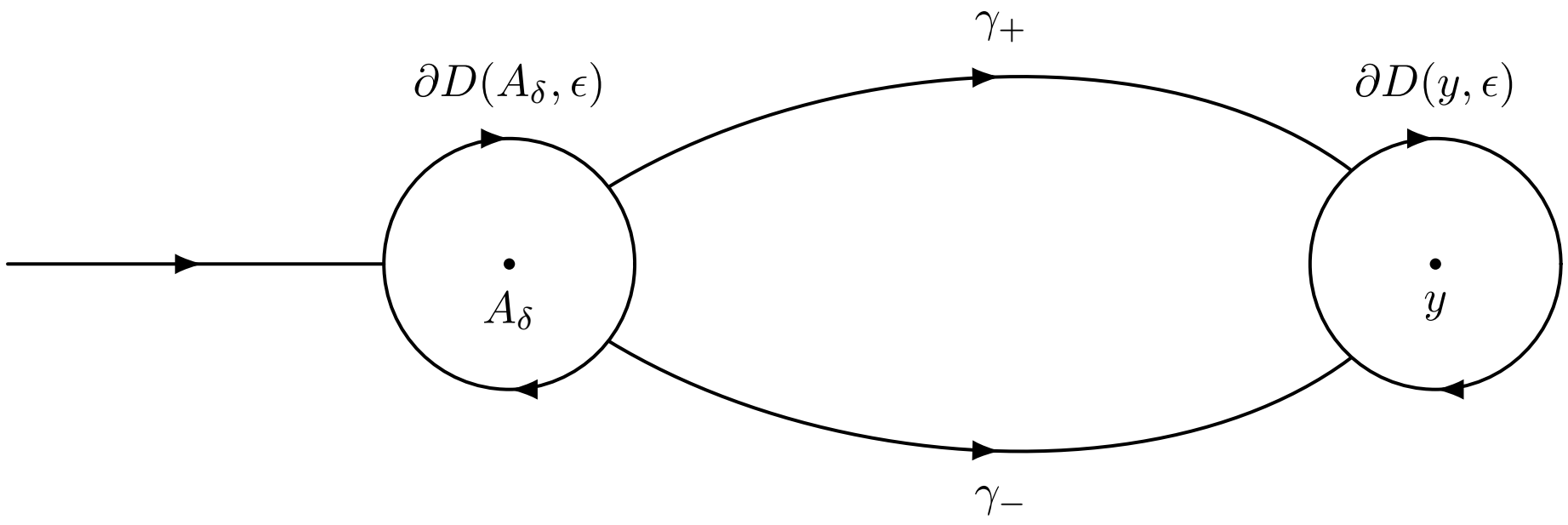

\pagebreak[3]
\begin{RHProblem}\label{rhp:R}
\hfill\break
\begin{enumerate}
    \item[(a)] $R_N$ is defined and analytic in $\mathbb C\setminus\Sigma_R$, where 
    \[
\Sigma_R:=\partial D(y,\epsilon)\cup\partial D(A_\delta,\epsilon)
\cup\bigl(\Sigma_S\setminus(D(y,\epsilon)\cup D(A_\delta,\epsilon)\cup[A_\delta,y])\bigr),
\]
as shown in Figure~\ref{fig:right-rhp-lens-error}.
    \item[(b)] For $z\in\Sigma_R$, we have
    \[
        R_{N,+}(z)=R_{N,-}(z)J_R(z),
    \]
    where
    \[
        J_R(z):=
        \begin{cases}
        P_N^{(y)}(z)(P_N^{(\infty)}(z))^{-1},&z\in\partial D(y,\epsilon),\\
        P_N^{(A_\delta)}(z)(P_N^{(\infty)}(z))^{-1},&z\in\partial D(A_\delta,\epsilon),\\
        P_{N,-}^{(\infty)}(z)J_S(z)
        (P_{N,+}^{(\infty)}(z))^{-1},
        &z\in\Sigma_R\setminus(\partial D(y,\epsilon)\cup\partial D(A_\delta,\epsilon)).
        \end{cases}
    \]
    where $J_S$ is defined in RH problem~\ref{rhp:S}.
    \item[(c)] As $z\to\infty$ with $z\in\mathbb C\setminus\Sigma_R$, we have
    \[
        R_N(z)=I+O(z^{-1}).
    \]
\end{enumerate}
\end{RHProblem}

The matching estimates \eqref{eq:upper-matching-norms} and
\eqref{eq:left-airy-matching}, together with the exponentially small
jumps on the rest of $\Sigma_R$, yield by a small-norm argument
\cite{Deift1999,DZ93}
\begin{equation}
 R_N(z)=I+O(N^{-1/3}),
\label{eq:R-bounds}
\end{equation}
uniformly for $z\in\mathbb C\setminus\Sigma_R$ and
$\delta=O(a_N\varepsilon_N)$.

\subsubsection{Representation of the resolvent kernel}
\label{subsec:right-recovery}

The Christoffel--Darboux formula and \eqref{eq:truncated-op-solution} give
\begin{equation}
 L_{N,y}(u,v)
 =-\frac{e^{-N(V(u)+V(v))/2}}{2\pi i}
 \frac{Y_{11}(u)Y_{21}(v)-Y_{21}(u)Y_{11}(v)}{u-v},
 \qquad u,v>y.
\label{eq:L-Y-recovery}
\end{equation}
Letting $v\to u$ gives
\begin{equation}
 L_{N,y}(u,u)=\frac{e^{-NV(u)}}{2\pi i}
  \left(Y_{11}(u)Y_{21}'(u)-Y_{21}(u)Y_{11}'(u)\right).
\label{eq:L-Y-diagonal}
\end{equation}

Tracing back the transformations $Y\to T_N\to S_N\to R_N$ in
\eqref{eq:finite-Y-to-T}, \eqref{eq:finite-T-to-S} and
\eqref{eq:finite-rh-error-definition}, we obtain from
\eqref{eq:L-Y-diagonal}
\begin{align}
 \varepsilon_NL_{N,y}(z,z)
 ={}&\varepsilon_N(\partial_z\zeta_{N,\mu})
 L_{-\mu_N}(\zeta_{N,\mu}-\mu_N,
 \zeta_{N,\mu}-\mu_N)
\notag\\
 &+\varepsilon_N\det\left((R_NE_N^{(y)})^{-1}
 \partial_z(R_NE_N^{(y)})v_{\mu_N}(\zeta_{N,\mu}),
 v_{\mu_N}(\zeta_{N,\mu})\right),
\label{eq:full-recovery-exact}
\end{align}
for $\delta=O(a_N\varepsilon_N)$ and $y<z<y+\epsilon$.
By Cauchy's formula applied to \eqref{eq:R-bounds} and the factorization
\eqref{eq:E-factorization}, the logarithmic derivative in
\eqref{eq:full-recovery-exact} satisfies, for $z$ near $y$,
\begin{equation}
 \varepsilon_N(R_NE_N^{(y)})^{-1}\partial_z(R_NE_N^{(y)})
 =\begin{pmatrix}
 O(N^{-2/3})&O(N^{-4/3})\\
 O(N^{-2/3})&O(N^{-2/3})
 \end{pmatrix}.
\label{eq:K-prefactor}
\end{equation}

\newpage
\section{Right-Tail Asymptotics and Conditional Limits}
\label{sec:right-proof}

\subsection{Proof of \eqref{eq:fixed-k-airy-right}}
\label{subsec:right-airy-tail}

We first derive asymptotic formulas for the determinant in the Airy Janossy
density and its integral.

\begin{Lemma}
\label{lem:airy-right-determinant}
Fix an integer $k\geq1$ and $K\Subset(0,\infty)$. Set $\mu=a^{1/4}u$
and write $\mathbf x=(x_1,\ldots,x_k)$ with $x_1>\cdots>x_k=a$.
As $a\to\infty$, for every fixed $L>0$,
\begin{equation}
\begin{aligned}
\det[L_{-\mu}(x_i-\mu,x_j-\mu)]_{i,j=1}^k
={}&2^{-k(2k+1)}\pi^{-k}a^{-k^2}\Delta(\mathbf x)^2\\
&\quad\times e^{-k\theta(a;-\mu)-2\sqrt a\sum_{i=1}^{k-1}(x_i-a)}
\left[1+o(1)\right],
\end{aligned}
\label{eq:fk-3-4}
\end{equation}
uniformly for $u\in K$ and $x_1\leq a+La^{-1/2}\log a$.
Moreover,
\begin{equation}
\begin{aligned}
&\int_{x_1>\cdots>x_{k-1}>a}
\det[L_{-\mu}(x_i-\mu,x_j-\mu)]_{i,j=1}^k
\prod_{i=1}^{k-1}\mathrm dx_i\\
&\qquad=c_ka^{-(3k^2-1)/2}e^{-k\theta(a;-\mu)}
\left[1+o(1)\right],
\end{aligned}
\label{eq:fk-3-14}
\end{equation}
uniformly for $u\in K$, where
\begin{equation}
c_k:=2^{-(3k^2+k-1)}\pi^{-k}
\prod_{j=0}^{k-2}j!(j+2)!.
\label{eq:fk-3-15}
\end{equation}
\end{Lemma}

\begin{proof}
For $v_\mu=(v_{\mu,1},v_{\mu,2})^{\mathsf T}$ defined in
\eqref{eq:airy-recovery-vector}, set $r_\mu:=v_{\mu,2}/v_{\mu,1}$.
Equations~\eqref{eq:airy-leading-vector} and \eqref{eq:airy-G-bounds} give
\begin{align}
 v_{\mu,1}(x)
 &=\frac{e^{-\theta(x;-\mu)/2}}{2\sqrt\pi\,x^{1/4}}[1+o(1)],
 \label{eq:fk-2-6}\\
 \left|\partial_x^j\left(r_\mu(x)+\sqrt{x}\right)\right|
 &=o(a^{1/2-j}),\qquad1\leq j\leq2k-1,
 \label{eq:fk-2-5}
\end{align}
uniformly for $u\in K$ and $a\leq x\leq a+La^{-1/2}\log a$, for each fixed $L>0$.
Equation~\eqref{eq:fk-2-3} gives
\[
\frac{\det[L_{-\mu}(x_i-\mu,x_j-\mu)]_{i,j=1}^k}
{\prod_{i=1}^k v_{\mu,1}(x_i)^2}
=\det\left[\frac{r_\mu(x_j)-r_\mu(x_i)}{x_i-x_j}\right]_{i,j=1}^k.
\]
Taking successive row and column divided differences and applying the
Hermite--Genocchi formula \cite[Section~9, equation~(52)]{deBoor2005}
together with \eqref{eq:fk-2-5}, we obtain
\[
\begin{aligned}
 &\frac{a^{k(2k-1)/2}}{\Delta(\mathbf x)^2}
 \det\left[\frac{r_\mu(x_j)-r_\mu(x_i)}{x_i-x_j}\right]_{i,j=1}^k\\
 &\qquad=\det\left[-\frac{a^{i+j-3/2}r_\mu^{(i+j-1)}(a)}{(i+j-1)!}\right]_{i,j=1}^k+o(1)\\
 &\qquad=\det\left[\binom{1/2}{i+j-1}\right]_{i,j=1}^k+o(1),
\end{aligned}
\]
uniformly for $u\in K$ and $x_1\leq a+La^{-1/2}\log a$.
The binomial determinant formula \cite[Theorem~26]{Krattenthaler1999} gives
\begin{equation}
 \det\left[\binom{1/2}{i+j-1}\right]_{i,j=1}^k=2^{-2k^2+k}>0.
\label{eq:fk-3-9}
\end{equation}
Noticing that
\begin{equation}
\theta(x;-\mu)-\theta(a;-\mu)
=2\sqrt a(x-a)+o(1),\qquad a\leq x\leq a+La^{-1/2}\log a,
\label{eq:fk-3-10}
\end{equation}
we obtain \eqref{eq:fk-3-4} from \eqref{eq:fk-2-6} and \eqref{eq:fk-3-9}.

The Laguerre integral formula \cite[Equation~(5.14.5)]{DLMF} gives
\begin{equation}
\begin{aligned}
&\int_{x_1>\cdots>x_{k-1}>a}
\Delta(\mathbf x)^2
e^{-2\sqrt a\sum_{i=1}^{k-1}(x_i-a)}
\prod_{i=1}^{k-1}\mathrm dx_i\\
&\qquad=(2\sqrt a)^{-(k^2-1)}\prod_{j=0}^{k-2}j!(j+2)!.
\end{aligned}
\label{eq:fk-5-9}
\end{equation}
By \eqref{eq:model-B-uniform} and
$\partial_x\theta(x;-\mu)=2\sqrt x-\mu/\sqrt x\geq\sqrt a$ for $x\geq a$,
\[
 L_{-\mu}(x-\mu,x-\mu)
 \leq Ca^{-1}e^{-\theta(a;-\mu)-\sqrt a(x-a)},\qquad x\geq a.
\]
For $k\geq2$, fix $L>3k(k-1)/2$. The Gram representation
\eqref{eq:conditional-resolvent-gram} and Hadamard's inequality give
\[
\begin{aligned}
 0\leq\det[L_{-\mu}(x_i-\mu,x_j-\mu)]_{i,j=1}^k
 &\leq\prod_{i=1}^k L_{-\mu}(x_i-\mu,x_i-\mu)\\
 &\leq Ca^{-k}e^{-k\theta(a;-\mu)-\sqrt a\sum_{i=1}^{k-1}(x_i-a)}.
\end{aligned}
\]
Hence,
\begin{equation}
\begin{aligned}
&a^{(3k^2-1)/2}e^{k\theta(a;-\mu)}
\int_{\substack{x_1>\cdots>x_{k-1}>a\\
x_1>a+La^{-1/2}\log a}}
\det[L_{-\mu}(x_i-\mu,x_j-\mu)]_{i,j=1}^k
\prod_{i=1}^{k-1}\mathrm dx_i\\
&\qquad\leq Ca^{3k(k-1)/2-L}=o(1).
\end{aligned}
\label{eq:fk-3-6}
\end{equation}
Combining \eqref{eq:fk-3-4}, \eqref{eq:fk-5-9} and
\eqref{eq:fk-3-6} gives \eqref{eq:fk-3-14}.
\end{proof}

\begin{proof}[Proof of \eqref{eq:fixed-k-airy-right}]
We argue by induction on $k\geq1$, assuming the result for all
positive indices smaller than $k$, and the case $k=1$ is included in the
argument below. Writing the wall as $-\mu$ and the particle positions as
$x_i-\mu$, we set $x_k=a$ and rewrite \eqref{eq:fk-3-23} as
\begin{align}
J_{a,k}(\mu)
&:=F_2(-\mu)\int_{x_1>\cdots>x_{k-1}>a}
\det[L_{-\mu}(x_i-\mu,x_j-\mu)]_{i,j=1}^k
\prod_{i=1}^{k-1}\mathrm dx_i,\label{eq:fk-airy-janossy-integrand}\\
\mathbb P(\Delta_{\mathrm{Ai},k}>a)
&=\int_{\mathbb R}J_{a,k}(\mu)\,\mathrm d\mu.\notag
\end{align}

The Tracy--Widom left-tail asymptotic \cite[Section~1.2]{BBD}
(see also \cite[Theorem~1]{DIK}) is
\begin{equation}
 \log F_2(-t)=-\frac{t^3}{12}-\frac18\log t
 +\frac{\log2}{24}+\zeta'(-1)+O(t^{-3}),\qquad t\to\infty.
\label{eq:TW-left-full}
\end{equation}
Lemma~\ref{lem:airy-right-determinant} and \eqref{eq:TW-left-full} give,
uniformly for $u$ in compact subsets of $(0,\infty)$,
\begin{equation}
\begin{aligned}
J_{a,k}(a^{1/4}u)
={}&2^{1/24}e^{\zeta'(-1)}c_k
a^{-(3k^2-1)/2-1/32}u^{-1/8}\\
&\quad\times e^{-\frac{4k}{3}a^{3/2}+a^{3/4}f_k(u)}[1+o(1)],
\end{aligned}
\label{eq:fk-airy-local-density}
\end{equation}
where
\[
 f_k(u):=2ku-\frac{u^3}{12}.
\]

Let $u_k:=2\sqrt{2k}$ be the maximizer of $f_k$ on $(0,\infty)$, and fix
$0<\delta<1$.
Applying Hadamard's inequality to
$\det[K_{\mathrm{Ai}}(x_i,x_j)]_{i,j=1}^k$ \cite{TW94} and using the Airy asymptotic
\cite[Equation~(9.7.5)]{DLMF}, we obtain
\begin{equation}
\begin{aligned}
\int_{-\infty}^0J_{a,k}(\mu)\,\mathrm d\mu
&\leq\mathbb P(\zeta_k>a)\\
&\leq\frac1{k!}\int_{(a,\infty)^k}
\det[K_{\mathrm{Ai}}(x_i,x_j)]_{i,j=1}^k
\prod_{i=1}^k\mathrm dx_i\\
&\leq\frac1{k!}\left(\int_a^\infty K_{\mathrm{Ai}}(x,x)\,\mathrm dx\right)^k\\
&=\frac{(\operatorname{Tr}K_{\mathrm{Ai},a})^k}{k!}
\leq C a^{-3k/2}e^{-\frac{4k}{3}a^{3/2}}.
\end{aligned}
\label{eq:fk-3-19}
\end{equation}

For $0\leq\mu\leq a$, applying Hadamard's inequality to the determinant
defining $J_{a,k}(\mu)$ in \eqref{eq:fk-airy-janossy-integrand} gives
\begin{align*}
J_{a,k}(\mu)
&\leq\frac{F_2(-\mu)}{(k-1)!}\,L_{-\mu}(a-\mu,a-\mu)\\
&\quad\times\left(\int_a^\infty
L_{-\mu}(x-\mu,x-\mu)\,\mathrm dx\right)^{k-1}.
\end{align*}
Using \eqref{eq:model-B-uniform} and \eqref{eq:TW-left-full} in this bound yields
\[
 0\leq J_{a,k}(\mu)
 \leq a^{C}\exp\left[-\frac{4k}{3}a^{3/2}
 +a^{3/4}f_k\left(\frac{\mu}{a^{1/4}}\right)\right].
\]
Since
\[
 f_k(u_k)-f_k(u)=\frac{(u-u_k)^2(u+2u_k)}{12},
\]
we obtain
\begin{equation}
\int_{\substack{0\leq\mu\leq a\\
|\mu/(a^{1/4}u_k)-1|>\delta}}J_{a,k}(\mu)\,\mathrm d\mu
\leq a^{C}e^{-\frac{4k}{3}a^{3/2}
+a^{3/4}(f_k(u_k)-c)}.
\label{eq:fk-airy-off-center}
\end{equation}

For the remaining range $\mu\geq a$, we use the bound implied by
\eqref{eq:TW-left-full}:
\begin{equation}
 \log F_2(-t)\leq C-ct^3,\qquad t\geq0.
\label{eq:TW-left-global}
\end{equation}
Fix $C_0>1$. Equations~\eqref{eq:model-B-uniform} and
\eqref{eq:TW-left-global} give
\begin{align}
\int_a^{C_0a}J_{a,k}(\mu)\,\mathrm d\mu
&\leq a^{C}e^{Ca^{3/2}}
\int_a^{C_0a}e^{-c\mu^3}\,\mathrm d\mu\notag\\
&\leq e^{-ca^3}.
\label{eq:fk-airy-intermediate-mu}
\end{align}
The identity $\zeta_1-\zeta_k=\sum_{j=1}^{k-1}\Delta_{\mathrm{Ai},j}$ gives
\begin{equation}
\begin{aligned}
\int_{C_0a}^\infty J_{a,k}(\mu)\,\mathrm d\mu
&\leq\mathbb P(\zeta_k\leq-(C_0-1)a)\\
&\leq
\begin{cases}
F_2(-(C_0-1)a),&k=1,\\
F_2(-(C_0-1)a/2)+\displaystyle\sum_{j=1}^{k-1}
\mathbb P\left(\Delta_{\mathrm{Ai},j}\geq\dfrac{(C_0-1)a}{2(k-1)}\right),&k\geq2.
\end{cases}
\end{aligned}
\label{eq:fk-airy-far-left-decomposition}
\end{equation}
By \eqref{eq:TW-left-global} and the induction hypothesis, for sufficiently large $C_0$, we have
\begin{equation}
\int_{C_0a}^\infty J_{a,k}(\mu)\,\mathrm d\mu
\leq e^{-(4k/3+1)a^{3/2}}.
\label{eq:fk-airy-far-left-bound}
\end{equation}

Combining \eqref{eq:fk-3-19}, \eqref{eq:fk-airy-off-center},
\eqref{eq:fk-airy-intermediate-mu} and \eqref{eq:fk-airy-far-left-bound}
with \eqref{eq:fk-airy-local-density}, we obtain
\begin{equation}
\label{eq:fk-airy-localization}
\frac{\displaystyle\int_{|\mu/(a^{1/4}u_k)-1|>\delta}
J_{a,k}(\mu)\,\mathrm d\mu}
{\displaystyle\int_{|\mu/(a^{1/4}u_k)-1|\leq\delta}
J_{a,k}(\mu)\,\mathrm d\mu}
=O(e^{-ca^{3/4}}).
\end{equation}
Since
\[
 f_k(u_k)=\frac{8\sqrt2}{3}k^{3/2},\qquad -f_k''(u_k)=\sqrt{2k},
\]
Laplace's method applied to \eqref{eq:fk-airy-local-density},
using \eqref{eq:fk-airy-localization}, gives
\begin{align*}
\mathbb P(\Delta_{\mathrm{Ai},k}>a)
&=a^{1/4}J_{a,k}(a^{1/4}u_k)
\sqrt{\frac{2\pi}{a^{3/4}\sqrt{2k}}}\,[1+o(1)]\\
&=C_ka^{-(48k^2-11)/32}
e^{-\frac{4k}{3}a^{3/2}+\frac{8\sqrt2}{3}k^{3/2}a^{3/4}}[1+o(1)],
\end{align*}
where \eqref{eq:fk-3-15} gives
\[
 C_k=2^{1/24}e^{\zeta'(-1)}c_ku_k^{-1/8}
 \sqrt{\frac{2\pi}{\sqrt{2k}}}.
\]
Taking logarithms proves \eqref{eq:fixed-k-airy-right}--\eqref{eq:def-CRk}.
\end{proof}

\subsection{Proof of \eqref{eq:fixed-k-finite-right}}
\label{subsec:right-finite-tail}
\label{subsec:right-truncation-saddle}

We first obtain a tail bound for the untruncated kernel, compare the
finite-$N$ and Airy largest-eigenvalue distributions, and derive asymptotic
formulas for the determinant in the Janossy density.

\begin{Lemma}
\label{lem:ordinary-right-edge-count}
Let $S=S(N)$ satisfy $1\leq S=o(N^{2/3})$. For the correlation kernel $K_N$
and the phase $\Phi$ in \eqref{eq:exterior-phase},
\begin{align}
    \int_{B+S\varepsilon_N}^{\infty}K_N(x,x)\,\dd x
    &=O\left(e^{-N\Phi(B+S\varepsilon_N)}\right).
\label{eq:ordinary-right-count-exact}
\end{align}
\end{Lemma}

\begin{proof}
The untruncated kernel representation and Airy asymptotics in
\cite[proof of Theorem~2.2, (5.2)]{LMZ26}
(see also \cite{DeiftGioev2007}) give, for a sufficiently small fixed $r>0$,
\[
 K_N(x,x)\leq CN\sqrt{x-B}\,e^{-N\Phi(x)},
 \qquad B+\varepsilon_N\leq x\leq B+r.
\]
Since $\Phi'(x)\asymp\sqrt{x-B}$ on $(B,B+r]$, integration yields
\[
 \int_{B+S\varepsilon_N}^{B+r}K_N(x,x)\,\dd x
 =O\left(e^{-N\Phi(B+S\varepsilon_N)}\right).
\]

By \cite[proof of Theorem~2.2, (5.17)--(5.18)]{LMZ26}, the integral over
$[B+r,\infty)$ is $O(e^{-cN})$;
noting that $N\Phi(B+S\varepsilon_N)=o(N)$, we conclude the proof.
\end{proof}

\begin{Lemma}
\label{prop:finite-comparison}
We have
\begin{align}
 \log D_N(B-\mu\varepsilon_N)
 =\log F_2(-\mu)+O\left(N^{-2/3}(1+\mu^4)\right),
\label{eq:det-comparison}
\end{align}
uniformly for $0\leq\mu=O(a_N)$.
\end{Lemma}

\begin{proof}
Differentiating the Airy Fredholm determinant gives
\[
 \frac{\dd}{\dd s}\log F_2(s)=L_s(s,s),\qquad s\in\R.
\]
Since $y=B-\mu\varepsilon_N$, this identity and
\eqref{eq:finite-endpoint-identity} yield
\[
 \frac{\dd}{\dd\mu}\log\frac{D_N(y)}{F_2(-\mu_N)}
 =-\varepsilon_NL_{N,y}(y,y)
   +\frac{\dd\mu_N}{\dd\mu}L_{-\mu_N}(-\mu_N,-\mu_N).
\]
Applying \eqref{eq:K-prefactor} together with Lemmas~\ref{lem:wall-Wronskian}
and~\ref{lem:local-change-of-variables} to \eqref{eq:full-recovery-exact}
as $z\downarrow y$, we obtain
\begin{align*}
 \frac{\dd}{\dd\mu}\log\frac{D_N(y)}{F_2(-\mu_N)}
 &=\frac{\delta_{\mathrm{eff}}'(\delta)-d_\delta(y)}{\kappa}
 L_{-\mu_N}(-\mu_N,-\mu_N)
 +O(N^{-2/3}(1+\mu_N))\\
 &=O\left(N^{-2/3}(1+\mu^3)\right).
\end{align*}
Corollary~2.3 of \cite{LMZ26} at $t=0$ gives
\begin{equation}
 |\log D_N(B)-\log F_2(0)|=O(N^{-2/3}).
\label{eq:det-anchor}
\end{equation}
Integrating from $0$ to $\mu$, with $\mu_N(0)=0$ and
\eqref{eq:det-anchor}, we find
$\log[D_N(y)/F_2(-\mu_N)]=O(N^{-2/3}(1+\mu^4))$.
Finally, \eqref{eq:muhat-expansion} and Lemma~\ref{lem:wall-Wronskian} give
\[
 \log\frac{F_2(-\mu_N)}{F_2(-\mu)}
 =\int_{\mu_N}^{\mu}L_{-t}(-t,-t)\,\dd t
 =O\left(N^{-2/3}\mu^2(1+\mu^2)\right),
\]
which completes the comparison \eqref{eq:det-comparison}.
\end{proof}

For $0\leq\mu=o(N^{2/3})$, write $y=B-\mu\varepsilon_N$.
With $\zeta_{N,\mu}$ and $v_{\mu_N}$ as in
\eqref{eq:moving-wall-coordinate} and \eqref{eq:airy-recovery-vector}, define
$w_{N,\mu}=(w_{N,\mu,1},w_{N,\mu,2})^{\mathsf T}$ by
\begin{equation}
\begin{aligned}
w_{N,\mu}(\xi):={}&\bigl[R_N(y)E_N^{(y)}(y)\bigr]^{-1}
R_N(y+\varepsilon_N\xi)E_N^{(y)}(y+\varepsilon_N\xi)\\
&\quad\times v_{\mu_N}\bigl(\zeta_{N,\mu}(y+\varepsilon_N\xi)\bigr).
\end{aligned}
\label{eq:right-recovery-vector}
\end{equation}
Substituting the RH representation of $Y$ into \eqref{eq:L-Y-recovery}
and using \eqref{eq:exact-normal-form}, we obtain
\begin{equation}
\varepsilon_NL_{N,y}
\bigl(y+\varepsilon_N\xi_1,y+\varepsilon_N\xi_2\bigr)
=\frac{\det(w_{N,\mu}(\xi_1),w_{N,\mu}(\xi_2))}{\xi_1-\xi_2}.
\label{eq:fk-5-3b}
\end{equation}

For $a_N>0$, an integer $k\geq1$ and $C>1$, define
\begin{equation}
I^{[C]}_{N,\mu,k}(a_N):=
\int_{a_N<\xi_{k-1}<\cdots<\xi_1\leq Ca_N}
\varepsilon_N^k
\det[L_{N,y}(y+\varepsilon_N\xi_i,y+\varepsilon_N\xi_j)]_{i,j=1}^k
\prod_{i=1}^{k-1}\mathrm d\xi_i,
\qquad\xi_k=a_N.
\label{eq:fk-truncated-janossy-integral}
\end{equation}

\begin{Lemma}
\label{lem:finite-right-determinant}
Fix an integer $k\geq1$,
$C>1$ and $K\Subset(0,\infty)$, and set $\mu=a_N^{1/4}u$.
For $\boldsymbol\xi=(\xi_1,\ldots,\xi_k)$ with
$\xi_1>\cdots>\xi_k=a_N$, we have
\begin{equation}
\begin{aligned}
&\varepsilon_N^k
\det[L_{N,y}(y+\varepsilon_N\xi_i,y+\varepsilon_N\xi_j)]_{i,j=1}^k\\
&\qquad=2^{-k(2k+1)}\pi^{-k}a_N^{-k^2}\Delta(\boldsymbol\xi)^2\\
&\qquad\quad\times
e^{-kN\Phi_{\mu\varepsilon_N}(y+a_N\varepsilon_N)-2\sqrt{a_N}\sum_{i=1}^{k-1}(\xi_i-a_N)}[1+o(1)],
\end{aligned}
\label{eq:fk-finite-marked-chamber}
\end{equation}
uniformly for $u\in K$ and $\xi_1-a_N=O(a_N^{-1/2})$.
Moreover,
\begin{equation}
I^{[C]}_{N,\mu,k}(a_N)
=c_ka_N^{-(3k^2-1)/2}e^{-kN\Phi_{\mu\varepsilon_N}(y+a_N\varepsilon_N)}[1+o(1)],
\label{eq:fk-5-8}
\end{equation}
uniformly for $u\in K$, with $c_k$ defined in \eqref{eq:fk-3-15}.
\end{Lemma}

\begin{proof}
Integration of \eqref{eq:K-prefactor} from $\xi=0$ and Gronwall's inequality give
\begin{equation}
\bigl[R_N(y)E_N^{(y)}(y)\bigr]^{-1}
R_N(y+\varepsilon_N\xi)E_N^{(y)}(y+\varepsilon_N\xi)-I
=a_NN^{-2/3}
\begin{pmatrix}
O(1)&O(a_N^{-1/2})\\
O(a_N^{1/2})&O(1)
\end{pmatrix},
\label{eq:fk-prefactor-estimate}
\end{equation}
on $|\xi|\leq(C+1)a_N$.
Recall from Lemma~\ref{lem:local-change-of-variables} that
\begin{align}
\zeta_{N,\mu}(y+\varepsilon_N\xi)&=\xi[1+O(a_NN^{-2/3})],\label{eq:fk-5-3a}\\
\partial_\xi\zeta_{N,\mu}(y+\varepsilon_N\xi)&=1+O(a_NN^{-2/3}),\label{eq:fk-coordinate-first-derivative}\\
\partial_\xi^j\zeta_{N,\mu}(y+\varepsilon_N\xi)&=O(N^{-2(j-1)/3}),
\qquad2\leq j\leq2k,\label{eq:fk-coordinate-higher-derivatives}\\
\mu_N(\mu\varepsilon_N)-\mu&=O(N^{-2/3}a_N^{1/2}),
\label{eq:fk-model-parameter-estimate}
\end{align}
uniformly for $u\in K$ and $\xi\in\mathbb C$ satisfying
$\operatorname{dist}(\xi,[a_N,Ca_N])\leq a_N/4$.
Combining \eqref{eq:fk-5-3a}--\eqref{eq:fk-model-parameter-estimate} and \eqref{eq:airy-G-bounds} with the
prefactor estimate \eqref{eq:fk-prefactor-estimate}, we obtain by Cauchy's formula
\begin{equation}
\label{eq:fk-prefactor-derivatives}
\begin{aligned}
&\partial_\xi^j\Bigl[
\bigl[R_N(y)E_N^{(y)}(y)\bigr]^{-1}
R_N(y+\varepsilon_N\xi)E_N^{(y)}(y+\varepsilon_N\xi)
G_{\mu_N}\bigl(\zeta_{N,\mu}(y+\varepsilon_N\xi)\bigr)-I\Bigr]\\
&\qquad=(a_NN^{-2/3}+a_N^{-1/2})a_N^{-j}
\begin{pmatrix}
O(1)&O(a_N^{-1/2})\\
O(a_N^{1/2})&O(1)
\end{pmatrix},\qquad0\leq j\leq2k-1,
\end{aligned}
\end{equation}
for $\operatorname{dist}(\xi,[a_N,Ca_N])\leq a_N/8$.
On the same domain, \eqref{eq:airy-leading-vector} and
\eqref{eq:fk-prefactor-derivatives} give
\begin{align}
\left|\partial_\xi^j\left(\frac{w_{N,\mu,2}(\xi)}{w_{N,\mu,1}(\xi)}
+\sqrt{\zeta_{N,\mu}(y+\varepsilon_N\xi)}\right)\right|
&=o(a_N^{1/2-j}),\qquad1\leq j\leq2k-1,
\label{eq:fk-5-4}\\
w_{N,\mu,1}(\xi)
&=\frac{e^{-N\Phi_{\mu\varepsilon_N}(y+\varepsilon_N\xi)/2}}{2\sqrt\pi\,\zeta_{N,\mu}(y+\varepsilon_N\xi)^{1/4}}
[1+o(1)].
\label{eq:fk-5-5}
\end{align}
For each fixed $T>0$, \eqref{eq:fk-5-3a}--\eqref{eq:fk-model-parameter-estimate} give
\[
N\bigl[\Phi_{\mu\varepsilon_N}(y+\varepsilon_N\xi)-\Phi_{\mu\varepsilon_N}(y+a_N\varepsilon_N)\bigr]
=2\sqrt{a_N}(\xi-a_N)+o(1),
\qquad a_N\leq\xi\leq a_N+Ta_N^{-1/2}.
\]
By \eqref{eq:fk-5-3b}, we factor $\prod_{i=1}^k w_{N,\mu,1}(\xi_i)^2$
out of the rescaled determinant in \eqref{eq:fk-finite-marked-chamber}
and take successive row and column divided differences of the remaining matrix.
The Hermite--Genocchi formula, as in the proof of \eqref{eq:fk-3-4}, together
with \eqref{eq:fk-5-3a}--\eqref{eq:fk-model-parameter-estimate} and
\eqref{eq:fk-5-4}--\eqref{eq:fk-5-5}, gives
\eqref{eq:fk-finite-marked-chamber} uniformly on the stated domain.
On $a_N\leq\xi\leq Ca_N$, we have
\[
N\partial_\xi\Phi_{\mu\varepsilon_N}(y+\varepsilon_N\xi)
=\left(2\sqrt{\zeta_{N,\mu}(y+\varepsilon_N\xi)}
-\frac{\mu_N}{\sqrt{\zeta_{N,\mu}(y+\varepsilon_N\xi)}}\right)
\partial_\xi\zeta_{N,\mu}(y+\varepsilon_N\xi)
\geq c\sqrt{a_N}.
\]
Using \eqref{eq:fk-5-4} to bound these divided differences and integrating
the preceding lower bound to control the factors in \eqref{eq:fk-5-5},
we obtain, for $\xi_1\leq Ca_N$ and some constant $C_1>0$,
\begin{equation}
\begin{aligned}
0\leq{}&\varepsilon_N^k
\det[L_{N,y}(y+\varepsilon_N\xi_i,y+\varepsilon_N\xi_j)]_{i,j=1}^k\\
\leq{}&C_1a_N^{-k^2}\Delta(\boldsymbol\xi)^2
e^{-kN\Phi_{\mu\varepsilon_N}(y+a_N\varepsilon_N)-c\sqrt{a_N}\sum_{i=1}^{k-1}(\xi_i-a_N)}.
\end{aligned}
\label{eq:fk-finite-marked-upper}
\end{equation}
After rescaling to $\sqrt{a_N}(\xi_i-a_N)$, we integrate
\eqref{eq:fk-finite-marked-chamber} by dominated convergence, using
\eqref{eq:fk-finite-marked-upper}.
The Laguerre integral \eqref{eq:fk-5-9}, with $a=a_N$, then gives
\eqref{eq:fk-5-8}.
\end{proof}

\begin{proof}[Proof of \eqref{eq:fixed-k-finite-right}]
As in Section~\ref{subsec:right-airy-tail}, we argue by induction on $k\geq1$.
In \eqref{eq:finite-fixed-gap-janossy}, set $y=B-\mu\varepsilon_N$ and
$x_i=y+\varepsilon_N\xi_i$, so that the endpoint
$x_k=y+a_N\varepsilon_N$ becomes $\xi_k=a_N$. We obtain
\begin{align}
\mathcal J_{N,k}(\mu)
&:=D_N(y)
\int_{\xi_1>\cdots>\xi_{k-1}>a_N}
\varepsilon_N^k
\det[L_{N,y}(y+\varepsilon_N\xi_i,y+\varepsilon_N\xi_j)]_{i,j=1}^k
\prod_{i=1}^{k-1}\mathrm d\xi_i,\label{eq:fk-finite-janossy-integrand}\\
\mathbb P(\Delta_{N,k}>a_N)
&=\int_{\mathbb R}\mathcal J_{N,k}(\mu)\,\mathrm d\mu.
\label{eq:fk-5-2}
\end{align}

Choose $C_0>1$ sufficiently large,
and then fix $C>C_0+1$ such that $(C-C_0)^{3/2}>k$.
For $\mu=a_N^{1/4}u$ with $u$ in a fixed compact set
$K\Subset(0,\infty)$, \eqref{eq:det-comparison} gives
\begin{equation}
\log\frac{D_N(B-\mu\varepsilon_N)}{F_2(-\mu)}
=O\bigl(N^{-2/3}(1+a_N)\bigr)=o(1).
\label{eq:fk-5-10}
\end{equation}
The first-variation calculation in \cite[proof of Lemma~2.2]{RZ} gives
$\partial_\delta\Phi_\delta'(z)|_{\delta=0}=0$ at fixed $z>B$.
By \eqref{eq:hard-phase} and $\Phi_\delta(B-\delta)=0$,
$\Phi_\delta(B-\delta+t^2)$ extends from $t>0$ to a function jointly
analytic in $(\delta,t)$ and odd in $t$. Hence
\[
\Phi_\delta(B-\delta+a_N\varepsilon_N)
=\Phi(B+a_N\varepsilon_N)-\delta\Phi'(B+a_N\varepsilon_N)
+O(\delta^2\sqrt{a_N\varepsilon_N}).
\]
Taking $\delta=\mu\varepsilon_N$ and using \eqref{eq:exterior-phase} gives
\begin{equation}
N\Phi_{\mu\varepsilon_N}(y+a_N\varepsilon_N)=N\Phi(B+a_N\varepsilon_N)
-a_N^{3/4}u\Gamma(a_N\varepsilon_N)+O(a_NN^{-2/3}).
\label{eq:fk-5-11}
\end{equation}

Set $u_{N,k}:=2\sqrt{k\Gamma(a_N\varepsilon_N)}$.
Lemma~\ref{lem:finite-right-determinant} together with
\eqref{eq:TW-left-full}, \eqref{eq:fk-5-10} and \eqref{eq:fk-5-11}
gives, uniformly for $u\in K$,
\begin{align*}
D_N(y)I^{[C]}_{N,\mu,k}(a_N)
={}&2^{1/24}e^{\zeta'(-1)}c_k
a_N^{-(3k^2-1)/2-1/32}u^{-1/8}\\
&\quad\times
e^{-kN\Phi(B+a_N\varepsilon_N)+a_N^{3/4}(k\Gamma(a_N\varepsilon_N)u-u^3/12)}[1+o(1)].
\end{align*}
The function $k\Gamma(a_N\varepsilon_N)u-u^3/12$ attains its unique
maximum on $(0,\infty)$ at $u=u_{N,k}$, where its second derivative is
$-\sqrt{k\Gamma(a_N\varepsilon_N)}$. The difference between its maximum
value and its value at $u$ factors as
\begin{equation}
\frac43k^{3/2}\Gamma(a_N\varepsilon_N)^{3/2}
-k\Gamma(a_N\varepsilon_N)u+\frac{u^3}{12}
=\frac{(u-u_{N,k})^2(u+2u_{N,k})}{12}.
\label{eq:fk-finite-phase-factorization}
\end{equation}
Fix $0<\eta<1$. Applying Laplace's method as in
Section~\ref{subsec:right-airy-tail} and using
$\Gamma(a_N\varepsilon_N)\to2$, we obtain
\begin{align}
&\int_{|\mu/(a_N^{1/4}u_{N,k})-1|\leq\eta}
D_N(y)I^{[C]}_{N,\mu,k}(a_N)\,\mathrm d\mu
\notag\\
&\quad=C_ka_N^{-(48k^2-11)/32}
e^{-kN\Phi(B+a_N\varepsilon_N)+\frac43k^{3/2}a_N^{3/4}\Gamma(a_N\varepsilon_N)^{3/2}}[1+o(1)].
\label{eq:fk-finite-central-integral}
\end{align}

It remains to estimate the contribution outside the central $\mu$-interval
in \eqref{eq:fk-finite-central-integral} and the error from truncating at
$\xi_1=Ca_N$.
As in \eqref{eq:fk-3-19}, we bound $\det[K_N(x_i,x_j)]_{i,j=1}^k$
by Hadamard's inequality and use
Lemma~\ref{lem:ordinary-right-edge-count} to obtain
\begin{equation}
\begin{aligned}
\int_{-\infty}^0\mathcal J_{N,k}(\mu)\,\mathrm d\mu
&\leq\Prob(\lambda_k>B+a_N\varepsilon_N)\\
&\leq\frac1{k!}\left(
\int_{B+a_N\varepsilon_N}^{\infty}K_N(x,x)\,\mathrm dx\right)^k
\leq C_1e^{-kN\Phi(B+a_N\varepsilon_N)}.
\end{aligned}
\label{eq:fk-finite-negative-mu}
\end{equation}
Applying Lemma~\ref{prop:outer-B} and
\eqref{eq:K-prefactor} to \eqref{eq:full-recovery-exact} gives
\begin{equation}
\varepsilon_NL_{N,y}
\bigl(y+\varepsilon_N\xi,y+\varepsilon_N\xi\bigr)
=O\left(a_N^{1/2}e^{-N\Phi_{\mu\varepsilon_N}(y+\varepsilon_N\xi)}\right),
\label{eq:full-kernel-upper}
\end{equation}
uniformly for $0\leq\mu\leq C_0a_N$ and $a_N\leq\xi\leq Ca_N$.
In this regime, Lemma~\ref{lem:local-change-of-variables} gives
\begin{equation}
N\partial_\mu\bigl[\Phi_{\mu\varepsilon_N}(B+(\xi-\mu)\varepsilon_N)\bigr]=O(\sqrt{a_N}).
\label{eq:full-hard-wall-derivative}
\end{equation}
Since $\Phi_0=\Phi$ and $\Phi(B+\xi\varepsilon_N)$ is increasing in $\xi$,
integration in $\mu$ yields
\[
N\Phi_{\mu\varepsilon_N}(y+\varepsilon_N\xi)\geq N\Phi(B+a_N\varepsilon_N)-C_1\mu\sqrt{a_N}.
\]
Hadamard's inequality and \eqref{eq:full-kernel-upper}, combined with
the preceding phase bound, \eqref{eq:TW-left-global} and
\eqref{eq:det-comparison}, give
\[
 D_N(y)I^{[C]}_{N,\mu,k}(a_N)
 \leq a_N^{C_1}
 e^{-kN\Phi(B+a_N\varepsilon_N)-c\mu^3+C_1\mu\sqrt{a_N}},
\]
uniformly for $0\leq\mu\leq C_0a_N$; the comparison error is absorbed
using $a_NN^{-2/3}\to0$.
Choose fixed $m>0$ sufficiently small and $U$ sufficiently large.
Combining this bound for $\mu/a_N^{1/4}\notin[m,U]$ with
Lemma~\ref{lem:finite-right-determinant} and
\eqref{eq:fk-5-10}--\eqref{eq:fk-finite-phase-factorization} for $\mu/a_N^{1/4}\in[m,U]$,
we obtain
\begin{equation}
\begin{aligned}
&\int_{\substack{0\leq\mu\leq C_0a_N\\
|\mu/(a_N^{1/4}u_{N,k})-1|>\eta}}
D_N(y)I^{[C]}_{N,\mu,k}(a_N)\,\mathrm d\mu\\
&\qquad\leq a_N^{C_1}
e^{-kN\Phi(B+a_N\varepsilon_N)
+a_N^{3/4}(\frac43k^{3/2}\Gamma(a_N\varepsilon_N)^{3/2}-c)}.
\end{aligned}
\label{eq:fk-finite-off-center}
\end{equation}

An event decomposition analogous to that in
\eqref{eq:fk-airy-far-left-decomposition} gives
\begin{equation}
\begin{aligned}
\int_{C_0a_N}^{\infty}\mathcal J_{N,k}(\mu)\,\mathrm d\mu
&\leq\Prob(\zeta_{N,k}\leq-(C_0-1)a_N)\\
&\leq
\begin{cases}
D_N(B-(C_0-1)a_N\varepsilon_N),&k=1,\\
D_N(B-(C_0-1)a_N\varepsilon_N/2)
+\displaystyle\sum_{j=1}^{k-1}
\Prob\left(\Delta_{N,j}\geq\dfrac{(C_0-1)a_N}{2(k-1)}\right),&k\geq2.
\end{cases}
\end{aligned}
\label{eq:fk-5-15}
\end{equation}
Equations~\eqref{eq:TW-left-global} and \eqref{eq:det-comparison},
together with the induction hypothesis, give,
for sufficiently large $C_0$,
\begin{equation}
\int_{C_0a_N}^{\infty}\mathcal J_{N,k}(\mu)\,\mathrm d\mu
\leq e^{-(4k/3+1)a_N^{3/2}}.
\label{eq:fk-finite-far-left}
\end{equation}

To bound the truncation error, we combine \eqref{eq:edge-scale} and
\eqref{eq:exterior-phase} with $\Gamma(s)=2+O(s)$ as $s\to0$ to obtain,
uniformly for $1\leq S\leq(C-C_0)a_N$,
\begin{equation}
\begin{aligned}
N\Phi(B+S\varepsilon_N)
&=\int_0^S\sqrt t\,\Gamma(t\varepsilon_N)\,\mathrm dt\\
&=\frac43S^{3/2}+O(\varepsilon_NS^{5/2}).
\end{aligned}
\label{eq:fk-finite-edge-phase}
\end{equation}
For $k\geq2$ and $0\leq\mu\leq C_0a_N$, the condition $\xi_1>Ca_N$
implies $\lambda_1>B+(C-C_0)a_N\varepsilon_N$.
Lemma~\ref{lem:ordinary-right-edge-count}
and \eqref{eq:fk-finite-edge-phase} therefore give
\begin{equation}
\begin{aligned}
&\int_0^{C_0a_N}\bigl[\mathcal J_{N,k}(\mu)
-D_N(y)I^{[C]}_{N,\mu,k}(a_N)\bigr]\,\mathrm d\mu\\
&\qquad\leq\Prob\bigl(\lambda_1>B+(C-C_0)a_N\varepsilon_N\bigr)\\
&\qquad\leq\int_{B+(C-C_0)a_N\varepsilon_N}^{\infty}
K_N(x,x)\,\mathrm dx\\
&\qquad\leq C_1 e^{-N\Phi(B+(C-C_0)a_N\varepsilon_N)}
=e^{-\frac43(C-C_0)^{3/2}a_N^{3/2}+o(a_N^{3/2})},
\end{aligned}
\label{eq:fk-finite-truncation-error}
\end{equation}
while for $k=1$, $\mathcal J_{N,1}(\mu)=D_N(y)I^{[C]}_{N,\mu,1}(a_N)$.

Comparing \eqref{eq:fk-finite-negative-mu}, \eqref{eq:fk-finite-off-center},
\eqref{eq:fk-finite-far-left}, and \eqref{eq:fk-finite-truncation-error}
with \eqref{eq:fk-finite-central-integral}, and using $(C-C_0)^{3/2}>k$, we obtain from
\eqref{eq:fk-5-2}
\[
\mathbb P(\Delta_{N,k}>a_N)
=\int_{|\mu/(a_N^{1/4}u_{N,k})-1|\leq\eta}
D_N(y)I^{[C]}_{N,\mu,k}(a_N)\,\mathrm d\mu\,[1+o(1)].
\]
Taking logarithms and applying \eqref{eq:fk-finite-central-integral} gives
\eqref{eq:fixed-k-finite-right} with the constant \eqref{eq:def-CRk}.
\end{proof}

\subsection{Proof of Corollary~\ref{cor:fixed-k-conditional-joint}}
\label{subsec:right-conditional-limits}

\begin{proof}
For the finite-$N$ ensemble, take $(z,\mathbf s)\in\mathbb R\times\mathcal S_k$
and set
\[
\begin{aligned}
 &\mu=a_N-m_{N,k}-\sigma_{N,k}z,\qquad y=B-\mu\varepsilon_N,\\
 &\xi_i=a_N+\frac{s_i}{\Gamma(a_N\varepsilon_N)a_N^{1/2}}
 \quad(1\leq i<k),\qquad \xi_k=a_N.
\end{aligned}
\]
Equations~\eqref{eq:fk-finite-janossy-integrand}--\eqref{eq:fk-5-2} give
\begin{equation}
\begin{aligned}
 \varrho_{N,k}(z,\mathbf s)
 ={}&\frac{\sigma_{N,k}[\Gamma(a_N\varepsilon_N)a_N^{1/2}]^{-(k-1)}
 D_N(y)}{\Prob(\Delta_{N,k}>a_N)}\\
 &\quad\times\varepsilon_N^k
 \det[L_{N,y}(y+\varepsilon_N\xi_i,y+\varepsilon_N\xi_j)]_{i,j=1}^k.
\end{aligned}
\label{eq:fk-finite-conditional-density}
\end{equation}
The above change of variables gives
\[
 u=\frac{\mu}{a_N^{1/4}}
 =u_{N,k}-(k\Gamma(a_N\varepsilon_N))^{-1/4}a_N^{-3/8}z.
\]
Substituting into \eqref{eq:fk-finite-phase-factorization} yields,
locally uniformly for $z\in\mathbb R$,
\[
 -\frac{a_N^{3/4}}{12}(u-u_{N,k})^2(u+2u_{N,k})
 =-\frac{z^2}{2}+o(1).
\]
Combining this expansion with \eqref{eq:TW-left-full},
\eqref{eq:fk-finite-marked-chamber} and
\eqref{eq:fk-5-10}--\eqref{eq:fk-5-11} gives the numerator asymptotics in
\eqref{eq:fk-finite-conditional-density}.
Using \eqref{eq:fixed-k-finite-right} for the denominator and
$\Gamma(a_N\varepsilon_N)\to2$, we obtain,
locally uniformly on $\mathbb R\times\mathcal S_k$,
\begin{equation}
 \varrho_{N,k}(z,\mathbf s)=\varrho_k(z,\mathbf s)[1+o(1)].
\label{eq:fk-conditional-local-density}
\end{equation}
Scheff\'e's lemma now gives
\[
 \|\varrho_{N,k}-\varrho_k\|_{L^1(\mathbb R\times\mathcal S_k)}\longrightarrow0.
\]
Since
\[
 \Gamma(a_N\varepsilon_N)a_N^{1/2}\sigma_{N,k}
 =k^{-1/4}\Gamma(a_N\varepsilon_N)^{3/4}a_N^{3/8}
 \longrightarrow\infty,
\]
the tightness of the $S_{N,i}$ under the conditional law gives
\[
 \max_{i<k}\left|
 \frac{\zeta_{N,i}-m_{N,k}}{\sigma_{N,k}}-Z_{N,k}\right|
 =\frac{\max_{i<k}S_{N,i}}
 {\Gamma(a_N\varepsilon_N)a_N^{1/2}\sigma_{N,k}}
 \xrightarrow{\mathbb P}0,
\]
which proves \eqref{eq:fixed-k-common-gaussian}.

For the Airy point process, take $(z,\mathbf s)\in\mathbb R\times\mathcal S_k$
and set
\[
 \mu=a-m_k(a)-\sigma_k(a)z,\qquad
 x_i=a+\frac{s_i}{2a^{1/2}}\quad(1\leq i<k),\qquad x_k=a.
\]
Equation~\eqref{eq:fk-3-23} gives
\[
 \varrho_{a,k}(z,\mathbf s)
 =\frac{\sigma_k(a)(2a^{1/2})^{-(k-1)}F_2(-\mu)}
 {\mathbb P(\Delta_{\mathrm{Ai},k}>a)}
 \det[L_{-\mu}(x_i-\mu,x_j-\mu)]_{i,j=1}^k.
\]
With this choice of $\mu$, we have, locally uniformly for $z\in\mathbb R$,
\[
 -\frac{\mu^3}{12}+2k\mu a^{1/2}
 =\frac{8\sqrt2}{3}k^{3/2}a^{3/4}
 -\frac{z^2}{2}+o(1).
\]
Equations~\eqref{eq:fixed-k-airy-right}--\eqref{eq:def-CRk},
\eqref{eq:fk-3-4}, and \eqref{eq:TW-left-full} yield
\[
 \varrho_{a,k}(z,\mathbf s)\longrightarrow\varrho_k(z,\mathbf s),
 \qquad (z,\mathbf s)\in\mathbb R\times\mathcal S_k.
\]
Scheff\'e's lemma gives convergence in $L^1(\mathbb R\times\mathcal S_k)$.
Since $2a^{1/2}\sigma_k(a)=2(2k)^{-1/4}a^{3/8}\to\infty$,
tightness of $2a^{1/2}(\zeta_i-\zeta_k)$, $1\leq i<k$, under the conditional law gives
\eqref{eq:airy-common-gaussian} by the argument following
\eqref{eq:fk-conditional-local-density}.
\end{proof}

\newpage
\section{Left-Tail Asymptotics and Conditional Limits}
\label{sec:left-proof}
We first evaluate the integral of the collision density $q_k$ in
\eqref{eq:def-left-qk}, which determines the cubic coefficient.

\begin{Lemma}
\label{lem:collision-density-mass}
For every fixed integer $k\ge1$, the function $q_k$ in
\eqref{eq:def-left-qk} satisfies
\begin{equation}
 \int_{\R}q_k(\mu)\,\dd\mu=\frac{k}{2}.
\label{eq:def-main-left-coefficient}
\end{equation}
\end{Lemma}

\begin{proof}
We compute the masses of $q_k$ by integrating their generating function.
For $0<t<1$, define the deformed gap probability and resolvent operator
\[
 F_t(\mu):=\det(I-(1-t) K_{\mathrm{Ai},\mu}),\qquad
 L_{t,\mu}:=(1-t) K_{\mathrm{Ai},\mu}
 (I-(1-t) K_{\mathrm{Ai},\mu})^{-1},
\]
where $L_{t,\mu}(x,y)$ denotes the resolvent kernel.

The parameterized Fredholm minor identity \cite[Section~2]{Nishigaki2021},
with parameter $1-t$, gives the following identity after the two prescribed
points coalesce at $\mu$:
\begin{equation}
 \sum_{j=1}^{\infty}t^{j-1}q_j(\mu)
 =(1-t)^{-2}F_t(\mu)\Bigl[
 L_{t,\mu}(\mu,\mu)\partial_1\partial_2L_{t,\mu}(\mu,\mu)
 -\bigl(\partial_1L_{t,\mu}(\mu,\mu)\bigr)^2\Bigr].
\label{eq:collision-mass-generating}
\end{equation}
To justify this limit, fix $\mu\in\mathbb R$ and take the prescribed points
at $\mu+h$ and $\mu$, with $0<h\leq1$. For $n\geq0$, the coefficient of
$t^n$ in the Fredholm expansion before coalescence is
\[
 \frac{F_2(\mu)}{n!}
 \int_{(\mu,\infty)^n}\det[L_\mu(x_i,x_j)]_{i,j=1}^{n+2}
 \prod_{i=1}^n\dd x_i,
 \qquad x_{n+1}=\mu+h,\quad x_{n+2}=\mu.
\]
We bound this coefficient divided by $h^2$ uniformly in $h$ using a Gram
representation of $L_\mu$. Define $\mathcal A_\mu$ on $L^2((0,\infty))$ by
\begin{equation}
 (\mathcal A_\mu f)(r):=\int_0^\infty
 \operatorname{Ai}(\mu+r+s)f(s)\,\mathrm ds.
\label{eq:airy-hankel-operator}
\end{equation}
The kernel of $\mathcal A_\mu^2$ is
\[
 \int_0^\infty\operatorname{Ai}(\mu+r+t)
 \operatorname{Ai}(\mu+s+t)\,\mathrm dt
 =K_{\mathrm{Ai}}(\mu+r,\mu+s),\qquad r,s>0.
\]
The operator $\mathcal A_\mu^2$ is a positive trace-class strict contraction.
Write $a_x:=\operatorname{Ai}(x+\cdot)$ in $L^2((0,\infty))$ and define,
with $\mu$ fixed,
\begin{equation}
\begin{aligned}
 w_x&:=(I-\mathcal A_\mu^2)^{-1/2}a_x,\\
 w_x^{(j)}&:=\frac{\mathrm d^j}{\mathrm dx^j}w_x
 =(I-\mathcal A_\mu^2)^{-1/2}\operatorname{Ai}^{(j)}(x+\cdot),
 \qquad j=0,1,2.
\end{aligned}
\label{eq:airy-resolvent-gram-vectors}
\end{equation}
Then
\[
 L_\mu(x,y)=\langle w_x,w_y\rangle,\qquad
 \|w_x^{(j)}\|^2=\partial_1^j\partial_2^jL_\mu(x,x).
\]

The integrals of the Gram determinants are expressed by elementary
symmetric sums. For a positive trace-class operator $T$ with eigenvalues
$\tau_1,\tau_2,\ldots$, counted with multiplicity, define
\begin{equation}
 e_m(T):=\sum_{i_1<\cdots<i_m}\tau_{i_1}\cdots\tau_{i_m},
 \qquad m\geq1,\qquad e_0(T):=1.
\label{eq:elementary-symmetric-sums}
\end{equation}
In the coefficient integral above, subtracting the last row and column
from their $(n+1)$st counterparts and applying the Gram inequality bounds
the coefficient divided by $h^2$ by
\[
 F_2(\mu)\|w_\mu\|^2
 \sup_{0\leq s\leq1}\|w_{\mu+s}'\|^2\,e_n(L_\mu).
\]
Since $(w_{\mu+h}-w_\mu)/h\to w_\mu'$ in $L^2((0,\infty))$,
the coefficient divided by $h^2$ converges to $q_{n+1}(\mu)$ as
$h\downarrow0$. The bound is summable because
$\sum_{n\geq0}t^ne_n(L_\mu)=\det(I+tL_\mu)<\infty$, so dominated
convergence gives \eqref{eq:collision-mass-generating}.
Replacing $(I-\mathcal A_\mu^2)^{-1/2}$ in the definition of $w_x$ by
$\sqrt{1-t}(I-(1-t)\mathcal A_\mu^2)^{-1/2}$ gives the Gram
representation of $L_{t,\mu}$. These representations also give $q_j\ge0$ and
\[
 0\le L_{t,\mu}(\mu,\mu)\partial_1\partial_2L_{t,\mu}(\mu,\mu)
 -\bigl(\partial_1L_{t,\mu}(\mu,\mu)\bigr)^2
 \le L_{t,\mu}(\mu,\mu)\partial_1\partial_2L_{t,\mu}(\mu,\mu).
\]

To integrate the right-hand side of \eqref{eq:collision-mass-generating},
let $u$ be the Ablowitz--Segur solution of Painlev\'e II satisfying
\begin{equation}
 u''=\mu u+2u^3,\qquad
 u(\mu)\sim\sqrt{1-t}\operatorname{Ai}(\mu)
 \quad(\mu\to+\infty),
\label{eq:collision-mass-pii}
\end{equation}
where primes denote derivatives with respect to $\mu$.
The resolvent identities in \cite[Sections~II--III]{TW94} give
\[
\begin{aligned}
 \frac{\mathrm d}{\mathrm d\mu}L_{t,\mu}(\mu,\mu)&=-u^2,\\
 \frac{\mathrm d}{\mathrm d\mu}\bigl[\partial_1L_{t,\mu}(\mu,\mu)\bigr]
 &=-u\bigl(u'+L_{t,\mu}(\mu,\mu)u\bigr).
\end{aligned}
\]
We integrate these identities using the Painlev\'e II equation and determine
the integration constants from the boundary conditions at $+\infty$.
This gives
\[
\begin{aligned}
 L_{t,\mu}(\mu,\mu)&=(u')^2-u^4-\mu u^2,\\
 \partial_1L_{t,\mu}(\mu,\mu)
 &=\frac{L_{t,\mu}(\mu,\mu)^2-u^2}{2}.
\end{aligned}
\]
To write the integrand as a total derivative, set
\[
 H:=\partial_1\partial_2L_{t,\mu}(\mu,\mu)
 +\frac12uu'+\frac34u^2L_{t,\mu}(\mu,\mu)
 -\frac14L_{t,\mu}(\mu,\mu)^3.
\]
Differentiating $F_tH$ with the aid of
\[
 F_t'=L_{t,\mu}(\mu,\mu)F_t,\qquad
 \frac{\mathrm d}{\mathrm d\mu}\bigl[\partial_1\partial_2L_{t,\mu}(\mu,\mu)\bigr]
 =-\bigl(u'+L_{t,\mu}(\mu,\mu)u\bigr)^2
\]
and the preceding identities yields
\begin{equation}
 F_t\Bigl[L_{t,\mu}(\mu,\mu)\partial_1\partial_2L_{t,\mu}(\mu,\mu)
 -\bigl(\partial_1L_{t,\mu}(\mu,\mu)\bigr)^2\Bigr]
 =\frac12F_t'+(F_t H)'.
\label{eq:collision-mass-primitive}
\end{equation}
For $t=0$, the corresponding integration by parts appears in
\cite[Appendix~B.2]{PerretSchehr}.

To integrate \eqref{eq:collision-mass-primitive} over $\R$, we verify
integrability and the vanishing of the boundary term $F_tH$.
For fixed $0<t<1$, the resolvent satisfies
\[
 \|(I-(1-t) K_{\mathrm{Ai},\mu})^{-1}\|\le t^{-1}.
\]
The Airy bounds give, as $\mu\to-\infty$,
\[
 \|a_\mu\|^2=O(|\mu|^{1/2}),\qquad
 \|a_\mu'\|^2=O(|\mu|^{3/2}),
\]
and the Gram representation of $L_{t,\mu}$ implies
\[
 L_{t,\mu}(\mu,\mu)=O_t(|\mu|^{1/2}),\qquad
 \partial_1\partial_2L_{t,\mu}(\mu,\mu)=O_t(|\mu|^{3/2}).
\]
To bound the remaining factors $u$ and $u'$ in $H$, we use the resolvent
representation
\[
\begin{aligned}
 Q_\mu(x)&:=\sqrt{1-t}\,
 [(I-(1-t) K_{\mathrm{Ai},\mu})^{-1}\operatorname{Ai}](x),\\
 u(\mu)&=Q_\mu(\mu),\qquad
 \partial_xQ_\mu(\mu)=u'+L_{t,\mu}(\mu,\mu)u.
\end{aligned}
\]
The function $Q_\mu$ satisfies
\[
 Q_\mu(x)=\sqrt{1-t}\operatorname{Ai}(x)
 +(1-t)\int_\mu^\infty K_{\mathrm{Ai}}(x,y)Q_\mu(y)\,\dd y,
 \qquad
 \|Q_\mu\|_2\le\frac{\sqrt{1-t}}{t}\|a_\mu\|_2.
\]
Since $\|\mathcal A_\mu\|\le1$, the identities
\[
 K_{\mathrm{Ai}}(\mu,\mu+\cdot)=\mathcal A_\mu a_\mu,\qquad
 \partial_1K_{\mathrm{Ai}}(\mu,\mu+\cdot)=\mathcal A_\mu a_\mu'
\]
bound the norms of these two vectors by $\|a_\mu\|$ and $\|a_\mu'\|$,
respectively.
Evaluating the integral equation and its derivative with respect to $x$
at $x=\mu$ and applying these bounds gives
\[
 u=O_t(|\mu|^{1/2}),\qquad
 u'+L_{t,\mu}(\mu,\mu)u=O_t(|\mu|),\qquad u'=O_t(|\mu|).
\]
Thus $H=O_t(|\mu|^{3/2})$ and
\[
 L_{t,\mu}(\mu,\mu)\partial_1\partial_2L_{t,\mu}(\mu,\mu)
 -\bigl(\partial_1L_{t,\mu}(\mu,\mu)\bigr)^2=O_t(|\mu|^2).
\]
This polynomial growth is dominated by the decay of $F_t$, since
\[
 0<F_t(\mu)\le
 \exp[-(1-t)\operatorname{Tr}K_{\mathrm{Ai},\mu}],\qquad
 \operatorname{Tr}K_{\mathrm{Ai},\mu}
 \sim\frac{2}{3\pi}|\mu|^{3/2}\quad(\mu\to-\infty).
\]
At $+\infty$, the same operator formulas and the Airy bounds give
exponential decay of $L_{t,\mu}(\mu,\mu)$,
$\partial_1\partial_2L_{t,\mu}(\mu,\mu)$, $u$ and $u'$, while $F_t\to1$.
These estimates imply that the right-hand side of
\eqref{eq:collision-mass-generating} is integrable on $\R$
and $F_tH$ vanishes at both endpoints, so integrating
\eqref{eq:collision-mass-primitive} yields
\[
 \int_{\R}F_t(\mu)\Bigl[
 L_{t,\mu}(\mu,\mu)\partial_1\partial_2L_{t,\mu}(\mu,\mu)
 -\bigl(\partial_1L_{t,\mu}(\mu,\mu)\bigr)^2\Bigr]\,\dd\mu=\frac12.
\]
Since $q_j\ge0$, Tonelli's theorem applied to
\eqref{eq:collision-mass-generating}, together with the preceding integral
evaluation, gives
\[
 \sum_{j=1}^{\infty}t^{j-1}\int_{\R}q_j(\mu)\,\dd\mu
 =\frac1{2(1-t)^2}
 =\frac12\sum_{j=1}^{\infty}j\,t^{j-1},\qquad 0<t<1.
\]
The series on the left is finite for every $0<t<1$, so uniqueness of
power-series coefficients proves \eqref{eq:def-main-left-coefficient}.
\end{proof}

\subsection{Proof of \eqref{eq:fixed-k-airy-left}}

\begin{proof}[Proof of \eqref{eq:fixed-k-airy-left}]
With $\mathbf s=(z_1,\ldots,z_{k-1},\mu+\varepsilon,\mu)$, the ordered Janossy
density gives
\begin{equation}
\begin{aligned}
\mathbb P(\Delta_{\mathrm{Ai},k}<b)
={}&\int_{\mathbb R}\int_0^b
F_2(\mu)\int_{z_1>\cdots>z_{k-1}>\mu+\varepsilon}
\det[L_\mu(s_i,s_j)]_{i,j=1}^{k+1}\\
&\hspace{31mm}\times
\prod_{i=1}^{k-1}\mathrm dz_i\,\mathrm d\varepsilon\,\mathrm d\mu\\
:={}&\int_{\mathbb R}\int_0^b
I_k(\mu,\varepsilon)\,\mathrm d\varepsilon\,\mathrm d\mu.
\end{aligned}
\label{eq:fk-4-1}
\end{equation}

We use the Gram vectors $w_x$ in \eqref{eq:airy-resolvent-gram-vectors}.
Fix $b_0>0$ sufficiently small, independently of $\mu$.
For $0<\varepsilon\leq b_0$, set
\[
 d_\varepsilon:=\frac{w_{\mu+\varepsilon}-w_\mu}{\varepsilon}
 =\int_0^1w_{\mu+t\varepsilon}'\,\mathrm dt,
 \qquad
 \partial_\varepsilon d_\varepsilon
 =\int_0^1t\,w_{\mu+t\varepsilon}''\,\mathrm dt.
\]
With $(v_1,\ldots,v_{k+1})
=(w_{z_1},\ldots,w_{z_{k-1}},d_\varepsilon,w_\mu)$,
row and column subtraction in the Gram determinant in \eqref{eq:fk-4-1} gives
\begin{equation}
\frac{I_k(\mu,\varepsilon)}{\varepsilon^2}
=F_2(\mu)
\int_{z_1>\cdots>z_{k-1}>\mu+\varepsilon}
\det[\langle v_i,v_j\rangle]_{i,j=1}^{k+1}
\prod_{i=1}^{k-1}\mathrm dz_i.
\label{eq:fk-4-5}
\end{equation}
For $z_1>\cdots>z_{k-1}>\mu$,
\begin{align*}
 &\det[\langle v_i,v_j\rangle]_{i,j=1}^{k+1}\\
 &\qquad\leq L_\mu(\mu,\mu)
 \sup_{0\leq t\leq b_0}\|w_{\mu+t}'\|^2
 \det[L_\mu(z_i,z_j)]_{i,j=1}^{k-1}.
\end{align*}
With $e_m$ defined in \eqref{eq:elementary-symmetric-sums}, the ordered
trace identity for $L_\mu$ on $L^2((\mu,\infty))$ is
\[
 \int_{z_1>\cdots>z_{k-1}>\mu}
 \det[L_\mu(z_i,z_j)]_{i,j=1}^{k-1}
 \prod_{i=1}^{k-1}\mathrm dz_i
 =e_{k-1}(L_\mu).
\]
Since $d_\varepsilon\to w_\mu'$ in $L^2((0,\infty))$,
dominated convergence and \eqref{eq:def-left-qk} give
\[
 \lim_{\varepsilon\downarrow0}
 \frac{I_k(\mu,\varepsilon)}{\varepsilon^2}
 =F_2(\mu)\int_{z_1>\cdots>z_{k-1}>\mu}
 \det\mathcal M_{k,\mu}(\mathbf z)
 \prod_{i=1}^{k-1}\mathrm dz_i=q_k(\mu),
\]
with
\begin{equation}
0\leq q_k(\mu)
\leq F_2(\mu)e_{k-1}(L_\mu)L_\mu(\mu,\mu)
\partial_1\partial_2L_\mu(\mu,\mu).
\label{eq:fk-4-8}
\end{equation}

To prove $q_k(\mu)>0$, we use the following Airy asymptotics for fixed $x<y$:
\[
 \frac{\operatorname{Ai}(y+t)}{\operatorname{Ai}(x+t)}
 \sim e^{-(y-x)\sqrt t},\qquad
 \frac{\operatorname{Ai}'(\mu+t)}{\operatorname{Ai}(\mu+t)}
 \sim-\sqrt t,\qquad t\to\infty.
\]
These distinct asymptotic orders imply linear independence of
$\operatorname{Ai}(z_1+\cdot),\ldots,
\operatorname{Ai}(z_{k-1}+\cdot),\operatorname{Ai}'(\mu+\cdot)$ and
$\operatorname{Ai}(\mu+\cdot)$.
Applying the injective operator $(I-\mathcal A_\mu^2)^{-1/2}$ preserves
this linear independence, so the Gram matrix
$\mathcal M_{k,\mu}(\mathbf z)$ is positive definite and
\[
 \det\mathcal M_{k,\mu}(\mathbf z)>0,
 \qquad q_k(\mu)>0.
\]

Differentiating \eqref{eq:fk-4-5} in $\varepsilon$ produces no boundary term
at $z_{k-1}=\mu+\varepsilon$ for $k\geq2$, since
$w_{\mu+\varepsilon}=w_\mu+\varepsilon d_\varepsilon$ makes the Gram
determinant vanish there. Differentiating its $k$th row and $k$th column
and applying Cauchy--Schwarz to the two resulting mixed Gram determinants yields
\begin{equation}
\begin{aligned}
\left|\partial_\varepsilon
\frac{I_k(\mu,\varepsilon)}{\varepsilon^2}\right|
&\leq2F_2(\mu)e_{k-1}(L_\mu)L_\mu(\mu,\mu)
 \|d_\varepsilon\|\,\|\partial_\varepsilon d_\varepsilon\|\\
&\leq CF_2(\mu)e_{k-1}(L_\mu)L_\mu(\mu,\mu)
 \prod_{j=1}^2\sup_{0\leq t\leq b_0}\|w_{\mu+t}^{(j)}\|
 =:g_k(\mu).
\end{aligned}
\label{eq:fk-4-6}
\end{equation}
Integrating \eqref{eq:fk-4-6} gives
\begin{equation}
 \left|I_k(\mu,\varepsilon)-\varepsilon^2q_k(\mu)\right|
 \leq\varepsilon^3g_k(\mu),\qquad0<\varepsilon\leq b_0.
\label{eq:fk-4-4}
\end{equation}

Let $N_I$ be the number of Airy particles in $I$. The Janossy trace identity is
\[
 F_2(\mu)e_{k-1}(L_\mu)
 =\mathbb P\{N_{(\mu,\infty)}=k-1\}\leq1.
\]
For $\mu>0$ and $k\geq2$,
$\zeta_1-\zeta_k=\sum_{j=1}^{k-1}\Delta_{\mathrm{Ai},j}$ gives
\begin{equation}
\begin{aligned}
\mathbb P\{N_{(-\mu,\infty)}=k-1\}
&\leq\mathbb P\{\zeta_k\leq-\mu\}\\
&\leq F_2(-\mu/2)
+\sum_{j=1}^{k-1}\mathbb P\left\{
\Delta_{\mathrm{Ai},j}\geq\frac{\mu}{2(k-1)}\right\}.
\end{aligned}
\label{eq:fk-4-11}
\end{equation}
Theorem~\ref{thm:fixed-k-right-tail} and \eqref{eq:TW-left-global} yield
\begin{equation}
F_2(-\mu)e_{k-1}(L_{-\mu})
\leq Ce^{-c\mu^{3/2}},\qquad \mu\geq0,
\label{eq:fk-4-13}
\end{equation}
including $k=1$, when $e_0(L_{-\mu})=1$.
Cauchy's estimates applied to \eqref{eq:fk-2-3}, with the Bessel
factorization in the proof of Lemma~\ref{lem:wall-Wronskian}, give
\[
 \sup_{0\leq s,t\leq b_0,\ j+\ell\leq4}
 |\partial_1^j\partial_2^\ell
 L_{-\mu}(-\mu+s,-\mu+t)|
 \leq Ce^{C\mu}.
\]
Thus \eqref{eq:fk-4-8}, \eqref{eq:fk-4-6} and \eqref{eq:fk-4-13} imply
\[
 q_k(-\mu)+g_k(-\mu)
 \leq Ce^{-c\mu^{3/2}+C\mu},\qquad \mu\geq0.
\]

Since $\operatorname{Tr}K_{\mathrm{Ai},\mu}\to0$ as $\mu\to\infty$,
the Airy asymptotics give, for sufficiently large $\mu$,
\begin{align*}
 \|(I-\mathcal A_\mu^2)^{-1}\|&\leq2,\\
 \sup_{0\leq t\leq b_0}\|w_{\mu+t}^{(j)}\|^2
 &\leq2\int_\mu^\infty|\operatorname{Ai}^{(j)}(r)|^2\,\mathrm dr
 \leq Ce^{-c\mu^{3/2}},\qquad j=0,1,2,\\
 q_k(\mu)+g_k(\mu)&\leq Ce^{-c\mu^{3/2}}.
\end{align*}
Trace-norm continuity of $\mathcal A_\mu^2$, $F_2(\mu)>0$ and
$L^2$-continuity of the translated Airy function and its first two
derivatives imply that $q_k$ and $g_k$ are locally bounded.
Together with the preceding tail bounds, this gives
$q_k,g_k\in L^1(\mathbb R)$. Using
\eqref{eq:def-main-left-coefficient} and \eqref{eq:fk-4-1}, we integrate
\eqref{eq:fk-4-4} to obtain
\begin{align*}
 \left|\mathbb P(\Delta_{\mathrm{Ai},k}<b)-\frac{k}{6}b^3\right|
 &\leq\int_0^b\varepsilon^3\,\mathrm d\varepsilon
 \int_{\mathbb R}g_k(\mu)\,\mathrm d\mu
 =\frac{b^4}{4}\|g_k\|_{L^1(\mathbb R)},\\
 \left|\int_{\mathbb R}I_k(\mu,\varepsilon)\,\mathrm d\mu
 -\frac{k}{2}\varepsilon^2\right|
 &\leq\varepsilon^3\|g_k\|_{L^1(\mathbb R)},
 \qquad0<b,\varepsilon\leq b_0.
\end{align*}
\end{proof}

\subsection{Proof of \eqref{eq:fixed-k-finite-left}}

For each integer $m\geq1$ and $t>0$, let $\mathbb P_{m,t}$ be the
law on $\lambda_1>\cdots>\lambda_m$ with density proportional to
$\Delta(\lambda)^2e^{-t\sum_{j=1}^mV(\lambda_j)}$, and let
$\mathbb E_{m,t}$ denote expectation with respect to $\mathbb P_{m,t}$.
For $1\leq n\leq N$, Andr\'eief's
identity \cite{Andreief1886} gives
\begin{equation}
 Z_{n,N}:=\int_{\lambda_1>\cdots>\lambda_n}
 \Delta(\lambda)^2e^{-N\sum_{j=1}^nV(\lambda_j)}\,\dd\lambda
 =\prod_{j=0}^{n-1}h_j(+\infty),
\label{eq:fk-general-partition-product}
\end{equation}
where $h_j(+\infty)$ is defined by \eqref{eq:truncated-op-norm}
with $y=+\infty$.

\begin{Lemma}
\label{lem:fixed-k-partition-coercivity}
Fix a positive integer $r$. For every integer $N\geq r+1$, with $m=N-r$,
\begin{align}
 \frac{Z_{m,N}}{Z_{N,N}}&\leq e^{CN},
\label{eq:fk-general-partition-ratio}\\
 \mathbb E_{m,N}\prod_{j=1}^m(1+|\lambda_j|)^{4r}
 &\leq e^{CN}.
\label{eq:fk-general-coercive-moment}
\end{align}
For each fixed $c_0>0$, there is an integer $N_0\geq r+1$ such that
\begin{equation}
 \mathbb P_{m,N}\{\lambda_1\leq B-c_0N^{-1/6}\}
 \leq Ce^{-cN^{3/2}},\qquad N\geq N_0.
\label{eq:fk-general-reduced-lower}
\end{equation}
\end{Lemma}

\begin{proof}
Fix a compact interval $I$ of positive length. The extremal property of
monic orthogonal polynomials and the monic Legendre norm give
\begin{align*}
 h_j(+\infty)
 &\geq e^{-N\max_I V}
 \inf_{\substack{p\ \mathrm{monic}\\\deg p=j}}\int_Ip(x)^2\,\dd x\\
 &=e^{-N\max_I V}\frac{|I|^{2j+1}}{(2j+1)\binom{2j}{j}^2}
 \geq e^{-CN},\qquad N-r\leq j<N,\\
 \frac{Z_{m,N}}{Z_{N,N}}
 &=\prod_{j=N-r}^{N-1}h_j(+\infty)^{-1}\leq e^{CN}.
\end{align*}
The Radon--Nikodym derivative is
\begin{equation}
 \frac{\dd\mathbb P_{m,N}}{\dd\mathbb P_{m,m}}
 =\frac{e^{-r\sum_jV(\lambda_j)}}
 {\mathbb E_{m,m}e^{-r\sum_jV(\lambda_j)}}.
\label{eq:fk-general-RN}
\end{equation}
Integration by parts gives
\begin{equation}
\begin{aligned}
 t\,\mathbb E_{m,t}\sum_{j=1}^m\lambda_jV'(\lambda_j)
 &=m+2\mathbb E_{m,t}\sum_{i<j}
 \left(\frac{\lambda_i}{\lambda_i-\lambda_j}
 +\frac{\lambda_j}{\lambda_j-\lambda_i}\right)\\
 &=m+2\binom m2=m^2.
\end{aligned}
\label{eq:fk-general-ibp}
\end{equation}
The growth of $V$ at infinity, \eqref{eq:fk-general-ibp} and Jensen's inequality give
\begin{align*}
 V(x)&\leq C_0+C_1xV'(x),\qquad
 \sup_x e^{-rV(x)}(1+|x|)^{4r}<\infty,\\
 \mathbb E_{m,m}e^{-r\sum_jV(\lambda_j)}
 &\geq e^{-r\mathbb E_{m,m}\sum_jV(\lambda_j)}\geq e^{-Cm},\\
 \mathbb E_{m,N}\prod_{j=1}^m(1+|\lambda_j|)^{4r}
 &=\frac{\mathbb E_{m,m}\prod_{j=1}^m
 [e^{-rV(\lambda_j)}(1+|\lambda_j|)^{4r}]}
 {\mathbb E_{m,m}e^{-r\sum_jV(\lambda_j)}}\leq e^{Cm}.
\end{align*}
Since $\inf_{\mathbb R}V> -\infty$, \eqref{eq:fk-general-RN} also gives,
for every event $E$,
\[
 \mathbb P_{m,N}(E)\leq e^{Cm}\mathbb P_{m,m}(E).
\]
Set $y_0=B-c_0N^{-1/6}$. By \cite[Theorem~1]{DIK} and
Lemma~\ref{prop:finite-comparison} at size $m$,
\begin{align*}
 \kappa m^{2/3}(B-y_0)&\asymp N^{1/2},\\
 \log D_m(y_0)&\leq-cN^{3/2}+CN^{4/3},\\
 \mathbb P_{m,N}\{\lambda_1\leq y\}
 &\leq e^{Cm}D_m(y_0)\leq Ce^{-cN^{3/2}},\qquad y\leq y_0,
\end{align*}
for all sufficiently large $N$.
\end{proof}

\begin{proof}[Proof of \eqref{eq:fixed-k-finite-left}]
Fix $b_0>0$ as in \eqref{eq:fk-4-4}, and define
\begin{equation}
 y_N(\mu):=B+\varepsilon_N\mu,\qquad F_N(\mu):=D_N(y_N(\mu)).
\label{eq:fk-6-1}
\end{equation}
For $x,y\geq\mu$, define
\[
 \widehat L_{N,\mu}(x,y)
 :=\varepsilon_NL_{N,y_N(\mu)}(y_N(x),y_N(y)).
\]
With $\mathbf s=(z_1,\ldots,z_{k-1},\mu+h,\mu)$, the endpoint form of the
ordered Janossy identity \eqref{eq:fixed-k-janossy-density}
\cite[Section~1.1]{SoshnikovCoupled} gives
\begin{equation}
\begin{aligned}
 \mathbb P(\Delta_{N,k}<b)
 &=\int_{\mathbb R}\int_0^b F_N(\mu)
 \int_{z_1>\cdots>z_{k-1}>\mu+h}
 \det[\widehat L_{N,\mu}(s_i,s_j)]_{i,j=1}^{k+1}
 \,\dd\mathbf z\,\dd h\,\dd\mu\\
 &=:\int_{\mathbb R}\int_0^b I_{N,k}(\mu,h)
 \,\dd h\,\dd\mu.
\end{aligned}
\label{eq:fk-6-2}
\end{equation}

Let $p_{j,N}(x):=\pi_j(x;+\infty)/\sqrt{h_j(+\infty)}$ be the
orthonormal polynomials for $e^{-NV(x)}\,\dd x$ on $\mathbb R$, and define
\begin{align*}
 \Phi_N(x)&:=e^{-NV(x)/2}
 (p_{0,N}(x),\ldots,p_{N-1,N}(x))^{\mathsf T},\\
 W_{N,\mu}(x)&:=\varepsilon_N^{1/2}
 \left(\int_{-\infty}^{y_N(\mu)}\Phi_N(t)\Phi_N(t)^*\,\dd t\right)^{-1/2}
 \Phi_N(y_N(x)).
\end{align*}
Then $\widehat L_{N,\mu}(x,y)=\langle W_{N,\mu}(x),W_{N,\mu}(y)\rangle$.
For $0<h\leq b_0$, set
\begin{equation}
 d_{N,\mu,h}:=\frac{W_{N,\mu}(\mu+h)-W_{N,\mu}(\mu)}h
 =\int_0^1W_{N,\mu}'(\mu+th)\,\dd t.
\label{eq:fk-finite-divided-difference}
\end{equation}
For $0<h\leq b_0$ and $(\mathbf z,\mu)\in\mathcal C_k$, set
\[
 (v_1,\ldots,v_{k+1})
 =(W_{N,\mu}(z_1),\ldots,W_{N,\mu}(z_{k-1}),d_{N,\mu,h},W_{N,\mu}(\mu))
\]
and define
\begin{align}
 f_{N,h}(\mathbf z,\mu)
 &:=F_N(\mu)\det[\langle v_i,v_j\rangle]_{i,j=1}^{k+1}
 \mathbf1_{\{z_{k-1}>\mu+h\}},\label{eq:fk-6-8a}\\
 f(\mathbf z,\mu)&:=F_2(\mu)\det\mathcal M_{k,\mu}(\mathbf z).
\label{eq:fk-left-limit-integrand}
\end{align}
For $k=1$, the indicator is omitted.
For $(\mathbf z,\mu)\notin\mathcal C_k$, set $f_{N,h}(\mathbf z,\mu)=f(\mathbf z,\mu)=0$.
Subtracting the last row and column from their $k$th counterparts in
\eqref{eq:fk-6-2} and factoring out $h^2$, we obtain from the definitions
of $f_{N,h}$, $f$ and $q_k$ that
\begin{align}
 \frac{I_{N,k}(\mu,h)}{h^2}
 &=\int_{z_1>\cdots>z_{k-1}>\mu}f_{N,h}(\mathbf z,\mu)\,\dd\mathbf z,
 \label{eq:fk-6-3}\\
 q_k(\mu)&=\int_{z_1>\cdots>z_{k-1}>\mu}f(\mathbf z,\mu)\,\dd\mathbf z.
 \notag
\end{align}
We prove, for every sequence $h_N\downarrow0$,
\begin{equation}
 \|f_{N,h_N}-f\|_{L^1(\mathcal C_k)}\longrightarrow0
 \qquad(N\to\infty),
 \label{eq:fk-6-7}
\end{equation}
and, for some sufficiently large $N_0$,
\begin{equation}
 \sup_{\substack{N\geq N_0\\0<h\leq b_0}}
 \int_{\mathbb R}\frac{I_{N,k}(\mu,h)}{h^2}\,\dd\mu<\infty.
 \label{eq:fk-finite-left-uniform-mass}
\end{equation}

Define
\[
 \widehat K_N(x,y):=\varepsilon_NK_N(y_N(x),y_N(y)).
\]
The analytic edge-kernel expansions in
\cite[(5.2), (5.5)--(5.6)]{LMZ26} give uniform convergence of
$\widehat K_N$ to $K_{\mathrm{Ai}}$ on fixed complex neighborhoods of compact
subsets of $\mathbb R^2$. Applying Cauchy's estimates to the analytic difference
$\widehat K_N-K_{\mathrm{Ai}}$ on smaller complex neighborhoods gives, for every
fixed $j,\ell\geq0$,
\[
 \partial_x^j\partial_y^\ell\widehat K_N(x,y)
 \longrightarrow\partial_x^j\partial_y^\ell K_{\mathrm{Ai}}(x,y)
 \quad\text{locally uniformly on $\mathbb R^2$}.
\]
Fix a sufficiently small $\delta>0$ so that the Airy construction is valid
in $D(B,2\delta)$, and set $r_N:=\delta/\varepsilon_N$.
The Airy sector bounds give, for $0\leq x,y\leq2r_N$,
\begin{equation}
 |\partial_x^j\partial_y^\ell\widehat K_N(x,y)|
 \leq Ce^{-c(x+y)}.
 \label{eq:fk-general-positive-kernel}
\end{equation}
The exterior phase bounds in
\cite[proof of Theorem~2.2, (5.17)--(5.19)]{LMZ26}
(see also \cite{DeiftGioev2007}) give, for every fixed integer $m\geq0$,
\[
 N^m(1+u)^m e^{-N\Phi(y_N(u))/2}
 \leq C e^{-cN}e^{-cu},\qquad u\geq r_N.
\]
In the mixed region $0\leq x\leq r_N$, $y\geq2r_N$, we have
$|x-y|\geq r_N$, so all fixed-order derivatives of $(x-y)^{-1}$ are
uniformly bounded; the same holds with $x,y$ interchanged.
Applying these bounds to the RH representation of $\widehat K_N$
extends \eqref{eq:fk-general-positive-kernel} to all $x,y\geq0$,
uniformly for all sufficiently large $N$.
For the Airy kernel,
\[
 |\partial_x^j\partial_y^\ell K_{\mathrm{Ai}}(x,y)|
 \leq\int_0^\infty
 |\Ai^{(j)}(x+t)\Ai^{(\ell)}(y+t)|\,\dd t
 \leq Ce^{-c(x+y)},\qquad x,y\geq0.
\]
For real $x,t$, the finite Gram representation gives
\[
 |\partial_x^j\widehat K_N(x,t)|^2
 \leq\left.\partial_x^j\partial_y^j\widehat K_N(x,y)\right|_{y=x}
 \widehat K_N(t,t).
\]
For fixed $Q>0$ and an integer $j\geq0$, this bound and its Airy analogue,
compact convergence and \eqref{eq:fk-general-positive-kernel} give
\[
 |\partial_x^j\widehat K_N(x,t)|^2
 +|\partial_x^jK_{\mathrm{Ai}}(x,t)|^2
 \leq C_Qe^{-c_Q(x+t)},\qquad x,t\geq-Q,
\]
for all sufficiently large $N$.

On $L^2((0,\infty))$, define
\[
 (\widehat{\mathbf K}_{N,\mu}g)(\xi)
 :=\int_0^\infty\widehat K_N(\mu+\xi,\mu+\eta)g(\eta)\,\dd\eta.
\]
This operator is a positive contraction, and
$\det(I-\widehat{\mathbf K}_{N,\mu})=F_N(\mu)$. For every fixed $Q\geq0$,
\cite[Corollary~2.3]{LMZ26}, with $t_0=-Q$, and the determinant product give
\begin{align}
 \sup_{|\mu|\leq Q}|F_N(\mu)-F_2(\mu)|&\longrightarrow0,
 \notag\\
 F_N(\mu)&\leq1-\|\widehat{\mathbf K}_{N,\mu}\|,
 \notag\\
 \sup_{\mu\geq-Q}\|(I-\widehat{\mathbf K}_{N,\mu})^{-1}\|
 &\leq F_N(-Q)^{-1}\leq2F_2(-Q)^{-1},
 \label{eq:fk-general-positive-inverse}\\
 \sup_{\mu\geq-Q}\|(I-\mathcal A_\mu^2)^{-1}\|
 &\leq F_2(-Q)^{-1},
 \notag
\end{align}
for all sufficiently large $N$, where $\mathcal A_\mu$ is defined in
\eqref{eq:airy-hankel-operator}.
For every fixed $L>0$, dominated convergence and the resolvent identity,
together with \eqref{eq:fk-general-positive-inverse}, give
\begin{align*}
 &\sup_{\substack{|\mu|\leq Q\\0\leq\xi\leq L}}
 \|\partial_1^j\widehat K_N(\mu+\xi,\mu+\cdot)
 -\partial_1^jK_{\mathrm{Ai}}(\mu+\xi,\mu+\cdot)\|_{L^2((0,\infty))}
 \longrightarrow0,\\
 &\sup_{|\mu|\leq Q}
 \|\widehat{\mathbf K}_{N,\mu}-\mathcal A_\mu^2\|_{\mathrm{HS}}
 \longrightarrow0,\\
 &\sup_{|\mu|\leq Q}
 \|(I-\widehat{\mathbf K}_{N,\mu})^{-1}-(I-\mathcal A_\mu^2)^{-1}\|
 \longrightarrow0.
\end{align*}
For $K=\widehat{\mathbf K}_{N,\mu}$ or $K=\mathcal A_\mu^2$,
\[
 K(I-K)^{-1}=K+K(I-K)^{-1}K.
\]
The preceding limits and the differentiated kernel identity give
\begin{equation}
\begin{aligned}
 &\sup_{\substack{|\mu|\leq Q\\0\leq\xi,\eta\leq L}}
 \Bigl|\partial_1^j\partial_2^\ell\widehat L_{N,\mu}
 (\mu+\xi,\mu+\eta)\\
 &\hspace{28mm}-\partial_1^j\partial_2^\ell L_\mu(\mu+\xi,\mu+\eta)\Bigr|
 \longrightarrow0,
\end{aligned}
 \label{eq:fk-compact-resolvent-convergence}
\end{equation}
for every fixed $j,\ell\geq0$. In \eqref{eq:fk-6-8a}, the entries involving
$d_{N,\mu,h}$ are averages of the first and mixed second derivatives of
$\widehat L_{N,\mu}$. By \eqref{eq:fk-compact-resolvent-convergence} and
the locally uniform convergence $F_N\to F_2$,
\[
 f_{N,h_N}\longrightarrow f
 \quad\text{locally uniformly on $\mathcal C_k$}.
\]

To deduce convergence of the total masses from local convergence, we need
uniform tail bounds. The Gram determinant in \eqref{eq:fk-6-8a} is bounded by
$\|d_{N,\mu,h}\|^2\|W_{N,\mu}(\mu)\|^2$ times its upper-left
$(k-1)\times(k-1)$ minor. Dropping the indicator, integrating this minor
over $z_1>\cdots>z_{k-1}>\mu$, and using
\eqref{eq:fk-finite-divided-difference}, we obtain, for $0<h\leq b_0$,
\begin{equation}
\begin{aligned}
 0\leq\frac{I_{N,k}(\mu,h)}{h^2}
 &\leq F_N(\mu)e_{k-1}(L_{N,y_N(\mu)})\\
 &\quad\times\|W_{N,\mu}(\mu)\|^2
 \sup_{0\leq t\leq b_0}\|W_{N,\mu}'(\mu+t)\|^2.
\end{aligned}
\label{eq:fk-6-4}
\end{equation}
Here $e_{k-1}$ is as in \eqref{eq:elementary-symmetric-sums}, with the
operator acting on $L^2((y_N(\mu),\infty))$.
For the finite-$N$ eigenvalue count $N_I$, the Janossy trace identity is
\[
 F_N(\mu)e_{k-1}(L_{N,y_N(\mu)})
 =\mathbb P\{N_{(y_N(\mu),\infty)}=k-1\}\leq1.
\]
For $|\mu|\leq Q$, \eqref{eq:fk-compact-resolvent-convergence} bounds the
two endpoint factors in \eqref{eq:fk-6-4}. For $k\geq2$, the remaining
$(k-1)$-point Janossy density is bounded by the corresponding correlation
density, so Hadamard's inequality gives
\begin{align*}
 f_{N,h}(\mathbf z,\mu)
 &\leq C_QF_N(\mu)
 \det[\widehat L_{N,\mu}(z_i,z_j)]_{i,j=1}^{k-1}\\
 &\leq C_Q\det[\widehat K_N(z_i,z_j)]_{i,j=1}^{k-1}
 \leq C_Q\prod_{i=1}^{k-1}\widehat K_N(z_i,z_i).
\end{align*}
Thus, for $Z\geq0$ and $|\mu|\leq Q$,
\begin{equation}
\begin{aligned}
 &\int_{\substack{z_1>\cdots>z_{k-1}>\mu\\z_1>Z}}
 f_{N,h}(\mathbf z,\mu)\,\dd\mathbf z\\
 &\qquad\leq\frac{C_Q}{(k-2)!}
 \left(\int_Z^\infty\widehat K_N(x,x)\,\dd x\right)
 \left(\int_{-Q}^\infty\widehat K_N(x,x)\,\dd x\right)^{k-2}\\
 &\qquad\leq C_Qe^{-cZ}.
\end{aligned}
\label{eq:fk-6-10}
\end{equation}
For $\mu\geq0$, \eqref{eq:fk-general-positive-kernel},
\eqref{eq:fk-general-positive-inverse} and the resolvent identity give
\[
 \|W_{N,\mu}(\mu)\|^2
 +\sup_{0\leq t\leq b_0}\|W_{N,\mu}'(\mu+t)\|^2
 \leq Ce^{-c\mu}.
\]
Substituting this bound into \eqref{eq:fk-6-4} yields
\begin{equation}
 \frac{I_{N,k}(\mu,h)}{h^2}\leq Ce^{-c\mu},\qquad\mu\geq0.
\label{eq:fk-6-11}
\end{equation}

By \cite[Theorem~1]{DIK} and Lemma~\ref{prop:finite-comparison},
for some fixed $s_0>0$,
\begin{align*}
 \log F_N(-s)&\leq-cs^3+CN^{-2/3}(1+s)^4,\\
 \sup_{s_0\leq s\leq N^{1/2}}
 \frac{N^{-2/3}(1+s)^4}{s^3}&=O(N^{-1/6}),\\
 \mathbb P\{\zeta_{N,1}\leq-s\}=F_N(-s)&\leq Ce^{-cs^3},
 \qquad s_0\leq s\leq N^{1/2}.
\end{align*}
For each fixed $j\geq1$ and $0<c_j<4j/3$, the finite-$N$ part of
Theorem~\ref{thm:fixed-k-right-tail} gives, along every joint sequence
$N,s\to\infty$ with $s\leq N^{1/2}$,
\begin{align*}
 N\Phi(B+s\varepsilon_N)
 &=\frac43s^{3/2}[1+O(sN^{-2/3})],\\
 e^{c_js^{3/2}}\mathbb P\{\Delta_{N,j}>s\}&\longrightarrow0.
\end{align*}
The sequential criterion therefore yields constants $\mu_j>0$ and
$N_j\geq j+1$ such that
\[
 \mathbb P\{\Delta_{N,j}>s\}\leq Ce^{-c_js^{3/2}},
 \qquad N\geq N_j,\quad\mu_j\leq s\leq N^{1/2}.
\]
For $k\geq2$, the identity
$\zeta_{N,1}-\zeta_{N,k}=\sum_{j=1}^{k-1}\Delta_{N,j}$ gives
\begin{align*}
 \mathbb P\{N_{(y_N(-\mu),\infty)}=k-1\}
 &\leq\mathbb P\{\zeta_{N,k}\leq-\mu\}\\
 &\leq F_N(-\mu/2)+\sum_{j=1}^{k-1}
 \mathbb P\left\{\Delta_{N,j}\geq\frac{\mu}{2(k-1)}\right\}.
\end{align*}
For $k=1$ the count probability equals $F_N(-\mu)$. Hence, for some
$\mu_0=\mu_0(k)>0$ and all sufficiently large $N$,
\begin{equation}
 F_N(-\mu)e_{k-1}(L_{N,y_N(-\mu)})
 \leq Ce^{-c\mu^{3/2}},
 \qquad\mu_0\leq\mu\leq N^{1/2}.
\label{eq:fk-6-14}
\end{equation}
In this range, set $y=B-\mu\varepsilon_N$.
Equations~\eqref{eq:moving-wall-coordinate}--\eqref{eq:muhat-expansion}
and integration of \eqref{eq:K-prefactor} from $\zeta=0$
give, on a fixed complex neighborhood $U$ of $[0,2b_0]$,
\begin{gather*}
 \sup_{\zeta\in U}\left\|\bigl[R_N(y)E_N^{(y)}(y)\bigr]^{-1}
 R_N(y+\varepsilon_N\zeta)E_N^{(y)}(y+\varepsilon_N\zeta)\right\|\leq C,\\
 \sup_{\zeta\in U}|\zeta_{N,\mu}(y+\varepsilon_N\zeta)|\leq C,
 \qquad\mu_N\asymp\mu.
\end{gather*}
The Bessel representation in the proof of Lemma~\ref{lem:wall-Wronskian},
the exponential bounds for $I_0$ and $I_1$, and Cauchy's formula on nested
neighborhoods yield, for the vector $w_{N,\mu}$ in \eqref{eq:right-recovery-vector},
\[
 \max_{0\leq j\leq3}\sup_{0\leq t\leq2b_0}
 \|w_{N,\mu}^{(j)}(t)\|\leq Ce^{C\mu}.
\]
By \eqref{eq:fk-5-3b}, for distinct $s,t\in[0,2b_0]$,
\begin{align*}
 \widehat L_{N,-\mu}(-\mu+s,-\mu+t)
 &=\frac{\det(w_{N,\mu}(s),w_{N,\mu}(t))}{s-t}\\
 &=\int_0^1\det\bigl(w_{N,\mu}'(t+\tau(s-t)),
 w_{N,\mu}(t)\bigr)\,\dd\tau.
\end{align*}
Differentiating the integral gives
\[
 \sup_{\substack{0\leq s,t\leq2b_0\\j+\ell\leq2}}
 |\partial_1^j\partial_2^\ell\widehat L_{N,-\mu}(-\mu+s,-\mu+t)|
 \leq Ce^{C\mu}.
\]
Substituting this derivative bound and \eqref{eq:fk-6-14} into
\eqref{eq:fk-6-4} gives
\begin{equation}
 \frac{I_{N,k}(-\mu,h)}{h^2}
 \leq Ce^{-c\mu^{3/2}+C\mu},
 \qquad\mu_0\leq\mu\leq N^{1/2}.
\label{eq:fk-6-15}
\end{equation}

For all sufficiently large $N$ and $\mu\leq-N^{1/2}$,
set $r=k+1$, $m=N-r$ and
$y=y_N(\mu)\leq B-\kappa^{-1}N^{-1/6}$.
In the original eigenvalue variables, write the joint density of the
$r$ largest eigenvalues at $\mathbf x=(z_1,\ldots,z_{k-1},y+s,y)$ as
$s^2\mathcal B_{N,k}(y,s;\mathbf z)$. Splitting the Vandermonde product
between these eigenvalues and the remaining $m$, and integrating out the
latter, gives
\begin{equation}
\begin{aligned}
 \mathcal B_{N,k}
 ={}&\frac{Z_{m,N}}{Z_{N,N}}e^{-N\sum_{i=1}^rV(x_i)}
 \frac{\Delta(\mathbf x)^2}{s^2}\\
 &\times\mathbb E_{m,N}\left[
 \mathbf1_{\{\lambda_1<y\}}
 \prod_{j=1}^m\prod_{i=1}^r(x_i-\lambda_j)^2\right].
\end{aligned}
\label{eq:fk-6-17}
\end{equation}
Cauchy--Schwarz, $|x-\lambda|\leq(1+|x|)(1+|\lambda|)$ and
Lemma~\ref{lem:fixed-k-partition-coercivity} yield
\begin{align*}
 &\mathbb E_{m,N}\left[
 \mathbf1_{\{\lambda_1<y\}}\prod_{j=1}^m\prod_{i=1}^r(x_i-\lambda_j)^2\right]\\
 &\quad\leq\mathbb P_{m,N}(\lambda_1<y)^{1/2}
 \left(\mathbb E_{m,N}\prod_{j=1}^m(1+|\lambda_j|)^{4r}\right)^{1/2}
 \prod_{i=1}^r(1+|x_i|)^{2m}\\
 &\quad\leq e^{-cN^{3/2}+CN}\prod_{i=1}^r(1+|x_i|)^{2m}.
\end{align*}
Since
\[
 \frac{\Delta(\mathbf x)^2}{s^2}
 =\Delta(\mathbf z)^2\prod_{i=1}^{k-1}(z_i-y-s)^2(z_i-y)^2,
 \qquad V(x)\geq c|x|^{\deg V}-C,
\]
the expectation bound and \eqref{eq:fk-general-partition-ratio} give,
uniformly for $0\leq s\leq b_0\varepsilon_N$,
\begin{equation}
 \mathcal B_{N,k}
 \leq e^{-cN^{3/2}+CN}
 \exp\left[-cN\sum_{i=1}^r|x_i|^{\deg V}\right].
\label{eq:fk-6-20}
\end{equation}
Changing to the original eigenvalue variables in \eqref{eq:fk-6-2} gives
\begin{align*}
 \frac{I_{N,k}(\mu,h)}{h^2}
 &=\varepsilon_N^4
 \int_{z_1>\cdots>z_{k-1}>y_N(\mu)+\varepsilon_Nh}
 \mathcal B_{N,k}(y_N(\mu),\varepsilon_Nh;\mathbf z)\,\dd\mathbf z,\\
 \dd\mu&=\varepsilon_N^{-1}\dd y.
\end{align*}
Consequently,
\begin{equation}
\begin{aligned}
 \int_{\mu\leq-N^{1/2}}\frac{I_{N,k}(\mu,h)}{h^2}\,\dd\mu
 &\leq\varepsilon_N^3e^{-cN^{3/2}+CN}
 \int_{\mathbb R^k}e^{-cN(|y|^{\deg V}+\sum_{i=1}^{k-1}|z_i|^{\deg V})}
 \,\dd\mathbf z\,\dd y\\
 &\leq e^{-cN^{3/2}+CN}.
\end{aligned}
\label{eq:fk-6-22}
\end{equation}

For $Q\geq\mu_0$, \eqref{eq:fk-6-11}, \eqref{eq:fk-6-15} and
\eqref{eq:fk-6-22} imply
\begin{align*}
 &\limsup_{N\to\infty}\sup_{0<h\leq b_0}
 \int_{|\mu|>Q}\frac{I_{N,k}(\mu,h)}{h^2}\,\dd\mu\\
 &\qquad\leq C\int_Q^\infty
 \bigl(e^{-c\mu}+e^{-c\mu^{3/2}+C\mu}\bigr)\,\dd\mu
 \longrightarrow0\qquad(Q\to\infty).
\end{align*}
The uniform bound for $f_{N,h}$ allows us to integrate the local limit over
bounded subsets of $\mathcal C_k$ by dominated convergence.
The preceding $\mu$-tail estimate and \eqref{eq:fk-6-10} for $k\geq2$
extend this convergence to the whole chamber:
\[
 \int_{\mathcal C_k}f_{N,h_N}\longrightarrow
 \int_{\mathcal C_k}f=\int_{\mathbb R}q_k(\mu)\,\dd\mu.
\]
Scheff\'e's lemma now proves \eqref{eq:fk-6-7}. The same tail bounds and
\eqref{eq:fk-6-4} on $|\mu|\leq Q$ give
\eqref{eq:fk-finite-left-uniform-mass}. Finally, for every $b_N\downarrow0$,
\begin{align*}
 \frac{\mathbb P(\Delta_{N,k}<b_N)}{b_N^3}
 &=\int_0^1t^2\int_{\mathbb R}
 \frac{I_{N,k}(\mu,b_Nt)}{(b_Nt)^2}\,\dd\mu\,\dd t\\
 &\longrightarrow\int_0^1t^2\,\dd t
 \int_{\mathbb R}q_k(\mu)\,\dd\mu=\frac{k}{6},
\end{align*}
by dominated convergence, using \eqref{eq:fk-6-3}, \eqref{eq:fk-6-7},
\eqref{eq:fk-finite-left-uniform-mass} and
\eqref{eq:def-main-left-coefficient}.
\end{proof}

\subsection{Proof of Corollary~\ref{cor:fixed-k-left-collision-law}}

\begin{proof}
For $b>0$, let $\rho_{k,b}$ be the conditional density on
$\mathcal C_k\times(0,1)$ of
\[
 \left(\zeta_1,\ldots,\zeta_{k-1},\zeta_{k+1},
 \frac{\Delta_{\mathrm{Ai},k}}{b}\right)
 \quad\text{given }\Delta_{\mathrm{Ai},k}<b.
\]
With $\mathbf s=(z_1,\ldots,z_{k-1},\mu+bt,\mu)$,
\eqref{eq:fk-4-1} gives
\[
 \rho_{k,b}(\mathbf z,\mu,t)
 =\frac{bF_2(\mu)}{\Prob(\Delta_{\mathrm{Ai},k}<b)}
 \det[L_\mu(s_i,s_j)]_{i,j=1}^{k+1}
 \mathbf1_{\{z_1>\cdots>z_{k-1}>\mu+bt\}},
\]
with the indicator omitted when $k=1$. Taking $\varepsilon=bt$ in the
Gram coalescence following \eqref{eq:fk-4-5} and using
\eqref{eq:fixed-k-airy-left}, we obtain, with $f$ from
\eqref{eq:fk-left-limit-integrand},
\begin{equation}
 \rho_{k,b}(\mathbf z,\mu,t)\longrightarrow
 \rho_k(\mathbf z,\mu,t):=\frac6k t^2f(\mathbf z,\mu).
\label{eq:airy-left-marked-limit-density}
\end{equation}
By \eqref{eq:collision-measure},
\[
 \rho_k(\mathbf z,\mu,t)\,\mathrm d\mathbf z\,\mathrm d\mu\,\mathrm dt
 =\mathrm d\nu_k(\mathbf z,\mu)\,3t^2\mathrm dt.
\]
For every sequence $b_N\downarrow0$ as $N\to\infty$, let
$\rho_{N,k,b_N}$ denote the joint density of
\[
 \left(\zeta_{N,1},\ldots,\zeta_{N,k-1},\zeta_{N,k+1},
 \frac{\Delta_{N,k}}{b_N}\right)
\]
conditional on $\{\Delta_{N,k}<b_N\}$.
Equation~\eqref{eq:fk-6-8a}, the local convergence following
\eqref{eq:fk-compact-resolvent-convergence} with $h_N=b_Nt$,
and \eqref{eq:fixed-k-finite-left} give, pointwise on
$\mathcal C_k\times(0,1)$,
\[
 \rho_{N,k,b_N}(\mathbf z,\mu,t)
 =\frac{t^2f_{N,b_Nt}(\mathbf z,\mu)}
 {\Prob(\Delta_{N,k}<b_N)/b_N^3}
 \longrightarrow\rho_k(\mathbf z,\mu,t).
\]
Scheff\'e's lemma gives, in $L^1(\mathcal C_k\times(0,1))$,
\[
 \rho_{k,b}\longrightarrow\rho_k\quad(b\downarrow0),
 \qquad\rho_{N,k,b_N}\longrightarrow\rho_k\quad(N\to\infty).
\]

Extend the densities by zero outside $\mathcal C_k\times(0,1)$.
Replacing $\zeta_{k+1}$ by $\zeta_k$, and similarly at finite $N$,
gives the conditional densities
\[
 \rho_{k,b}(\mathbf z,\mu-bt,t),\qquad
 \rho_{N,k,b_N}(\mathbf z,\mu-b_Nt,t).
\]
By translation invariance of the $L^1$ norm, the preceding convergence
and $L^1$-continuity of translations of $f$, uniformly for shifts of
size at most $b$ or $b_N$, give convergence of these densities to
$\rho_k$ in $L^1(\mathbb R^k\times(0,1))$. This proves
\eqref{eq:airy-left-collision-law} and \eqref{eq:finite-left-collision-law}.
\end{proof}

\clearpage
\appendix
	\section{The Airy parametrix}\label{airy}
	The Airy parametrix $\Phi^{({\Ai})}$ is the unique solution of the following RH problem.
	\begin{figure}[H]
		\begin{center}
			
			\tikzset{every picture/.style={line width=0.75pt}} 
			
			\begin{tikzpicture}[x=0.75pt,y=0.75pt,yscale=-1,xscale=1]
				
				\draw    (52,120) -- (203,120) -- (333,120) ;
				\draw [shift={(132.5,120)}, rotate = 180] [fill={rgb, 255:red, 0; green, 0; blue, 0 }] (7.14,-3.43) -- (0,0) -- (7.14,3.43) -- (4.74,0) -- cycle ;  
				\draw [shift={(273,120)}, rotate = 180] [fill={rgb, 255:red, 0; green, 0; blue, 0 }  ](7.14,-3.43) -- (0,0) -- (7.14,3.43) -- (4.74,0) -- cycle    ;
				\draw    (116.5,20) -- (216.5,120) ;
				\draw [shift={(170.04,73.54)}, rotate = 225] [fill={rgb, 255:red, 0; green, 0; blue, 0 }  ](7.14,-3.43) -- (0,0) -- (7.14,3.43) -- (4.74,0) -- cycle    ;
				\draw    (117.5,208) -- (216.5,120) ;
				\draw [shift={(170.74,160.68)}, rotate = 138.37] [fill={rgb, 255:red, 0; green, 0; blue, 0 }  ](7.14,-3.43) -- (0,0) -- (7.14,3.43) -- (4.74,0) -- cycle    ;
				\draw  [fill={rgb, 255:red, 0; green, 0; blue, 0 }  ,fill opacity=1 ] (216.67,119.83) .. controls (216.67,119.19) and (216.15,118.66) .. (215.5,118.66) .. controls (214.85,118.66) and (214.33,119.19) .. (214.33,119.83) .. controls (214.33,120.48) and (214.85,121) .. (215.5,121) .. controls (216.15,121) and (216.67,120.48) .. (216.67,119.83) -- cycle ;
				
				\draw (218.5,123) node [anchor=north west][inner sep=0.75pt]   [align=left] {$0$};
				\draw (343,117) node [anchor=north west][inner sep=0.75pt]   [align=left] {$\Sigma_1$};
				\draw (98,14) node [anchor=north west][inner sep=0.75pt]   [align=left] {$\Sigma_2$};
				\draw (33,113) node [anchor=north west][inner sep=0.75pt]   [align=left] {$\Sigma_3$};
				\draw (95,205) node [anchor=north west][inner sep=0.75pt]   [align=left] {$\Sigma_4$};

			\end{tikzpicture}
			\caption{The jump contours of the RH problem for $\Phi^{({\Ai})}$.}
			\label{fig:Airy}
		\end{center}
	\end{figure}
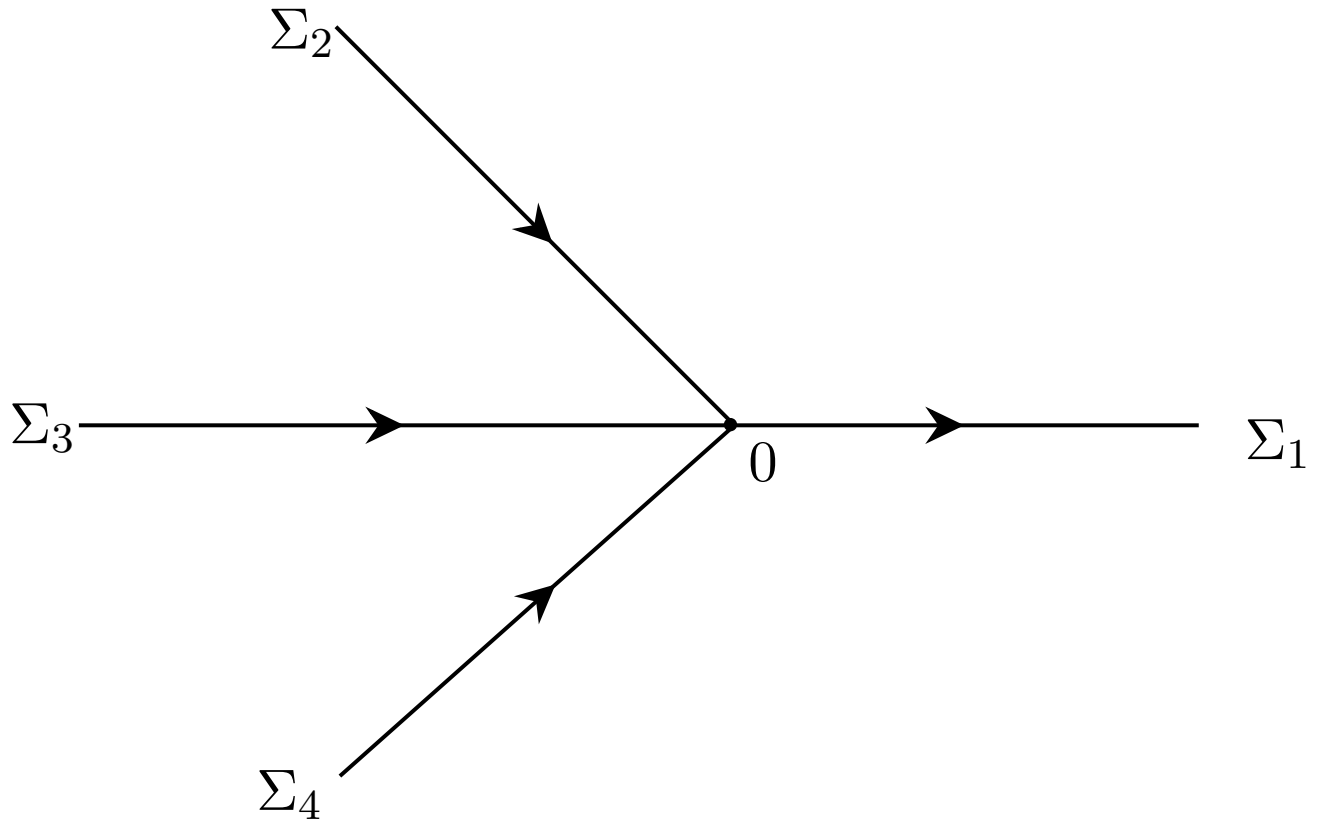
	\begin{rhp}
		\hfill
		\begin{itemize}
			\item[\rm(a)] $\Phi^{(\mathrm{Ai})}(z)$ is analytic in $\mathbb{C} \setminus \{\cup_{j=1}^4 \Sigma_j \cup \{0\}\}$, where the contours $\Sigma_j$, $j=1, 2, 3, 4$, are shown in Figure~\ref{fig:Airy}.
			\item[\rm(b)] For $z \in \cup_{j=1}^{4} \Sigma_j$,
			\begin{align}\label{eq:airy-jump-matrices}
				\Phi^{({\Ai})}_+(z)=\Phi^{({\Ai})}_-(z)\begin{cases}
					\begin{pmatrix}
						1 & 1\\
						0 & 1
					\end{pmatrix}, & \qquad z \in \Sigma_1,\\
					\begin{pmatrix}
						1 & 0\\
						1& 1
					\end{pmatrix}, & \qquad z \in \Sigma_2 \cup \Sigma_4,\\
					\begin{pmatrix}
						0 & 1\\
						-1 & 0
					\end{pmatrix}, & \qquad z \in \Sigma_3.
				\end{cases}
			\end{align}
			\item[\rm(c)] As $z \to \infty$,
			\begin{align}\label{infty:Ai}
				\Phi^{({\Ai})}(z) = \frac{1}{\sqrt{2}} \begin{pmatrix}
					z^{-\frac 14} & 0\\
					0 & z^{\frac 14}
				\end{pmatrix} \begin{pmatrix}
					1 & \ii\\
					\ii & 1
				\end{pmatrix}\left(I + \Boh(z^{-\frac 32})\right)e^{-\frac 23 z^{3/2} \sigma_3}.
			\end{align}
			\item[\rm(d)] $\Phi^{({\Ai})}(z)$ is bounded near the origin.
		\end{itemize}
	\end{rhp}

	Let $\omega:= e^{ 2 \pi \ii/3}$. The solution is given by (cf. \cite{Deift1999})
	\begin{align*}
		\Phi^{({\Ai})}(z) = \sqrt{2 \pi} \begin{cases}
			\begin{pmatrix}
				\Ai (z) & - \omega^2 \Ai (\omega^2 z)\\
				-\ii \Ai'(z) & \ii \omega \Ai' (\omega^2 z)
			\end{pmatrix}, & \arg z \in \left(0, \frac{3 \pi}{4} \right),\\
			\begin{pmatrix}
				-\omega \Ai (\omega z) & - \omega^2 \Ai (\omega^2 z)\\
				\ii \omega^2 \Ai'(\omega z) & \ii \omega \Ai' (\omega^2 z)
			\end{pmatrix}, & \arg z \in \left(\frac{3 \pi}{4}, \pi \right),\\
			\begin{pmatrix}
				-\omega^2 \Ai (\omega^2 z) & \omega \Ai (\omega z)\\
				\ii \omega \Ai'(\omega^2 z) & -\ii \omega^2 \Ai' (\omega z)
			\end{pmatrix}, & \arg z \in \left(-\pi, -\frac{3 \pi}{4} \right),\\
			\begin{pmatrix}
				\Ai (z) &  \omega \Ai (\omega z)\\
				-\ii \Ai'(z) & -\ii \omega^2 \Ai' (\omega z)
			\end{pmatrix}, & \arg z \in \left(-\frac{3 \pi}{4}, 0\right).
		\end{cases}
	\end{align*}
{}

\clearpage
\section{The Bessel parametrix}
\label{A}
The Bessel parametrix of order zero, $\Phi^{(\Bes)}$, is the unique
solution of the following RH problem.
		\begin{figure}[H]
			\centering

\tikzset{every picture/.style={line width=0.5pt}} 

\begin{tikzpicture}[x=0.75pt,y=0.75pt,yscale=-1,xscale=1]

\draw    (100,130) -- (269.5,131) ;
\draw [shift={(188.55,130.52)}, rotate = 180.34] [fill={rgb, 255:red, 0; green, 0; blue, 0 }  ][line width=0.08]  [draw opacity=0] (7.14,-3.43) -- (0,0) -- (7.14,3.43) -- cycle    ;
\draw   (169.12,29.61) -- (269.5,131) ;
\draw [shift={(221.98,83.01)}, rotate = 225.29] [fill={rgb, 255:red, 0; green, 0; blue, 0 }  ][line width=0.08]  [draw opacity=0] (7.14,-3.43) -- (0,0) -- (7.14,3.43) -- cycle    ;
\draw     (269.5,131) -- (171.12,231.61) ;
\draw [shift={(224.01,177.52)}, rotate = 134.36] [fill={rgb, 255:red, 0; green, 0; blue, 0 }  ][line width=0.08]  [draw opacity=0] (7.14,-3.43) -- (0,0) -- (7.14,3.43) -- cycle    ;
\draw  [fill={rgb, 255:red, 0; green, 0; blue, 0 }  ,fill opacity=1 ][line width=0.75]  (267.5,131.5) .. controls (267.5,130.67) and (268.17,130) .. (269,130) .. controls (269.83,130) and (270.5,130.67) .. (270.5,131.5) .. controls (270.5,132.33) and (269.83,133) .. (269,133) .. controls (268.17,133) and (267.5,132.33) .. (267.5,131.5) -- cycle ;

\draw (281,80) node [anchor=north west][inner sep=0.75pt]   [align=left] {I};
\draw (271.5,133) node [anchor=north west][inner sep=0.75pt]   [align=left] {0};
\draw (145,78) node [anchor=north west][inner sep=0.75pt]   [align=left] {II};
\draw (144,156) node [anchor=north west][inner sep=0.75pt]   [align=left] {III};
\draw (175,13) node [anchor=north west][inner sep=0.75pt]   [align=left] {$\Gamma_1$};
\draw (80,122) node [anchor=north west][inner sep=0.75pt]   [align=left] {$\Gamma_2$};
\draw (181,232) node [anchor=north west][inner sep=0.75pt]   [align=left] {$\Gamma_3$};

\end{tikzpicture}
			\caption{The jump contours $\Gamma_j$, $j=1,2,3$, and the sectors I--III in the RH problem for $\Phi^{(\Bes)}$.}
			\label{figure-Bessel}
		\end{figure}
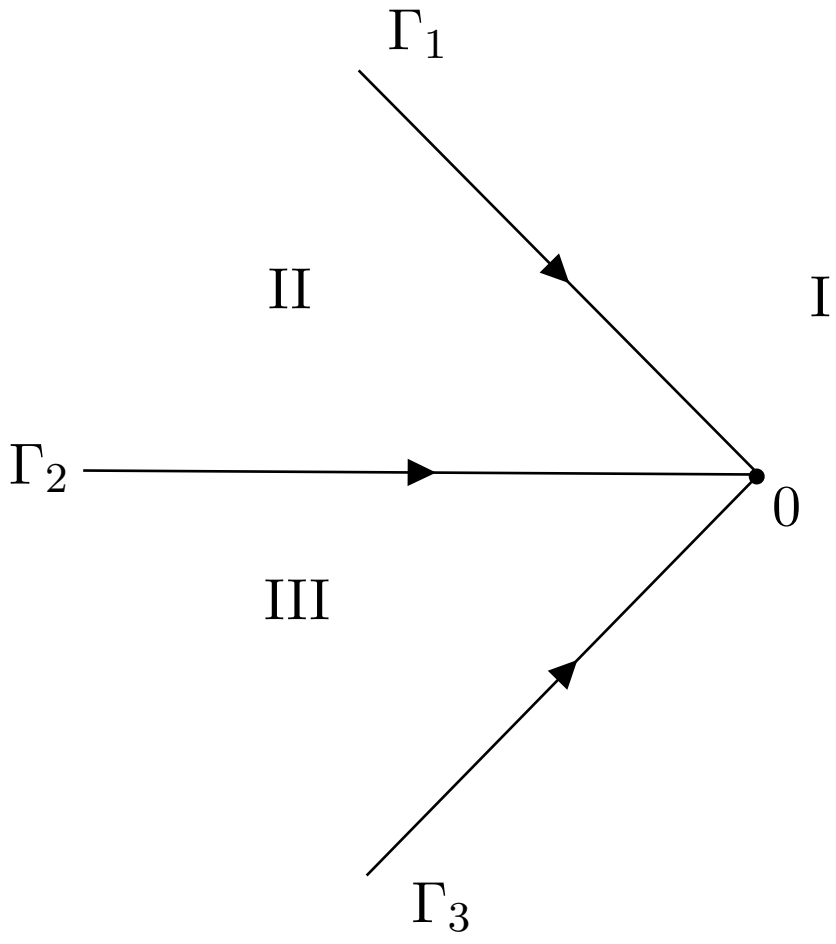

\begin{rhp}\label{rhp:Bessel}
\hfill
\begin{itemize}
 \item[\rm(a)] $\Phi^{(\Bes)}$ is analytic in
 $\mathbb C\setminus(\bigcup_{j=1}^3\Gamma_j\cup\{0\})$, where the contours
 $\Gamma_j$ and sectors I--III are shown in Figure~\ref{figure-Bessel}.
 \item[\rm(b)] For $z\in\bigcup_{j=1}^3\Gamma_j$,
 \begin{align*}
  \Phi_+^{(\Bes)}(z)=\Phi_-^{(\Bes)}(z)\begin{cases}
   \begin{pmatrix}
    1&0\\
    1&1
   \end{pmatrix},&\qquad z\in\Gamma_1\cup\Gamma_3,\\
   \begin{pmatrix}
    0&1\\
    -1&0
   \end{pmatrix},&\qquad z\in\Gamma_2.
  \end{cases}
 \end{align*}
 \item[\rm(c)] As $z\to\infty$,
 \begin{align}\label{infty:Bes}
  \Phi^{(\Bes)}(z)
  =\frac{(\pi^2z)^{-\sigma_3/4}}{\sqrt2}
  \begin{pmatrix}1&\ii\\\ii&1\end{pmatrix}
  \left(I+\frac1{8z^{1/2}}
  \begin{pmatrix}-1&-2\ii\\-2\ii&1\end{pmatrix}
  +\Boh(z^{-1})\right)e^{z^{1/2}\sigma_3}.
 \end{align}
 \item[\rm(d)] As $z\to0$,
 \begin{equation*}
  \Phi^{(\Bes)}(z)=\Boh(\log|z|).
 \end{equation*}
\end{itemize}
\end{rhp}
Let $I_0$ and $K_0$ denote the modified Bessel functions of order zero
\cite[Chapter~10]{DLMF}.  By \cite[Section~6]{KMVV04}, the solution is
\begin{equation}\label{def:Bespara}
 \Phi^{(\Bes)}(z)=
 \begin{pmatrix}
  I_0(z^{1/2})&\dfrac{\ii}{\pi}K_0(z^{1/2})\\
  \pi\ii z^{1/2}I_0'(z^{1/2})&-z^{1/2}K_0'(z^{1/2})
 \end{pmatrix}
 \begin{cases}
  I,&z\in\mathrm I,\\
  \begin{pmatrix}1&0\\-1&1\end{pmatrix},&z\in\mathrm{II},\\
  \begin{pmatrix}1&0\\1&1\end{pmatrix},&z\in\mathrm{III}.
 \end{cases}
\end{equation}

\clearpage

\section*{Statements and Declarations}

\paragraph{Competing interests}
The author declares no competing interests.

\paragraph{Data availability}
No datasets were generated or analysed during this study.


\begin{thebibliography}{99}

\bibitem{Andreief1886}
C. Andr\'eief,
\emph{Note sur une relation entre les int\'egrales d\'efinies des produits
des fonctions},
M\'emoires de la Soci\'et\'e des Sciences Physiques et Naturelles de Bordeaux
(3) \textbf{2} (1886), 1--14.

\bibitem{Baik2006}
J. Baik,
\href{https://doi.org/10.1215/S0012-7094-06-13321-5}
{\emph{Painlev\'e formulas of the limiting distributions for nonnull complex
sample covariance matrices}},
Duke Mathematical Journal \textbf{133} (2006), no.~2, 205--235.

\bibitem{BBP05}
J. Baik, G. Ben Arous, and S. P\'ech\'e,
\href{https://doi.org/10.1214/009117905000000233}
{\emph{Phase transition of the largest eigenvalue for nonnull complex sample
covariance matrices}},
The Annals of Probability \textbf{33} (2005), no.~5, 1643--1697.

\bibitem{BBD}
J. Baik, R. Buckingham, and J. DiFranco,
\href{https://doi.org/10.1007/s00220-008-0433-5}
{\emph{Asymptotics of Tracy--Widom distributions and the total integral of a
Painlev\'e II function}},
Communications in Mathematical Physics \textbf{280} (2008), 463--497.

\bibitem{BR00}
J. Baik and E. M. Rains,
\href{https://doi.org/10.1023/A:1018615306992}
{\emph{Limiting distributions for a polynuclear growth model with external
sources}},
Journal of Statistical Physics \textbf{100} (2000), nos.~3--4, 523--541.

\bibitem{BCEHK25}
Z. Bao, G. Cipolloni, L. Erd\H{o}s, J. Henheik, and O. Kolupaiev,
\emph{Law of fractional logarithm for random matrices},
\href{https://arxiv.org/abs/2503.18922v2}{arXiv:2503.18922v2} (2025).

\bibitem{Baslingker26}
J. Baslingker,
\emph{Optimal upper tail estimates for the edge eigenvalues of
$\beta$-ensembles},
\href{https://arxiv.org/abs/2610.11239v2}{arXiv:2610.11239v2} (2026).

\bibitem{BBBK25}
J. Baslingker, R. Basu, S. Bhattacharjee, and M. Krishnapur,
\href{https://doi.org/10.2140/pmp.2025.6.1379}
{\emph{Optimal tail estimates in $\beta$-ensembles and applications to last
passage percolation}},
Probability and Mathematical Physics \textbf{6} (2025), no.~4, 1379--1442.

\bibitem{BS}
A. Borodin and A. Soshnikov,
\href{https://doi.org/10.1023/A:1026025003309}
{\emph{Janossy densities. I. Determinantal ensembles}},
Journal of Statistical Physics \textbf{113} (2003), 595--610.

\bibitem{ClaeysDoeraene}
T. Claeys and A. Doeraene,
\href{https://doi.org/10.1111/sapm.12209}
{\emph{The generating function for the Airy point process and a system of
coupled Painlev\'e II equations}},
Studies in Applied Mathematics \textbf{140} (2018), no.~4, 403--437.

\bibitem{CorwinHammond14}
I. Corwin and A. Hammond,
\href{https://doi.org/10.1007/s00222-013-0462-3}
{\emph{Brownian Gibbs property for Airy line ensembles}},
Inventiones mathematicae \textbf{195} (2014), no.~2, 441--508.

\bibitem{CorwinSun14}
I. Corwin and X. Sun,
\href{https://doi.org/10.1214/ECP.v19-3504}
{\emph{Ergodicity of the Airy line ensemble}},
Electronic Communications in Probability \textbf{19} (2014), no.~49, 1--11.

\bibitem{DVJ1988}
D. J. Daley and D. Vere-Jones,
\emph{An Introduction to the Theory of Point Processes},
New York: Springer, 1988.

\bibitem{deBoor2005}
C. de Boor,
\href{https://surveys-in-approximation-theory.com/papers/2/2.pdf}
{\emph{Divided differences}},
Surveys in Approximation Theory \textbf{1} (2005), 46--69.

\bibitem{Deift1999}
P. Deift,
\emph{Orthogonal Polynomials and Random Matrices: A Riemann--Hilbert Approach},
Courant Lecture Notes in Mathematics, vol.~3,
Courant Institute of Mathematical Sciences, New York, and
American Mathematical Society, Providence, RI, 1999.

\bibitem{DeiftGioev2007}
P. Deift and D. Gioev,
\href{https://doi.org/10.1002/cpa.20164}
{\emph{Universality at the edge of the spectrum for unitary, orthogonal,
and symplectic ensembles of random matrices}},
Communications on Pure and Applied Mathematics \textbf{60} (2007), no.~6,
867--910.

\bibitem{DIK}
P. Deift, A. Its, and I. Krasovsky,
\href{https://doi.org/10.1007/s00220-007-0409-x}
{\emph{Asymptotics of the Airy-kernel determinant}},
Communications in Mathematical Physics \textbf{278} (2008), 643--678.

\bibitem{DKMVZ99}
P. Deift, T. Kriecherbauer, K. T.-R. McLaughlin, S. Venakides, and X. Zhou,
\href{https://doi.org/10.1002/(SICI)1097-0312(199911)52:11<1335::AID-CPA1>3.0.CO;2-1}
{\emph{Uniform asymptotics for polynomials orthogonal with respect to varying
exponential weights and applications to universality questions in random
matrix theory}},
Communications on Pure and Applied Mathematics \textbf{52} (1999), no.~11,
1335--1425.

\bibitem{DeiftTrogdon}
P. Deift and T. Trogdon,
\href{https://doi.org/10.1002/cpa.21715}
{\emph{Universality for the Toda algorithm to compute the largest eigenvalue
of a random matrix}},
Communications on Pure and Applied Mathematics \textbf{71} (2018), no.~3,
505--536.

\bibitem{DZ93}
P. Deift and X. Zhou,
\href{https://doi.org/10.2307/2946540}
{\emph{A steepest descent method for oscillatory Riemann--Hilbert problems.
Asymptotics for the MKdV equation}},
Annals of Mathematics (2) \textbf{137} (1993), 295--368.

\bibitem{DV13}
L. Dumaz and B. Vir\'ag,
\href{https://doi.org/10.1214/11-AIHP475}
{\emph{The right tail exponent of the Tracy--Widom $\beta$ distribution}},
Annales de l'Institut Henri Poincar\'e, Probabilit\'es et Statistiques
\textbf{49} (2013), no.~4, 915--933.

\bibitem{EKS16}
P. Eichelsbacher, T. Kriecherbauer, and K. Sch\"uler,
\href{https://doi.org/10.3842/SIGMA.2016.093}
{\emph{Precise deviations results for the maxima of some determinantal
point processes: the upper tail}},
Symmetry, Integrability and Geometry: Methods and Applications
\textbf{12} (2016), Paper No.~093, 18~pp.

\bibitem{EX23}
L. Erd\H{o}s and Y. Xu,
\href{https://doi.org/10.3150/22-BEJ1490}
{\emph{Small deviation estimates for the largest eigenvalue of Wigner matrices}},
Bernoulli \textbf{29} (2023), no.~2, 1063--1079.

\bibitem{FIK92}
A. S. Fokas, A. R. Its and A. V. Kitaev,
\href{https://doi.org/10.1007/BF02096594}
{\emph{The isomonodromy approach to matrix models in 2d quantum gravity}},
Comm. Math. Phys. \textbf{147} (1992), 395--430.

\bibitem{Forrester93}
P. J. Forrester,
\href{https://doi.org/10.1016/0550-3213(93)90126-A}
{\emph{The spectrum edge of random matrix ensembles}},
Nuclear Physics B \textbf{402} (1993), no.~3, 709--728.

\bibitem{JMU81}
M. Jimbo, T. Miwa, and K. Ueno,
\href{https://doi.org/10.1016/0167-2789(81)90013-0}
{\emph{Monodromy preserving deformation of linear ordinary differential
equations with rational coefficients}},
Physica D: Nonlinear Phenomena \textbf{2} (1981), no.~2, 306--352.

\bibitem{Krattenthaler1999}
C. Krattenthaler,
\href{https://www.mat.univie.ac.at/~slc/wpapers/s42kratt.pdf}
{\emph{Advanced determinant calculus}},
S\'eminaire Lotharingien de Combinatoire \textbf{42} (1999), Article B42q, 67~pp.

\bibitem{KMVV04}
A. B. J. Kuijlaars, K. T.-R. McLaughlin, W. Van Assche, and M. Vanlessen,
\href{https://doi.org/10.1016/j.aim.2003.08.015}
{\emph{The Riemann--Hilbert approach to strong asymptotics for orthogonal
polynomials on $[-1,1]$}},
Advances in Mathematics \textbf{188} (2004), 337--398.

\bibitem{Ledoux07}
M. Ledoux,
\href{https://doi.org/10.1007/978-3-540-72053-9_10}
{\emph{Deviation inequalities on largest eigenvalues}},
in \emph{Geometric Aspects of Functional Analysis},
Lecture Notes in Mathematics, vol.~1910,
Springer, Berlin, 2007, pp.~167--219.

\bibitem{LR10}
M. Ledoux and B. Rider,
\href{https://doi.org/10.1214/EJP.v15-798}
{\emph{Small deviations for beta ensembles}},
Electronic Journal of Probability \textbf{15} (2010), no.~41, 1319--1343.

\bibitem{LMZ26}
J. Liu, Q. Mao, and L. Zhang,
\emph{Edgeworth expansions of extreme eigenvalue distributions for random
unitary ensembles},
\href{https://arxiv.org/abs/2609.39169}{arXiv:2609.39169} (2026).

\bibitem{Nishigaki2021}
S. M. Nishigaki,
\href{https://doi.org/10.1093/ptep/ptab123}
{\emph{Tracy--Widom method for J\'anossy density and joint distribution of
extremal eigenvalues of random matrices}},
Progress of Theoretical and Experimental Physics \textbf{2021} (2021),
no.~11, Art.~ID 113A01.

\bibitem{DLMF}
F. W. J. Olver, A. B. Olde Daalhuis, D. W. Lozier, B. I. Schneider,
R. F. Boisvert, C. W. Clark, B. R. Miller, B. V. Saunders, H. S. Cohl,
and M. A. McClain, eds.,
\emph{NIST Digital Library of Mathematical Functions},
\url{https://dlmf.nist.gov/}, Release 1.2.8 of 2026-09-15.

\bibitem{PZ17}
E. Paquette and O. Zeitouni,
\href{https://doi.org/10.1214/16-AOP1161}
{\emph{Extremal eigenvalue correlations in the GUE minor process and a law
of fractional logarithm}},
The Annals of Probability \textbf{45} (2017), no.~6A, 4112--4166.

\bibitem{PS03}
L. Pastur and M. Shcherbina,
\href{https://www.mathnet.ru/eng/jmag255}
{\emph{On the edge universality of the local eigenvalue statistics of matrix models}},
Matematicheskaya Fizika, Analiz, Geometriya \textbf{10} (2003), no.~3, 335--365.

\bibitem{PerretSchehr}
A. Perret and G. Schehr,
\href{https://doi.org/10.1007/s10955-014-1044-5}
{\emph{Near-extreme eigenvalues and the first gap of Hermitian random matrices}},
Journal of Statistical Physics \textbf{156} (2014), 843--876.

\bibitem{RRV11}
J. A. Ram\'{\i}rez, B. Rider, and B. Vir\'ag,
\href{https://doi.org/10.1090/S0894-0347-2011-00703-0}
{\emph{Beta ensembles, stochastic Airy spectrum, and a diffusion}},
Journal of the American Mathematical Society \textbf{24} (2011), no.~4,
919--944.

\bibitem{RZ}
B. Rider and X. Zhou,
\href{https://doi.org/10.1093/imrn/rnn037}
{\emph{Janossy densities for unitary ensembles at the spectral edge}},
International Mathematics Research Notices \textbf{2008} (2008), Art.~ID rnn037, 51~pp.

\bibitem{Scheffe1947}
H. Scheff\'e,
\href{https://doi.org/10.1214/aoms/1177730390}
{\emph{A useful convergence theorem for probability distributions}},
The Annals of Mathematical Statistics \textbf{18} (1947), no.~3, 434--438.

\bibitem{Simon1981}
B. Simon,
\href{https://doi.org/10.1090/S0002-9939-1981-0619977-2}
{\emph{Convergence in trace ideals}},
Proceedings of the American Mathematical Society \textbf{83} (1981), 39--43.

\bibitem{SoshnikovCoupled}
A. Soshnikov,
\href{https://doi.org/10.1007/s00220-004-1177-5}
{\emph{Janossy densities of coupled random matrices}},
Communications in Mathematical Physics \textbf{251} (2004), 447--471.

\bibitem{TW93}
C. A. Tracy and H. Widom,
\href{https://doi.org/10.1016/0370-2693(93)91114-3}
{\emph{Level-spacing distributions and the Airy kernel}},
Physics Letters B \textbf{305} (1993), 115--118.

\bibitem{TW94}
C. A. Tracy and H. Widom,
\href{https://doi.org/10.1007/BF02100489}
{\emph{Level-spacing distributions and the Airy kernel}},
Communications in Mathematical Physics \textbf{159} (1994), 151--174.

\bibitem{WBF}
N. S. Witte, F. Bornemann, and P. J. Forrester,
\href{https://doi.org/10.1088/0951-7715/26/6/1799}
{\emph{Joint distribution of the first and second eigenvalues at the soft edge
of unitary ensembles}},
Nonlinearity \textbf{26} (2013), 1799--1822.

\end{thebibliography}
\end{document}